\documentclass[11pt]{article}

\usepackage{a4,amsfonts,exscale,hyperref}

\usepackage{srcltx}

\long\def\forget#1{}

\usepackage{amsmath}
\usepackage{amssymb}
\usepackage{amscd}
\usepackage{textcomp}

\usepackage{amsthm}
\theoremstyle{plain}
\newtheorem{theorem}{Theorem}[section]

\newtheorem{lemma}[theorem]{Lemma}

\newtheorem{corollary}[theorem]{Corollary}
\newtheorem{proposition}[theorem]{Proposition}
\newtheorem{conjecture}[theorem]{Conjecture}
       
\theoremstyle{definition}
\newtheorem{definition}[theorem]{Definition}

\newtheorem{example}[theorem]{Example}
\newtheorem{remark}[theorem]{Remark}

\theoremstyle{remark}

\usepackage{xy}
\xyoption{all}

\newdir^{ (}{{}*!/-3pt/\dir^{(}}    
\newdir^{  }{{}*!/-3pt/\dir^{}}    
\newdir_{ (}{{}*!/-3pt/\dir_{(}}    
\newdir_{  }{{}*!/-3pt/\dir_{}}

\newcounter{zahl}

\def\labelenumi{\theenumi} 
\def\theenumi{(\alph{enumi})}

\def\p@enumii{\theenumi}

\newcommand{\DS}{\displaystyle}
\newcommand{\TS}{\textstyle}
\newcommand{\SC}{\scriptstyle}
\newcommand{\SSC}{\scriptscriptstyle}

\DeclareMathOperator{\Aut}{Aut}

\DeclareMathOperator{\End}{End}

\DeclareMathOperator{\Frob}{Frob}
\DeclareMathOperator{\Gal}{Gal}
\DeclareMathOperator{\GL}{GL}
\DeclareMathOperator{\Grass}{Grass}

\DeclareMathOperator{\Koh}{H}

\DeclareMathOperator{\Hom}{Hom}
\newcommand{\CHom}{{\cal H}om}
\DeclareMathOperator{\Id}{Id}

\DeclareMathOperator{\Isom}{Isom}
\DeclareMathOperator{\Lie}{Lie}

\DeclareMathOperator{\Rep}{\ul{Rep}}

\DeclareMathOperator{\Spm}{Sp}

\DeclareMathOperator{\Spf}{Spf}
\DeclareMathOperator{\Stab}{Stab}

\DeclareMathOperator{\Var}{V}

\newcommand{\ad}{{\rm ad}}
\newcommand{\alg}{{\rm alg}}
\newcommand{\an}{{\rm an}}

\newcommand{\con}{{\rm con}}
\newcommand{\cris}{{\rm cris}}
\newcommand{\der}{{\rm der}}

\DeclareMathOperator{\diag}{diag}
\newcommand{\dR}{{\rm dR}}
\newcommand{\et}{{\rm\acute{e}t}}

\DeclareMathOperator{\id}{\,id}
\DeclareMathOperator{\im}{im}

\renewcommand{\mod}{{\rm\,mod\,}}

\newcommand{\rig}{{\rm rig}}
\DeclareMathOperator{\rk}{rk}

\newcommand{\topol}{{\rm top}}

\DeclareMathOperator{\weight}{wt}

\def\ulD{{\underline{D\hspace{-0.1em}}\,}{}}
\def\ulTD{{\underline{\wt{D}\hspace{-0.1em}}\,}{}}

\renewcommand{\phi}{\varphi}
\renewcommand{\epsilon}{\varepsilon}

\usepackage{amsfonts}
\newcommand{\BOne} {{\mathchoice{\hbox{\rm1\kern-2.7pt l\kern.9pt}}
                              {\hbox{\rm1\kern-2.7pt l\kern.9pt}}
                              {\hbox{\scriptsize\rm1\kern-2.3pt l\kern.4pt}}
                              {\hbox{\scriptsize\rm1\kern-2.4pt l\kern.5pt}}}}

\newcommand{\BA}{{\mathbb{A}}}

\newcommand{\BC}{{\mathbb{C}}}
\newcommand{\BD}{{\mathbb{D}}}

\newcommand{\BF}{{\mathbb{F}}}
\newcommand{\BG}{{\mathbb{G}}}

\newcommand{\BN}{{\mathbb{N}}}

\newcommand{\BP}{{\mathbb{P}}}
\newcommand{\BQ}{{\mathbb{Q}}}
\newcommand{\BR}{{\mathbb{R}}}

\newcommand{\BZ}{{\mathbb{Z}}}

\newcommand{\bA}{{\mathbf{A}}}
\newcommand{\bB}{{\mathbf{B}}}

\newcommand{\bD}{{\mathbf{D}}}
\newcommand{\bE}{{\mathbf{E}}}
\newcommand{\be}{{\mathbf{e}}}

\newcommand{\bM}{{\mathbf{M}}}
\newcommand{\bN}{{\mathbf{N}}}

\newcommand{\bX}{{\mathbf{X}}}

\newcommand{\CalD}{{\cal{D}}}
\newcommand{\CE}{{\cal{E}}}
\newcommand{\CF}{{\cal{F}}}
\newcommand{\CG}{{\cal{G}}}
\newcommand{\CH}{{\cal{H}}}

\newcommand{\CJ}{{\cal{J}}}

\newcommand{\CL}{{\cal{L}}}
\newcommand{\CM}{{\cal{M}}}
\newcommand{\CN}{{\cal{N}}}
\newcommand{\CO}{{\cal{O}}}

\newcommand{\CQ}{{\cal{Q}}}

\newcommand{\CS}{{\cal{S}}}

\newcommand{\CV}{{\cal{V}}}

\newcommand{\rinj}{ \mbox{\mathsurround=0pt \;\raisebox{0.63ex}{\small $\subset$} \hspace{-1.07em} $\longrightarrow$\;}}

\newcommand{\rsurj}{\mbox{\mathsurround=0pt \;$\longrightarrow \hspace{-0.7em} \to$\;}}

\let\setminus\smallsetminus

\newcommand{\es}{\enspace}

\newcommand{\dual}{^{\SSC\!\lor}}

\newcommand{\mal}{^{\SSC\times}}

\newcommand{\ul}[1]{{\underline{#1}}}
\newcommand{\ol}[1]{{\overline{#1}}}

\newcommand{\wt}[1]{{\widetilde{#1}}}

\usepackage{ifthen}

\newcommand{\invlim}[1][]{\ifthenelse{\equal{#1}{}}
{\DS \lim_{\longleftarrow}}
{\DS \lim_{\underset{#1}{\longleftarrow}}}
}

\newcommand{\dirlim}[1][]{\ifthenelse{\equal{#1}{}}
{\DS \lim_{\longrightarrow}}
{\DS \lim_{\underset{#1}{\longrightarrow}}}
}

\newcommand{\dbl}{{\mathchoice{\mbox{\rm [\hspace{-0.15em}[}}
                              {\mbox{\rm [\hspace{-0.15em}[}}
                              {\mbox{\scriptsize\rm [\hspace{-0.15em}[}}
                              {\mbox{\tiny\rm [\hspace{-0.15em}[}}}}
\newcommand{\dbr}{{\mathchoice{\mbox{\rm ]\hspace{-0.15em}]}}
                              {\mbox{\rm ]\hspace{-0.15em}]}}
                              {\mbox{\scriptsize\rm ]\hspace{-0.15em}]}}
                              {\mbox{\tiny\rm ]\hspace{-0.15em}]}}}}

\newcommand{\PLoc}{\BQ_p\mbox{-}\ul{\rm Loc}}

\newcommand{\dotBD}{\vbox{\hbox{\kern2pt\bf.}\vskip-4.5pt\hbox{$\BD$}}}

\def\?{\ 
???\ \immediate\write16{}
\immediate\write16{Warning: There was still a question mark . . . }
\immediate\write16{}}

\DeclareMathOperator{\Nilp}{\CN \!{\it ilp}}
\DeclareMathOperator{\Sets}{\CS \!{\it ets}}
\newcommand{\HeckeTower}{{\breve\CE}}
\newcommand{\hhh}{n}
\newcommand{\nnn}{\lambda}
\newcommand{\nn}{h}
\newcommand{\LL}{L}
\newcommand{\LF}{F}

\def\longto{\longrightarrow}
\def\into{\hookrightarrow}
\def\onto{\rsurj}
\def\isoto{\stackrel{}{\mbox{\hspace{1mm}\raisebox{+1.4mm}{$\SC\sim$}\hspace{-3.5mm}$\longrightarrow$}}}

\newbox\mybox
\def\arrover#1{\mathrel{
       \setbox\mybox=\hbox spread 1.4em{\hfil$\scriptstyle#1$\hfil}
       \vbox{\offinterlineskip\copy\mybox
             \hbox to\wd\mybox{\rightarrowfill}}}}

\begin{document}


\author{Urs Hartl
\footnote{The author acknowledges support of the Deutsche Forschungsgemeinschaft in form of DFG-grant HA3006/2-1 and SFB 878.}
}

\title{On a Conjecture of Rapoport and Zink}

\maketitle

\begin{abstract}
\noindent
In their book, 
Rapoport and Zink constructed rigid analytic period spaces $\CF^{wa}$ for Fontaine's filtered isocrystals, and period morphisms from PEL moduli spaces of $p$-divisible groups to some of these period spaces. They conjectured the existence of an \'etale bijective morphism $\CF^a\to\CF^{wa}$ of rigid analytic spaces and of a universal local system of $\BQ_p$-vector spaces on $\CF^a$. Such a local system would give rise to a tower of \'etale covering spaces $\HeckeTower_{\wt K}$ of $\CF^a$, equipped with a Hecke-action, and an action of the automorphism group $J(\BQ_p)$ of the isocrystal with extra structure.

For Hodge-Tate weights $n-1$ and $n$ we construct in this article an intrinsic Berkovich open subspace $\CF^0$ of $\CF^{wa}$ and the universal local system on $\CF^0$. We show that only in exceptional cases $\CF^0$ equals all of $\CF^{wa}$ and when the Shimura group is $\GL_n$ we determine all these cases. We conjecture that the rigid-analytic space associated with $\CF^0$ is the maximal possible $\CF^a$, and that $\CF^0$ is connected. We give evidence for these conjectures. For those period spaces possessing PEL period morphisms, we show that $\CF^0$ equals the image of the period morphism. Then our local system is the rational Tate module of the universal $p$-divisible group and carries a $J(\BQ_p)$-linearization. We construct the tower $\HeckeTower_{\wt K}$ of \'etale covering spaces, and we show that it is canonically isomorphic in a Hecke and $J(\BQ_p)$-equivariant way to the tower constructed by Rapoport and Zink using the universal $p$-divisible group.

\end{abstract}

\tableofcontents

\bigskip

%
%

\section{Introduction}

The conjecture of Rapoport and Zink, we discuss in this article, is concerned with $p$-adic period spaces, the existence of local systems of $\BQ_p$-vector spaces on them and their resulting \'etale covering spaces. Rapoport and Zink fix a reductive group $G$ over $\BQ_p$, an element $b\in G(K_0)$ for $K_0:=W(\BF_p^{\,\alg})[\frac{1}{p}]$, and a conjugacy class $\{\mu\}$ of cocharacters of $G$. They consider the projective variety $\breve\CF$ over a finite extension $\breve E/K_0$ parametrizing the cocharacters belonging to this conjugacy class. They show that the set of weakly admissible cocharacters is a rigid analytic subspace $(\breve\CF^{wa}_b)^\rig$ of $\breve\CF$, which is called a \emph{$p$-adic period space}; see \cite[Proposition 1.36]{RZ}. On \cite[p.~29]{RZ} they conjecture the existence of an \'etale morphism $(\breve\CF^a_b)^\rig\to(\breve\CF^{wa}_b)^\rig$ of rigid analytic spaces which is bijective on rigid analytic points, and the existence of a tensor functor $\ul\CV$ from $\BQ_p$-rational representations $\rho:G\to\GL(V)$ to local systems of $\BQ_p$-vector spaces on $(\breve\CF^a_b)^\rig$, such that at every point $\mu$ of $(\breve\CF^a_b)^\rig$ the fiber of the local system $\ul\CV(\rho)$ is the crystalline Galois representation, which by the Colmez-Fontaine Theorem~\cite{CF} corresponds to the filtration on the isocrystal $(V\otimes_{\BQ_p}K_0\,,\,\rho(b)\cdot\phi)$ given by $\rho\circ\mu$. Here $\phi$ is the $p$-Frobenius on $K_0$.

The interest in such a tensor functor comes from the fact that it defines a tower $(\HeckeTower_{\wt K})_{\wt K\subset\wt G(\BQ_p)}$ of \'etale covering spaces of $(\breve\CF_b^a)^\rig$, on which an inner form $\wt G(\BQ_p)$ of $G$ acts through Hecke-correspondences, and the group $J(\BQ_p)=\{g\in G(K_0):g^{-1}b\phi(g)=b\}$ acts horizontally; see Remark~\ref{Rem1.8} for details. The first instance of such a tower was constructed by Drinfeld~\cite{Drinfeld76}. There $G$ is the group of units in the central division algebra over $\BQ_p$ with Hasse invariant $-\tfrac{1}{n}$, the group $J(\BQ_p)$ equals $\GL_n(\BQ_p)$, and the period space is the \emph{Drinfeld upper halfspace} $\Omega$ which is the complement of all $\BQ_p$-rational hyperplanes in $\BP^{n-1}$; see \cite[1.44--1.46]{RZ}. It has been shown that the $\ell$-adic cohomology of Drinfeld's tower realizes local Langlands and Jaquet-Langlands correspondences \cite{Carayol90,Harris97,Harris-Taylor}. Inspired by this case, Rapoport \cite[Hope 4.2]{Rapoport95} asked whether such towers also exist over general $p$-adic period spaces. For those period spaces possessing PEL period morphisms, Rapoport and Zink \cite[5.32--5.39]{RZ} constructed a tower and conjectured that it consists of \'etale covering spaces of the image of the period morphism \cite[1.37]{RZ}. We prove this conjecture in Theorem~\ref{ThmPEL} by determining the image of the period morphism, by constructing the tensor functor and the associated tower $(\HeckeTower_{\wt K})_{\wt K\subset\wt G(\BQ_p)}$ of \'etale covering spaces over this image, and by showing that the latter tower coincides with the tower of \cite[5.32--5.39]{RZ}. We show in Example~\ref{Example6.4} and Theorem~\ref{Thm8.2} that the distinction between $(\breve\CF_b^{wa})^\rig$ and $(\breve\CF_b^a)^\rig$ is crucial for this, because the image of the period morphism is smaller than $(\breve\CF_b^{wa})^\rig$ except in a few low-dimensional cases. Therefore one should in general not expect that the tensor functor and the tower exists over all of $(\breve\CF_b^{wa})^\rig$. Let us give more details.

\medskip

De Jong~\cite{dJ} pointed out that to study local systems and \'etale covering spaces, it is best to work in the category of Berkovich's $\breve E$-analytic spaces rather than rigid analytic spaces. For convenience of the reader we review Berkovich's definition and the relation with rigid analytic spaces in Appendix~\ref{AppBerkovichSpaces}. We show that in fact $(\breve\CF^{wa}_b)^\rig$ is the rigid analytic space associated with an open $\breve E$-analytic subspace $\breve\CF^{wa}_b$ of $\breve\CF$ (Proposition~\ref{Prop1.2b}). In Section~\ref{SectPeriodSpaces} we recall the definition of $(\breve\CF^{wa}_b)^\rig$ and the precise formulation of the conjecture. We slightly strengthen the conjecture and reformulate it in terms of an $\breve E$-analytic space $\breve\CF_b^a$ in Conjecture~\ref{ConjRZ}. We also explain its interpretation in terms of the \'etale fundamental group, and we construct the tower $(\HeckeTower_{\wt K})_{\wt K\subset\wt G(\BQ_p)}$ of \'etale covering spaces of $\breve\CF_b^a$ in Section~\ref{SectPeriodSpaces}.

We further contribute to the conjecture in the situation where $G$ has a faithful representation $\rho$ such that all Hodge-Tate weights of $\rho\circ\mu$ are $n-1$ and $n$ for some integer $n$. Namely, we construct an intrinsic open $\breve E$-analytic subspace $\breve\CF_b^0$ of $\breve\CF_b^{wa}$ (Theorem~\ref{Thm2.1}) and the universal local system on $\breve\CF_b^0$ (Theorem~\ref{ThmDOR}). We conjecture that $\breve\CF^0_b$ is connected (Conjecture~\ref{Conj5.3c}) and equals the largest possible $\breve\CF_b^a$ as above (Conjecture~\ref{Conj5.10}). We provide some evidence for these conjectures in Section~\ref{SectConjectures}. The construction of $\breve\CF_b^0$ proceeds as follows. With any analytic point $\mu$ of $\breve\CF$ (these are the ones of which $\breve E$-analytic spaces consist; see Definition~\ref{DefAnalyticPoint}) we associate in Section~\ref{SectMinuscule} a $\phi$-module $\bM_\mu$ over the ring $\wt\bB^\dagger_\rig$ and we let $\breve\CF^0_b$ be the set of those $\mu$ for which $\bM_\mu$ is isoclinic of slope zero. This explicit condition is approachable by combinatorial properties of the Newton slopes of $b$ as we exemplify by our computations in Example~\ref{Example6.4} and Theorem~\ref{Thm8.2}. The ring $\wt\bB^\dagger_\rig$, which is a kind of maximal unramified extension of the Robba ring, is defined in Section~\ref{SectFontaine'sRings}, and the notion of $\phi$-modules and their elementary properties are recalled in Section~\ref{SectPhiModules}. Our construction is inspired by Berger's~\cite{Berger08} construction which associates with any \emph{rigid analytic point} $\mu\in\breve\CF$ (the ones whose residue field $\CH(\mu)$ is \emph{finite} over $\breve E$) a $(\phi,\Gamma)$-module over the Robba ring. Due to the restriction to rigid analytic points, Berger's approach works even without the above assumption on the Hodge-Tate weights. However, the techniques we use in Section~\ref{SectOpenness}, where we prove that $\breve\CF^0_b$ is an open $\breve E$-analytic subspace of $\breve\CF$, require in an essential way to work with $\breve E$-analytic spaces. These also have points whose residue field is not finite over $\breve E$. Moreover, they are topological spaces in the classical sense and we show that the complement $\breve\CF\setminus\breve\CF^0_b$ is the image of a compact set under a continuous map (Theorem~\ref{Thm2.1}). It is worth noting that this crude map is in fact only continuous and not a morphism of $\breve E$-analytic or rigid analytic spaces.

For Hodge-Tate weights $-1$ and $0$ Rapoport and Zink also study a period morphism $\breve\pi:\breve\CG^\an\to\breve\CF_b^{wa}$ starting at the $\breve E$-analytic space $\breve\CG^\an$ associated with a Rapoport-Zink space $\breve\CG$, that is, a PEL moduli space of $p$-divisible groups in a fixed isogeny class. The problem to determine the image of such period morphisms was already mentioned by Grothendieck~\cite{Grothendieck2}. By the Colmez-Fontaine Theorem~\cite{CF} it contains every classical rigid analytic point. When $G=\GL(V)$, de Jong~\cite[\S\S6 and 7]{dJ} observed that the rational Tate module of the universal $p$-divisible group over $\breve\CG^\an$ descends to a local system $\CV$ of $\BQ_p$-vector spaces on the image $\breve\pi(\breve\CG^\an)$ of $\breve\pi$; compare Theorem~\ref{ThmC} where we also show that $\breve\pi(\breve\CG^\an)=\breve\CF^0_b$ and that $\breve\CG^\an$ is identified with the space of $\BZ_p$-lattices inside $\CV$. This uses a recent result of Faltings~\cite{Faltings08}, or alternatively a result of Scholze and Weinstein~\cite{ScholzeWeinstein}. In the general PEL-situation we prove in Theorem~\ref{ThmPEL} that $\breve\pi(\breve\CG^\an)=\breve\CF^0_b$ and that the rational Tate-module of the universal $p$-divisible group over $\breve\CG^\an$ induces a tensor functor $\ul\CV$ from $\BQ_p$-rational representations of $G$ to local systems of $\BQ_p$-vector spaces on $\breve\CF^0_b$ as in the conjecture of Rapoport and Zink. We conjecture that $\breve\CF^0_b$ is the largest open $\breve E$-analytic subset on which such a tensor functor exists (Conjectures~\ref{Conj6.5} and \ref{Conj7.6}) and we give some evidence for this in Section~\ref{SectConjectures}. We further show in Theorems~\ref{ThmC} and \ref{ThmPEL} that the tensor functor $\ul\CV$ carries a canonical $J(\BQ_p)$-linearization which gives rise to a horizontal $J(\BQ_p)$-action on the associated tower $(\HeckeTower_{\wt K})_{\wt K\subset\wt G(\BQ_p)}$ of \'etale covering spaces. And we prove that $(\HeckeTower_{\wt K})_{\wt K\subset\wt G(\BQ_p)}$ is canonically and $J(\BQ_p)\times\wt G(\BQ_p)$-equivariantly isomorphic to the tower constructed by Rapoport and Zink \cite[5.32--5.39]{RZ} using the universal $p$-divisible group on $\breve\CG^\an$.

Note that on the associated rigid analytic spaces, the inclusion $\breve\pi(\breve\CG^\an)=\breve\CF_b^0\subset\breve\CF_b^{wa}$ induces an \'etale morphism $\breve\pi(\breve\CG^\an)^\rig\to (\breve\CF_b^{wa})^\rig$ which is bijective on classical rigid analytic points. It was noticed in \cite{RZ} and \cite{dJ} as a peculiarity that this \'etale morphism might nevertheless fail to be an isomorphism. We give a natural explanation of this phenomenon. Indeed, we show that in general $\breve\pi(\breve\CG^\an)=\breve\CF^0_b$ is strictly contained in $\breve\CF^{wa}_b$ (Example~\ref{Example6.4}) and hence $\breve\pi(\breve\CG^\an)^\rig\to(\breve\CF_b^{wa})^\rig$ in general is not an isomorphism (Proposition~\ref{Prop2.1b}). For $G=\GL_n$ we even compile a complete list of those $b\in\GL_n(K_0)$ for which $\breve\pi(\breve\CG^\an)=\breve\CF_b^{wa}$ in Theorem~\ref{Thm8.2}. The fact that the inclusion $\breve\pi(\breve\CG^\an)\subset\breve\CF_b^{wa}$ may be strict while it induces a bijection on rigid analytic points by the Colmez-Fontaine Theorem~\cite{CF} again shows that one must take into considerations also the points of $\breve\CF_b^{wa}$ whose residue field is not finite over $\breve E$. In other words one should work with Berkovich's $\breve E$-analytic spaces rather than rigid analytic spaces.

To end this introduction let us mention that some of the results presented here were announced in \cite{HartlCR-RZ} and already summarized in \cite[\S11.4]{DOR}. The ideas in this article are inspired by our analogous theory in equal characteristic~\cite{HartlPSp}, or \cite[\S6.2]{HartlDict}, where we were able to prove the analog of the Rapoport-Zink Conjecture and the surjectivity of the period morphism onto $\breve\CF^0_b$ for arbitrary Hodge-Tate weights. In the beginning we had hoped to be able to extend our approach here also to Hodge-Tate weights other than $n-1$ and $n$, but we did not succeed. Not even the construction of the $\phi$-module $\bM_\mu$ over $\wt\bB^\dagger_\rig$, let alone the openness result could be established along this line. The reason lies in the fact that for analytic points $\mu\in\breve\CF$ whose residue field $\CH(\mu)$ is not finite over $\breve E$ the ring $\wt\bB^\dagger_\rig$ is not a $\CH(\mu)$-algebra (like Fontaine's field $\bB_\dR$ is not a $\BC_p$-algebra). So Berger's construction could not be adapted. We also mention that Kedlaya~\cite{Kedlaya10} has announced a strategy to overcome these difficulties and to extend our results to arbitrary Hodge-Tate weights; see also Remarks~\ref{RemarkKedlaya} and \ref{Rem1.4}(a).

\medskip

\emph{Acknowledgments.} The author would like to thank the anonymous referee for his comments and explanations. They helped to improve the exposition. He also would like to thank F.~Andreatta, L.~Berger, V.~Berko\-vich, O.~Brinon, J.-M.~Fontaine, A.~Genestier, U.~G\"ortz, M.~Kisin, V.~Lafforgue, J.~Pottharst, M.~Rapo\-port, E.~Viehmann, and T.~Zink for many helpful discussions and for their interest in this work.

%
%

\section{The Conjecture of Rapoport and Zink} \label{SectPeriodSpaces}
\setcounter{equation}{0}

In this section we recall the construction of Rapoport's and Zink's $p$-adic period spaces. We strengthen their conjecture mentioned above and interpret it in terms of the \'etale fundamental group. We also explain its relation with towers of \'etale covering spaces. Let us first recall the definition of filtered isocrystals. We denote by $\BF_p^{\,\alg}$ an algebraic closure of the finite field with $p$ elements and by $K_0:=W(\BF_p^{\,\alg})[\frac{1}{p}]$ the fraction field of the ring of Witt vectors over $\BF_p^{\,\alg}$. Let $\phi=W(\Frob_p)$ be the Frobenius lift on $K_0$. 

\begin{definition}\label{DefIsoc}
An \emph{$F$-isocrystal over $\BF_p^{\,\alg}$} is a finite dimensional $K_0$-vector space $D$ equipped with a $\phi$-linear automorphism $\phi_D$. If $\LL$ is a field extension of $K_0$ and $Fil^\bullet D_\LL$ is an exhaustive separated decreasing filtration of $D_\LL:=D\otimes_{K_0}\LL$ by $\LL$-subspaces we say that $\ulD=(D,\phi_D,Fil^\bullet D_\LL)$ is a \emph{filtered isocrystal over $\LL$}. The integers $i$ for which $Fil^{-i}D_\LL\ne Fil^{-i+1}D_\LL$ are called the \emph{Hodge-Tate weights} of $\ulD$. We let $t_N(\ulD)$ be the $p$-adic valuation of $\det\phi_D$ (with respect to any basis of $D$) and we let 
\[
t_H(\ulD)\es=\es\sum_{i\in\BZ}i\cdot\dim_\LL gr^i_{Fil^\bullet}(D_\LL)\,.
\]
The filtered isocrystal $\ulD$ is called \emph{weakly admissible}\footnote{This used to be the terminology until Colmez--Fontaine~\cite{CF} showed that for $\LL/K_0$ finite \emph{weakly admissible implies admissible}. Since we consider also infinite extensions $\LL/K_0$ for which the Colmez--Fontaine Theorem fails we stick to the old terminology.} if 
\begin{equation}\label{Eq2.1}
t_H(\ulD)\es=\es t_N(\ulD)\qquad\text{and}\qquad t_H(\ulD')\es\le\es t_N(\ulD')\es
\end{equation}
for any subobject $\ulD'=\bigl(D',\phi_D|_{D'},Fil^\bullet D'_\LL\bigr)$ of
$\ulD$, where $D'$ is any $\phi_D$-stable $K_0$-subspace of $D$ which is equipped with
the induced filtration $Fil^i D'_\LL=D'_\LL\cap Fil^i D_\LL$ on $D'_\LL:=D'\otimes_{K_0}\LL$.
\end{definition}

\medskip

To construct period spaces let $G$ be a reductive linear algebraic group over $\BQ_p$. Fix a conjugacy class $\{\mu\}$ of cocharacters
\[
\mu:\BG_m\to G
\]
defined over subfields of $\BC_p$ (the $p$-adic completion of an algebraic closure of $K_0$). Let $E$ be the field of definition of the conjugacy class. It is a finite extension of $\BQ_p$. Two cocharacters in this conjugacy class are called \emph{equivalent} if they induce the same weight filtration on the category $\Rep_{\BQ_p}\!\!\!G$ of finite dimensional $\BQ_p$-rational representations of $G$. There is a projective variety $\CF$ over $E$ whose $\BC_p$-valued points are in bijection with the equivalence classes of cocharacters (from the fixed conjugacy class $\{\mu\}$). Namely let $V$ in $\Rep_{\BQ_p}\!\!\!G$ be any faithful representation of $G$. With a cocharacter $\mu$ defined over a field $\LL$ one associates the filtration $Fil^i_\mu V_\LL:=\bigoplus_{j\ge i}V_{\LL,j}$ of $V_\LL:=V\otimes_{\BQ_p}\LL$ given by the weight spaces
\[
V_{\LL,j}\es:=\es\bigl\{\,v\in V_\LL:\es\mu(z)\cdot v\;=\;z^jv\text{ for all }z\in\BG_m(\LL)\,\bigr\}\,.
\]
This defines a closed embedding of $\CF$ into a partial flag variety of $V$
\begin{equation}\label{Eq2.0}
\CF\es\rinj\es\CF lag(V)\otimes_{\BQ_p}E\,,
\end{equation}
where the points of $\CF lag(V)$ with value in a $\BQ_p$-algebra $R$ are the filtrations $F^i$ of $V\otimes_{\BQ_p} R$ by $R$-submodules which are direct summands such that $\rk_R gr^i_{F^\bullet}$ is the multiplicity of the weight $i$ of the conjugacy class $\{\mu\}$ on $V$.

\bigskip

Now let $b\in G(K_0)$. For any $V$ in $\Rep_{\BQ_p}\!\!\!G$ one obtains an $F$-isocrystal $(V\otimes_{\BQ_p}K_0\,,\, b\!\cdot\!\phi)$. Its automorphism group contains the group from \cite[Proposition 1.12]{RZ}
\begin{equation}\label{EqJ}
J(\BQ_p)\;=\;\bigl\{\,g\in G(K_0):g^{-1}b\phi(g)=b\,\bigr\}\,.
\end{equation}

A pair $(b,\mu)$ with an element $b\in G(K_0)$ and a cocharacter $\mu:\BG_m\to G$ defined over $\LL\supset K_0$ is called \emph{weakly admissible} if for some faithful representation $\rho:G\into\GL(V)$ in $\Rep_{\BQ_p}\!\!\!G$ the filtered isocrystal
\begin{equation}\label{EqIsocrystal}
\ulD_{b,\mu}(V)\es:=\es(V\otimes_{\BQ_p}K_0\,,\,\rho(b)\!\cdot\!\phi\,,\,Fil^\bullet_\mu V_\LL)
\end{equation}
is weakly admissible. In fact this then holds for any $V$ in $\Rep_{\BQ_p}\!\!\!G$; see \cite[1.18]{RZ}. 
If $(b,\mu)$ is weakly admissible and $\LL/K_0$ is finite then the filtered isocrystal $\ulD_{b,\mu}(V)$ is \emph{admissible}, that is it arises from a crystalline Galois-representation $\Gal(\LL^\alg/\LL)\to\GL(U)$ via Fontaine's covariant functor
\[
\ulD_{b,\mu}(V)\es\cong\es\bD_{\cris}(U)\es:=\es (U\otimes_{\BQ_p}\bB_{\cris})^{\Gal(\LL^\alg/\LL)}
\]
which by the Colmez-Fontaine theorem~\cite{CF} is an equivalence of categories from crystalline representations of $\Gal(\LL^\alg/\LL)$ to weakly admissible filtered isocrystals over $\LL$. We write $V_\cris(\ulD_{b,\mu}(V)):=U$ for Fontaine's covariant inverse functor. The assignment
\begin{equation}\label{Eq2.1a}
\Rep_{\BQ_p}\!\!\!G\es\longto\es\Rep_{\BQ_p}\!\!\bigl(\Gal(\LL^\alg/\LL)\bigr)\,,\quad V\es\mapsto\es V_\cris(\ulD_{b,\mu}(V))
\end{equation}
defines a tensor functor from $\Rep_{\BQ_p}\!\!\!G$ to the category of continuous $\Gal(\LL^\alg/\LL)$-represen\-tations in finite dimensional $\BQ_p$-vector spaces.

Let $\breve E=EK_0$ be the completion of the maximal unramified extension of $E$. In what follows we consider cocharacters $\mu$ defined over complete extensions $\LL$ of $\breve E$. Let $\breve\CF^\rig$ be the rigid analytic space over $\breve E$ associated with the variety $\breve\CF=\CF\otimes_E\breve E$. Rapoport and Zink define the \emph{$p$-adic period space} associated with $(G,b,\{\mu\})$ as
\[
(\breve\CF_b^{wa})^\rig\es:=\es\bigl\{\,\mu\in\breve\CF^\rig:\es(b,\mu)\text{ is weakly admissible}\,\bigr\}\,.
\]
It is acted on by the group $J(\BQ_p)$ from \eqref{EqJ} via $\mu\mapsto\gamma\mu\gamma^{-1}$ for $\gamma\in J(\BQ_p)$.
Rapoport and Zink show that $(\breve\CF^{wa}_b)^\rig$ is an admissible open rigid analytic subspace of $\breve\CF^\rig$; see \cite[Proposition 1.36]{RZ}. For the proof they reduce to the case where $b$ is \emph{decent with some positive integer $s$}, that is
\begin{equation}\label{EqDecent}
b\cdot\phi(b)\cdot\ldots\cdot\phi^{s-1}(b)\es=\es s\nu(p)\,.
\end{equation}
where $\nu\in\Hom_{K_0}(\BG_m,G)\otimes_\BZ \BQ$ is the slope quasi-cocharacter \cite[1.7]{RZ} and $s$ satisfies $s\nu\in\Hom_{K_0}(\BG_m,G)$. If $b$ is decent with integer $s$ then $(\breve\CF^{wa}_b)^\rig$ has a natural structure of rigid analytic subspace of $(\CF\otimes_E E_s)^\rig$ over the field $E_s=E\cdot W(\BF_{p^s})[\frac{1}{p}]$ from which it arises by base change to $\breve E$. In \cite[p.~29]{RZ} Rapoport and Zink make a conjecture on the existence of a universal local system of $\BQ_p$-vector spaces on $(\breve\CF^{wa}_b)^\rig$ which we refine in \ref{ConjRZ} below.

Stimulated by \cite{RZ}, the notion of local systems of $\BQ_p$-vector spaces on rigid analytic spaces was studied by de Jong~\cite[\S 4]{dJ} who pointed out that this is best done working with Berkovich's $\breve E$-analytic spaces rather than rigid analytic spaces. So let $\breve\CF^\an$ be the $\breve E$-analytic space associated with the variety $\breve\CF$. See Appendix~\ref{AppBerkovichSpaces} for the notion of $\breve E$-analytic space and the relation with rigid analytic spaces. In fact the argument of \cite[Proposition 1.36]{RZ} shows that there exists an open $\breve E$-analytic subspace $\breve\CF^{wa}_b$ of $\breve\CF^\an$ whose associated rigid analytic space is the period space $(\breve\CF^{wa}_b)^\rig$; we explain the proof of this fact in Proposition~\ref{Prop1.2b}.

Local systems of $\BQ_p$-vector spaces on an $\breve E$-analytic space $Y$ are related to the category $\Rep_{\BQ_p}\!\!\bigl(\pi_1^\et(Y,\bar y)\bigr)$ of finite dimensional continuous $\BQ_p$-linear representations of the \'etale fundamental group $\pi_1^\et(Y,\bar y)$ of $Y$ as follows. This group was defined by de Jong~\cite{dJ}. 

\begin{proposition} \label{Prop1.2c}
(\cite[Theorem 4.2]{dJ}.) For any geometric base point $\bar y$ of $Y$ there is a natural $\BQ_p$-linear tensor functor
\[
\omega_{\bar y}:\es\PLoc_Y\es\to\es \Rep_{\BQ_p}\!\!\bigl(\pi_1^\et(Y,\bar y)\bigr)
\]
which assigns to a local system $\CV$ the $\pi_1^\et(Y,\bar y)$-representation $\CV_{\bar y}$. It is an equivalence if $Y$ is connected.
\end{proposition}

In particular, for any point $y\in Y(\LL)$ and geometric base point $\bar y$ lying above $y$, any local system $\CV$ of $\BQ_p$-vector spaces on $Y$ gives rise to a Galois representation
\begin{equation}\label{EqGalRep}
\Gal(\LL^\alg/\LL)\;=\;\pi_1^\et(y,\bar y)\;\to\;\pi_1^\et(Y,\bar y)\;\to\;\GL(\CV_{\bar y})\,.
\end{equation}
We denote this Galois representation by $\CV_y$.
Now we propose the following refinement of the conjecture of Rapoport and Zink \cite[p.~29]{RZ}; see also \cite[Conjecture 11.4.4]{DOR}.

\begin{conjecture}\label{ConjRZ} There exists a unique largest arcwise connected dense open $\breve E$-analytic subspace $\breve\CF^a_b\subset\breve\CF^{wa}_b$ invariant under $J(\BQ_p)$ with $\breve\CF^a_b(\LL)=\breve\CF^{wa}_b(\LL)$ for all finite extensions $\LL/\breve E$, and a tensor functor $\ul\CV$ from $\Rep_{\BQ_p}\!\!\!G$ to the category $\PLoc_{\breve\CF^a_b}$ of local systems of $\BQ_p$-vector spaces on $\breve\CF^a_b$ with the following property:
\begin{equation} \label{EqRZ}
\hspace{0.07\textwidth}\parbox{0.67\textwidth}{For any point $\mu\in\breve\CF^a_b(\LL)$ with $\LL/\breve E$ finite, the tensor functor
$\Rep_{\BQ_p}\!\!\!G \to \Rep_{\BQ_p}\!\!\bigl(\Gal(\LL^\alg/\LL)\bigr)\,,\;V\mapsto V_\cris(\ulD_{b,\mu}(V))$ from (\ref{Eq2.1a}) is isomorphic to the tensor functor $V\mapsto\ul\CV(V)_{\mu}$ from \eqref{EqGalRep} which associates with a representation $V\in\Rep_{\BQ_p}\!\!\!G$ the fiber at $\mu$ of the corresponding local system $\ul\CV(V)$.}\quad\qquad\qquad\qquad\quad \eqref{EqRZ}\nonumber
\end{equation}
\addtocounter{equation}{1}
\end{conjecture}

\begin{remark}\label{Rem1.4}
(a) There seem to be reasons to expect that the conjecture is false if $\mu$ is not minuscule; see 
\cite{ScholzeHodge,ScholzeWeinstein}.

\medskip\noindent
(b) Note that any $\breve E$-analytic space is locally arcwise connected by \cite[Theorem 3.2.1]{Berkovich1}. Moreover, every open $\breve E$-analytic subspace of $\breve\CF^\an$ is paracompact by Lemma~\ref{LemmaParacompact}.

\medskip\noindent
(c) Forgetting the action of $\Gal(\LL^\alg/\LL)$, Condition (\ref{EqRZ}) implies that
\begin{equation}\label{EqRZb}
\hspace{0.07\textwidth}\parbox{0.8\textwidth}{for any point $\mu\in\breve\CF^a_b(\LL)$ with $\LL/\breve E$ finite, the fiber functors $\Rep_{\BQ_p}\!\!\!G \to (\BQ_p\text{-vector spaces})\,,\;V\mapsto V_\cris(\ulD_{b,\mu}(V))$ and $V\mapsto\ul\CV(V)_{\mu}$ are isomorphic.}\qquad\quad (\ref{EqRZb})\nonumber
\addtocounter{equation}{1}
\end{equation}
But note that \eqref{EqRZb} for $G=\GL_n$ simply says that $\dim_{\BQ_p}\ul\CV(V)=\dim_{\BQ_p}V$ and also in general \eqref{EqRZb} is weaker than Condition~\eqref{EqRZ}.

\medskip\noindent
(d) In a preprint version of this article (Remark~7.2 in \href{http://arxiv.org/abs/math/0605254v1}{arXiv:math.NT/0605254v1}) we stated (and \cite[Remark~11.4.5]{DOR} quoted this) that the existence of a unique largest dense open $\breve E$-analytic subspace $\breve\CF^a_b\subset\breve\CF^{wa}_b$ satisfying \eqref{EqRZ} is automatic. We actually do not know whether this is correct. The problem is that on the same open subspace $U\subset\breve\CF^{wa}_b$ there might be two different tensor functors satisfying \eqref{EqRZ}, because the groups $\Gal(\CH(x)^\alg/\CH(x))$ are not necessarily dense in $\pi_1^\et(U,\bar y)$ when $x$ ranges over all points of $U$, and even less when in addition $\CH(x)/\breve E$ is finite. Indeed, these groups lie in the kernel of $\pi_1^\et(U,\bar y)\to\pi_1^\topol(U,\bar y)$ which however has dense image; see \cite[Theorem~2.10(iv)]{dJ}. Therefore the tensor functors on two different open subspaces do not automatically glue to a tensor functor on the union of the open subspaces.
\end{remark}

To derive the original conjecture of Rapoport and Zink from Conjecture~\ref{ConjRZ} we in addition note the following

\begin{proposition}\label{Prop2.1b}
If Conjecture~\ref{ConjRZ} holds, then the open immersion $\breve\CF^a_b\subset\breve\CF^{wa}_b$ induces an \'etale morphism of rigid analytic spaces $(\breve\CF^a_b)^\rig\to(\breve\CF^{wa}_b)^\rig$ which is bijective on rigid analytic points. It is an isomorphism if and only if $\breve\CF^a_b=\breve\CF^{wa}_b$.
\end{proposition}

\begin{proof}
The functor $(\es)^\rig$ takes \'etale morphisms to \'etale morphisms. The rest is a consequence of Theorem~\ref{ThmFormalRigBerkovich} since $\breve\CF^a_b$ and $\breve\CF^{wa}_b$ are paracompact by Lemma~\ref{LemmaParacompact}.
\end{proof}

The nature of such a morphism $(\breve\CF^a_b)^\rig\to(\breve\CF^{wa}_b)^\rig$ may seem somewhat obscure since it is a non-quasi-compact refinement of the Grothendieck topology on the same underlying set of points $\breve\CF^a_b(\breve E^\alg)=\breve\CF^{wa}_b(\breve E^\alg)$. It is less mysterious on the level of $\breve E$-analytic spaces because there one can see the points in the complement $\breve\CF^{wa}_b\setminus\breve\CF^a_b$; compare Example~\ref{Example6.4} below.

After recalling some of Fontaine's rings and some facts on $\phi$-modules in the next two sections we will construct in Sections~\ref{SectMinuscule} and \ref{SectOpenness} for Hodge-Tate weights $n-1$ and $n$ an open $\breve E$-analytic  subspace $\breve\CF_b^0$ of $\breve\CF^{wa}_b$ which is a candidate for the space searched for in Conjecture~\ref{ConjRZ}.

Let us end this section by explaining the significance of Conjecture~\ref{ConjRZ} in terms of the \'etale fundamental group and in terms of \'etale covering spaces. Let $\mu\in\breve\CF_b^a(\LL)$ be a point with values in a finite extension $\LL/\breve E$ and let $\bar\mu$ be a geometric base point of $\breve\CF_b^a$ lying above $\mu$. Then Rapoport and Zink~\cite[1.19]{RZ} consider a cohomology class
\[
cl(b,\mu)\;=\;\kappa(b)-\mu^\#\;\in\;\Koh^1(\BQ_p,G)
\]
where $\kappa(b)\in \pi_1(G)_{\Gal(\ol\BQ_p/\BQ_p)}$ is the Kottwitz point~\cite[\S7]{Kottwitz97} and $\mu^\#$ is the image of $\mu$ under the natural map $X_\ast(T)\to\pi_1(G)\to\pi_1(G)_{\Gal(\ol\BQ_p/\BQ_p)}$. By the weak admissibility of the pair $(b,\mu)$, the difference $\kappa(b)-\mu^\#$ lies in $(\pi_1(G)_{\Gal(\ol\BQ_p/\BQ_p)})_{tors}$ which Kottwitz~\cite[3.2]{Kottwitz97} identifies with $\Koh^1(\BQ_p,G)$. 

Let $\omega_0$ be the forgetful fiber functor $\Rep_{\BQ_p}\!\!\!G\to(\BQ_p\text{-vector spaces})$. The cohomology class $cl(b,\mu)$ defines a $G$-torsor over $\BQ_p$ which by \cite[Theorem 3.2]{DM} is of the form $\Isom^\otimes(\omega_0,\wt\omega)$ for a fiber functor $\wt\omega:\Rep_{\BQ_p}\!\!\!G\to(\BQ_p\text{-vector spaces})$. Consider the image of $cl(b,\mu)$ in $\Koh^1(\BQ_p,G^\ad)$ and let $\wt G$ be the associated inner form of $G$. This yields a non-canonical isomorphism of $\BQ_p$-group schemes $\wt G\cong\Aut^\otimes(\wt\omega)$ by \cite[Proposition 2.8]{DM} which we fix in the sequel.

\begin{corollary}\label{Cor1.7}
Let $Y\subset\breve\CF_b^{wa}$ be an open connected $\breve E$-analytic subspace, and let $\bar\mu$ be a geometric base point of $Y$. Then the set of isomorphism classes of tensor functors $\ul\CV:\Rep_{\BQ_p}\!\!\!G \to\PLoc_{Y}$ satisfying (\ref{EqRZb}) is in bijection with the set of isomorphism classes of continuous group homomorphisms
\[
\pi_1^\et(Y,\bar \mu)\;\to\;\wt G(\BQ_p)
\]
where two homomorphisms are isomorphic if they differ by composition with an inner automorphism of the $\BQ_p$-group scheme $\wt G$.
\end{corollary}

\begin{proof}
Let $\ul\CV:\Rep_{\BQ_p}\!\!\!G \to\PLoc_{Y}$ be a tensor functor and consider the fiber functor from $\Rep_{\BQ_p}\!\!\!G$ to $(\BQ_p\text{-vector spaces})$ obtained as the composition $\omega_{\bar\mu}\circ\ul\CV$ with the forgetful fiber functor $\omega_{\bar\mu}$ on $\Rep_{\BQ_p}\!\!\bigl(\pi_1^\et(Y,\bar \mu)\bigr)$ from Proposition~\ref{Prop1.2c}. By \cite[Theorem 2.9]{dJ} the fiber functors $\omega_{\bar\mu}$ and $\omega_{\bar\mu'}$ at any two geometric base points $\bar\mu$ and $\bar\mu'$ are isomorphic. So without changing the isomorphism class of $\omega_{\bar\mu}\circ\ul\CV$ we may assume that $\bar\mu$ maps to a point $\mu\in Y(\LL)$ with values in a finite extension $\LL/\breve E$.
Then by  a theorem of Wintenberger~\cite[Corollary to Proposition 4.5.3]{Wintenberger97} (see also \cite[Proposition on page 4]{CF}) Condition~(\ref{EqRZb}) at $\mu$ is equivalent to the assertion that the cohomology class of the $G$-torsor $\Isom^\otimes(\omega_0,\omega_{\bar\mu}\circ\ul\CV)$ equals $cl(b,\mu)\in\Koh^1(\BQ_p,G)$. This yields a non-canonical isomorphism of fiber functors $\beta_0:\wt\omega\isoto\omega_{\bar\mu}\circ\ul\CV$ by \cite[Theorem 3.2]{DM}, and a non-canonical isomorphism of $\BQ_p$-group schemes $\wt G\cong\Aut^\otimes(\wt\omega)\cong\Aut^\otimes(\omega_{\bar\mu}\circ\ul\CV)$ which is uniquely determined only up to an inner automorphism of $\wt G$. Any element $\gamma\in\pi_1^\et(Y,\bar \mu)$ yields a tensor automorphism of $\omega_{\bar\mu}\circ\ul\CV$. This defines a group homomorphism $f:\pi_1^\et(Y,\bar \mu)\to\wt G(\BQ_p)$ which is unique up to composition with an inner automorphism of the group scheme $\wt G$. Since for all $\rho\in\Rep_{\BQ_p}\!\!\!G$ the induced homomorphism \mbox{$\pi_1^\et(Y,\bar \mu)\to\wt G(\BQ_p)\to\GL(\omega_{\bar\mu}\circ\ul\CV(\rho))(\BQ_p)$} is continuous, also $f$ is continuous.

Conversely let $f:\pi_1^\et(Y,\bar \mu)\to\wt G(\BQ_p)$ be a continuous group homomorphism. Then we define a tensor functor $\Rep_{\BQ_p}\!\!\!G \to\Rep_{\BQ_p}\!\!\bigl(\pi_1^\et(Y,\bar \mu)\bigr)$ by sending a representation $\rho$ of $G$ to the representation
\[
\rho':\;\pi_1^\et(Y,\bar\mu)\longto\GL\bigl(\wt\omega(\rho)\bigr)(\BQ_p)\;,\quad \gamma\mapsto\wt\omega(\rho)(f(\gamma))\,.
\]
Here $\wt\omega(\rho)(f(\gamma))$ is the automorphism by which $f(\gamma)\in\wt G(\BQ_p)\cong\Aut^\otimes(\wt\omega)(\BQ_p)$ acts on the $\BQ_p$-vector space $\wt\omega(\rho)$. Note that $\rho'$ is continuous because $\wt G(\BQ_p)\to\GL\bigl(\wt\omega(\rho)\bigr)(\BQ_p)$ is continuous. Let $\ul\CV(\rho):=\omega_{\bar\mu}^{-1}(\rho')$ be the local system of $\BQ_p$-vector spaces on $Y$ induced from $\rho'$ via Proposition~\ref{Prop1.2c}. This defines a tensor functor $\ul\CV:\Rep_{\BQ_p}\!\!\!G \to\PLoc_{Y}$. The composition $\omega_{\bar\mu}\circ\ul\CV$ with the forgetful fiber functor $\omega_{\bar\mu}$ equals $\wt\omega$. Since the cohomology class of $\Isom^\otimes(\omega_0,\wt\omega)$ equals $cl(b,\mu)$, Wintenberger's result mentioned above shows that the tensor functor $\ul\CV$ satisfies condition (\ref{EqRZb}) at $\mu$. Moreover, Wintenberger's theorem implies that the isomorphism class of the fiber functor \eqref{Eq2.1a} is constant on $Y$. Also for any two geometric base points $\bar\mu$ and $\bar\mu'$ of $Y$ the fiber functors $\omega_{\bar\mu}$ and $\omega_{\bar\mu'}$ are isomorphic by \cite[Theorem 2.9]{dJ} since $Y$ is connected. Therefore the isomorphism class of the fiber functor $\omega_{\bar\mu}\circ\ul\CV$ is likewise constant on $Y$ and condition \eqref{EqRZb} holds at any point $\mu\in\breve\CF^a_b(\LL)$ with $\LL/\breve E$ finite. Clearly the assignments $\ul\CV\mapsto f$ and $f\mapsto\ul\CV$ are inverse to each other.
\end{proof}

\begin{remark}\label{Rem1.8}
(a) In the situation of Corollary~\ref{Cor1.7} any tensor functor $\ul\CV:\Rep_{\BQ_p}\!\!\!G \to\PLoc_{Y}$ induces a tower of \'etale covering spaces of $Y$ with Hecke action. More precisely, let $\wt K\subset\wt G(\BQ_p)$ be a compact open subgroup and let $\HeckeTower_{\wt K}$ be the space representing \emph{$\wt K$-level structures} on $\ul\CV$, that is residue classes modulo $\wt K$ of tensor isomorphisms 
\[
(\beta:\wt\omega\isoto\omega_{\bar\mu}\circ\ul\CV)\;\mod \wt K\es\in\es\Isom^\otimes(\wt\omega,\omega_{\bar\mu}\circ\ul\CV)/\wt K
\]
such that the class $\beta\wt K$ is invariant under the \'etale fundamental group of $\HeckeTower_{\wt K}$. Then $\HeckeTower_{\wt K}$ is the \'etale covering space of $Y$ corresponding to the discrete $\pi_1^\et(Y,\bar\mu)$-set $\Isom^\otimes(\wt\omega,\omega_{\bar\mu}\circ\ul\CV)/\wt K$. Any choice of a fixed tensor isomorphism $\beta_0\in\Isom^\otimes(\wt\omega,\omega_{\bar\mu}\circ\ul\CV)$ associates with $\ul\CV$ a representation $\pi_1^\et(Y,\bar\mu)\to\wt G(\BQ_p)$ as in the proof of Corollary~\ref{Cor1.7}, and induces an identification of the $\pi_1^\et(Y,\bar\mu)$-sets $\Isom^\otimes(\wt\omega,\omega_{\bar\mu}\circ\ul\CV)/\wt K$ and $\wt G(\BQ_p)/\wt K$.

On the tower $(\HeckeTower_{\wt K})_{\wt K\subset \wt G(\BQ_p)}$ the group $\wt G(\BQ_p)$ acts via Hecke correspondences: Let $g\in\wt G(\BQ_p)$ and let $\wt K,\wt K'\subset\wt G(\BQ_p)$ be compact open subgroups. Then the Hecke correspondence $\pi(g)_{\wt K,\wt K'}$ is given by the diagram
\begin{equation}\label{EqHeckeAction}
\xymatrix @C=0pc {
& \HeckeTower_{\wt K'\cap g^{-1}\wt Kg}\ar[dl] \ar[dr] & & & &\beta\mod (\wt K'\cap g^{-1}\wt Kg)\ar@{|->}[dl] \ar@{|->}[dr]\\
\HeckeTower_{\wt K'}& & \HeckeTower_{\wt K}\ar@{-->}[ll] &\quad & \beta\mod\wt K' & & \beta\,g^{-1}\mod \wt K
}
\end{equation}
The whole construction does not depend on the choice of the base point $\bar\mu$ by \cite[Theorem 2.9]{dJ} and hence also applies if $Y$ is not connected.

\medskip\noindent
(b) Assume moreover, that a group $\Gamma$ acts on $Y$ and let $\ul\CV:\Rep_{\BQ_p}\!\!\!G \to\PLoc_{Y}$ be a tensor functor which carries a \emph{$\Gamma$-linearization}, that is, for every $\gamma\in\Gamma$ an isomorphism $\phi_\gamma:\gamma^*\ul\CV\isoto\ul\CV$ of tensor functors (where $\gamma^*\ul\CV$ is the pullback of $\ul\CV$ under the morphism $\gamma:Y\to Y$), satisfying the cocycle condition $\phi_\gamma\circ\gamma^*\phi_\delta=\phi_{\delta\gamma}$ for all $\gamma,\delta\in\Gamma$.

Then the tower of \'etale covering spaces $\HeckeTower_{\wt K}$ inherits an action of $\Gamma$ over $Y$ as follows. Recall that $\HeckeTower_{\wt K}$ corresponds to the discrete $\pi_1^\et(Y,\bar\mu)$-set $\Isom^\otimes(\wt\omega,\omega_{\bar\mu}\circ\ul\CV)/\wt K$. By functoriality of the fundamental group, $\gamma^*\HeckeTower_{\wt K}:=\HeckeTower_{\wt K}\times_{Y,\gamma}Y$ corresponds to the discrete $\pi_1^\et(Y,\gamma^{-1}\bar\mu)$-set $\Isom^\otimes(\wt\omega,\omega_{\gamma^{-1}\bar\mu}\circ\gamma^*\ul\CV)/\wt K$ which as a set equals $\Isom^\otimes(\wt\omega,\omega_{\bar\mu}\circ\ul\CV)/\wt K$. By \cite[Theorem~2.9]{dJ} there exists an isomorphism of fiber functors $\psi_\gamma:\omega_{\gamma^{-1}\bar\mu}\isoto\omega_{\bar\mu}$ providing a change of base point isomorphism
\[
\Isom^\otimes(\wt\omega,\omega_{\gamma^{-1}\bar\mu}\circ\gamma^*\ul\CV)/\wt K \isoto \Isom^\otimes(\wt\omega,\omega_{\bar\mu}\circ\gamma^*\ul\CV)/\wt K\,,\quad \beta'\wt K \longmapsto (\psi_\gamma\circ\gamma^*\ul\CV)\circ\beta'\wt K\,.
\]
This shows that the covering space $\gamma^*\HeckeTower_{\wt K}$ of $Y$ corresponds to the discrete $\pi_1^\et(Y,\bar\mu)$-set $\Isom^\otimes(\wt\omega,\omega_{\bar\mu}\circ\gamma^*\ul\CV)/\wt K$. Then the isomorphism
\begin{equation}\label{EqGammaAction}
\Isom^\otimes(\wt\omega,\omega_{\bar\mu}\circ\ul\CV)/\wt K \isoto \Isom^\otimes(\wt\omega,\omega_{\bar\mu}\circ\gamma^*\ul\CV)/\wt K\,,\quad \beta\wt K \longmapsto \omega_{\bar\mu}(\phi_\gamma^{-1})\circ\beta\wt K 
\end{equation}
of  $\pi_1^\et(Y,\bar\mu)$-sets corresponds to an isomorphism $\HeckeTower_{\wt K}\isoto\gamma^\ast\HeckeTower_{\wt K}$ of \'etale covering spaces of $Y$. We compose the latter with the projection onto $\HeckeTower_{\wt K}$ to obtain an isomorphism $\gamma:=\gamma_{\HeckeTower_{\wt K}}$ of $\HeckeTower_{\wt K}$ which makes the following diagram commutative
\[
\xymatrix @R=1pc @C+1pc {
\HeckeTower_{\wt K} \ar[r]^{\DS \gamma_{\HeckeTower_{\wt K}}} \ar[d] & \HeckeTower_{\wt K} \ar[d] \\
Y \ar[r]^{\DS\gamma} & \;\,Y \,.
}
\]
These isomorphisms for varying $\wt K$ are compatible with the projection morphisms of the tower. Furthermore, from the explicit descriptions in \eqref{EqHeckeAction} and \eqref{EqGammaAction} one sees that the action of $\Gamma$ on the tower $(\HeckeTower_{\wt K})_{\wt K\subset \wt G(\BQ_p)}$ commutes with the action of $\wt G(\BQ_p)$ through Hecke correspondences.
\end{remark}

%
%

\section{Fontaine's Rings} \label{SectFontaine'sRings}
\setcounter{equation}{0}

We recall from Colmez~\cite{Colmez} some of the rings used in $p$-adic Hodge theory \cite{Fontaine,Fontaine2,Fontaine94,FW1,FW2}.
Let $\CO_\LL$ be a complete valuation ring of rank one which is an extension of $\BZ_p$ and let $\LL$ be its fraction field. Let $v_p$ be the valuation on $\CO_\LL$ which we assume to be normalized so that $v_p(p)=1$.
In $p$-adic Hodge theory it is usually assumed that $\CO_\LL$ is discretely valued with perfect residue field. However, we make neither of these assumptions here, because in Section~\ref{SectMinuscule} we need this more general situation. Let $C$ be the completion of an algebraic closure $\LL^\alg$ of $\LL$ and let $\CO_C$ be the valuation ring of $C$. For $n\in\BN$ let $\epsilon^{(n)}$ be a primitive $p^n$-th root of unity chosen in such a way that $(\epsilon^{(n+1)})^p=\epsilon^{(n)}$ for all $n$. We  define

\begin{tabbing}
$\DS\wt\bE^+\es:=\es\wt \bE^+(C)\es:= \es\bigl\{\,x=(x^{(n)})_{n\in\BN_0}:\es x^{(n)}\in\CO_{C}\,,\,(x^{(n+1)})^p=x^{(n)}\es\text{for all }n\,\bigr\}$\,.
\end{tabbing}

\noindent
Fix the elements $\epsilon:=(1,\epsilon^{(1)},\epsilon^{(2)},\ldots)$ and $\ol\pi:=\epsilon-1$ of $\wt\bE^+$. 
With the 
multiplication $xy:=\bigl(x^{(n)}y^{(n)}\bigr)_{n\in\BN_0}$, the 
addition $x+y:=\bigl(\lim_{m\to\infty}(x^{(m+n)}+y^{(m+n)})^{p^m}\bigr)_{n\in\BN_0}$, and 
the 
valuation $v_\bE(x):=v_p(x^{(0)})$, $\wt \bE^+(C)$ becomes a complete valuation ring of rank one with algebraically closed fraction field, called $\wt\bE:=\wt\bE(C)$, of characteristic $p$. 
%
%
%
%
%
Next we define

\begin{tabbing}
$\wt\bA^+\es$\=$:=\es\wt\bA^+(C)$\es\=\kill
$\wt\bA^+$\>$:=\es\wt\bA^+(C)$\>$:=\es W\bigl(\wt\bE^+(C)\bigr)$ and \\[2mm]
$\wt\bA$\>$:=\es\wt\bA(C)$\>$:=\es W\bigl(\wt\bE(C)\bigr)$ the rings of Witt vectors,\\[2mm]
$\wt\bB^+$\>$:=\es\wt\bB^+(C)$\>$:=\es\wt\bA^+(C)[\frac{1}{p}]$ and \\[2mm]
$\wt\bB$\>$:=\es\wt\bB(C)$\>$:=\es\wt\bA(C)[\frac{1}{p}]$ the fraction field of $\wt\bA(C)$.
\end{tabbing}

By Witt vector functoriality there is the Frobenius lift $\phi:=W(\Frob_p)$ 
 on these four rings. For $x\in\wt\bE(C)$ we let $[x]\in\wt\bA(C)$ denote the Teichm\"uller lift. Set $\pi:=[\epsilon]-1$. 
%
%
If $x=\sum_{i=0}^\infty p^i[x_i]\in\wt\bA(C)$ then we set $w_k(x):=\min\{\,v_\bE(x_i):i\le k\,\}$. For $r>0$ let

\begin{tabbing}
$\wt\bA^{(0,r]}\es$\=$:=\es\wt \bA^{(0,r]}(C)$\es\=\kill
$\wt\bA^{(0,r]}$\>$:=\es\wt\bA^{(0,r]}(C)$\>$\DS:=\es\bigl\{\,x\in\wt\bA(C):\es\lim_{k\to+\infty}w_k(x)+{\TS\frac{k}{r}}=+\infty\,\bigr\}$\,,\\[2mm]
$\wt\bB^{(0,r]}$\>$:=\es\wt\bB^{(0,r]}(C)$\>$:=\es\wt\bA^{(0,r]}(C)[\frac{1}{p}]$\,,\\[2mm]
$\wt\bB^\dagger$\>$:=\es\wt\bB^\dagger(C)$\>$\DS:=\es\bigcup_{r\to 0}\wt\bB^{(0,r]}(C)$\,.
\end{tabbing}

\noindent
One has $\wt\bA^+(C)\subset\wt\bA^{(0,r]}(C)\subset\wt\bA^{(0,s]}(C)$ for all $r\ge s>0$.
%
%
On $\wt\bB^{(0,r]}(C)$ there is a valuation defined for $x=\sum_{i\gg-\infty}^\infty p^i[x_i]$ as 

\begin{tabbing}
$\wt \bA^{]0,r]}$\es\=\kill
$v^{(0,r]}(x)\es:=\es\min\{\,w_k(x)+\frac{k}{r}:\,k\in\BZ\,\}\;=\;\min\{\,v_\bE(x_i)+\frac{i}{r}:\,i\in\BZ\,\}$\,.
\end{tabbing}

\noindent
Now one defines $\wt\bB^{]0,r]}(C)$ 
as the Fr\'echet completion of $\wt\bB^{(0,r]}(C)$ 
with respect to the family of semi-valuations $v^{[s,r]}(x):=\min\{\,v^{(0,s]}(x),v^{(0,r]}(x)\,\}$ for $0<s\le r$. This means in concrete terms that a sequence of elements $x_n\in\wt\bB^{(0,r]}(C)$ converges in $\wt \bB^{]0,r]}(C)$ if and only if $\lim_{n\to\infty}v^{[s,r]}(x_{n+1}-x_n)=+\infty$ for all $0<s\le r$. 
Also if $r\ge s$ we let $\wt\bB^{[s,r]}(C)$ be the completion of $\wt\bB^{(0,r]}(C)$ with respect to $v^{[s,r]}$. We view $\wt\bB^{]0,r]}(C)$ as a subring of $\wt\bB^I(C)$ for any closed subinterval $I\subset(0,r]$. 
Let

\begin{tabbing}
$\wt\bB^\dagger_\rig\es$\=$:=\es\wt\bB^\dagger_\rig(C)\es$\=\kill
$\wt\bB^\dagger_\rig$\>$:=\es\wt\bB^\dagger_\rig(C)$\>$\DS:=\es\bigcup_{r\to 0}\wt\bB^{]0,r]}(C)$ 
.
\end{tabbing}


The morphism $\phi$ extends to bicontinuous isomorphisms of rings $\wt\bB^{]0,pr]}(C)\to\wt\bB^{]0,r]}(C)$ and $\wt\bB^{[ps,pr]}(C)\to\wt\bB^{[s,r]}(C)$ defining an automorphism of $\wt\bB^\dagger_\rig(C)$. The rings $\wt\bB^{]0,r]}(C)$ and $\wt\bB^\dagger_\rig(C)$ were studied in detail by Berger~\cite[\S2]{Berger02} who called the first  $\wt\bB_\rig^{\dagger,(p-1)/pr}$ instead. Note that the ring $\wt\bB^\dagger_\rig(C)$ is denoted $\Gamma^\alg_{\an,\con}$ by Kedlaya~\cite[Definition 2.4.8, Remark 2.4.13]{Kedlaya} and is an integral domain.

There is a homomorphism $\theta:\wt\bB^{(0,1]}(C)\to C$ sending $\sum_{i\gg-\infty}^\infty p^i[x_i]$ to $\sum_{i\gg-\infty}^\infty p^ix_i^{(0)}$ which extends by continuity to $\wt\bB^{]0,1]}(C)$ and even to $\wt\bB^{[1,1]}(C)$.
The series 

\begin{tabbing}
$t\es:=\es\log[\epsilon]\es=\es\sum_{n=1}^\infty{\TS\frac{(-1)^{n-1}}{n}}\pi^n\es\in\es\wt\bB^{]0,1]}(C)$ 
\end{tabbing}

\noindent 
satisfies $\phi(t)=pt$ and $t\,\wt\bB^{[1,1]}(C)=\ker\bigl(\theta:\wt\bB^{[1,1]}(C)\to C\bigr)$.
%
%
Let

\begin{tabbing}
$\bB^+_\cris\es$\=$:=\es\bB^+_\cris(C)\es$\= \kill
$\bB^+_\cris$\>$:=\es\bB^+_\cris(C)$ \>be the $p$-adic completion of $\wt\bB^+(C)\bigl[\frac{w^n}{n!}:\; w\in\ker\theta,n\in\BN\bigr]$ and \\[2mm]
$\bB_\cris$\>$:=\es\bB_\cris(C)$\>$:=\es\bB^+_\cris(C)[\frac{1}{t}]$.
\end{tabbing}
The morphism $\phi$ extends to endomorphisms of $\bB^+_\cris(C)$ and $\bB_\cris(C)$. Let

\begin{tabbing}
$\bB^+_\cris\es$\=$:=\es\bB^+_\cris(C)\es$\= \kill
$\DS \wt\bB^+_\rig$\>$:=\es\wt\bB^+_\rig(C)\es:=\es\bigcap_{n\in\BN_0}\phi^n\bB^+_\cris(C)$ and \\[2mm]
$\DS \wt\bB_\rig$\>$:=\es\wt\bB_\rig(C)\es:=\es\wt\bB^+_\rig[\tfrac{1}{t}]$\,.
\end{tabbing}

\noindent
One easily sees that $\wt\bB^+_\rig(C)\subset\wt\bB^{]0,r]}(C)$ for any $r>0$. More precisely, $\wt\bB^{]0,1]}(C)$ equals the $p$-adic completion of $\wt\bB^+_\rig(C)[\frac{p}{[\bar\pi]}]$ and hence is a flat $\wt\bB^+_\rig(C)$-algebra.

%
%

\section{\texorpdfstring{$\phi$}{phi}-Modules} \label{SectPhiModules}
\setcounter{equation}{0}

Let $\LL$ and $C$ be as in the previous section. To shorten notation we will drop the denotation $(C)$ from the rings introduced there.
We recall some definitions and facts from Kedlaya~\cite{Kedlaya}. Let $a$ be a positive integer. 

\begin{definition}\label{Def1.1}
A \emph{$\phi^a$-module} over $\wt\bB^\dagger_\rig$ is a finite free $\wt\bB^\dagger_\rig$-module $\bM$ with a bijective $\phi^a$-semilinear map $\phi_{\bM}:\bM\to\bM$. The rank of $\bM$ as a $\wt\bB^\dagger_\rig$-module is denoted $\rk\bM$. A \emph{morphism of $\phi^a$-modules} is a morphism of the underlying $\wt\bB^\dagger_\rig$-modules which commutes with the $\phi_\bM$'s. We denote the set of morphism between two $\phi^a$-modules $\bM$ and $\bM'$ by $\Hom_{\phi^a}(\bM,\bM')$.
\end{definition}

\begin{example}\label{Def1.3}
Let $c,d\in\BZ$ with $d>0$ and $(c,d)=1$. Define the $\phi^a$-module $\bM(c,d)$ over $\wt\bB^\dagger_\rig$ as $\bM(c,d)=\bigoplus_{i=1}^d \wt\bB^\dagger_\rig \be_i$ equipped with 
\begin{equation}\label{EqStandardBasis}
\phi_\bM(\be_1)\es=\es \be_2\,,\quad\ldots\;,\quad\phi_\bM(\be_{d-1})\es=\es \be_d\,,\quad \phi_\bM(\be_d)\es=\es p^c\be_1\,.
\end{equation}
By abuse of notation we also write $\bM(nc,nd):=\bM(c,d)^{\oplus n}$ for $n\in\BN_{>0}$.
\end{example}

The category of $\phi^a$-modules over $\wt\bB^\dagger_\rig$ is a rigid additive tensor category with unit object $\bM(0,1)$. For a positive integer $b$ there is a \emph{restriction of Frobenius functor} $[b]_\ast$ from $\phi^a$-modules over $\wt\bB^\dagger_\rig$ to $\phi^{ab}$-modules over $\wt\bB^\dagger_\rig$ sending $(\bM,\phi_\bM)$ to $(\bM,\phi_\bM^b)$. Kedlaya proved the following structure theorem.

\begin{theorem}\label{Thm1.4}
(\cite[Theorem 4.5.7]{Kedlaya}) \es Any $\phi^a$-module $\bM$ over $\wt\bB^\dagger_\rig$ is isomorphic to a direct sum of $\phi^a$-modules $\bM(c_i,d_i)$ for uniquely determined pairs $(c_i,d_i)$ up to permutation with $\gcd(c_i,d_i)=1$. It satisfies $\wedge^d\bM\cong\bM(c,1)$ where $c=\sum_ic_i$ and $d=\rk \bM=\sum_i d_i$.
\end{theorem}

\begin{definition}\label{Def2.9a}
If $\det\bM:=\wedge^{\rk\bM}\,\bM\cong\bM(c,1)$ we define the \emph{degree} and the \emph{weight} of $\bM$ as $\deg\bM:=c$ and $\weight\bM:=\frac{\deg\bM}{\rk\bM}$.
\end{definition}

\begin{proposition}\label{Prop2.9}
(\cite[Lemma 3.4.9]{Kedlaya}) \es Every $\phi^a$-submodule $\bM'\subset\bM(c,d)^{\oplus n}$ over $\wt\bB^\dagger_\rig$ satisfies $\weight\bM'\ge\weight\bM(c,d)^{\oplus n}=\frac{c}{d}$.
\end{proposition}

Let us record some facts about the $\bM(c,d)$.

\begin{lemma}\label{Lemma1.3a}
(\cite[Lemmas 4.1.2 and 3.2.4]{Kedlaya}) \es
The $\phi^a$-modules $\bM(c,d)$ over $\wt\bB^\dagger_\rig$ satisfy
\begin{enumerate}
\item $\bM(c,d)\otimes\bM(c',d')\cong\bM(cd'+c'd,dd')$,
\item $\bM(c,d)\dual\;\cong\;\bM(-c,d)$ and
\item $[d]_\ast\bM(c,d)\;\cong\;\bM(c,1)^{\oplus d}$.
\end{enumerate}
\end{lemma}

\begin{definition}\label{Def1.1a}
For a $\phi^a$-module $\bM$ over $\wt\bB^\dagger_\rig$ we define the set of \emph{$\phi^a$-invariants} as
\[
\Koh^0_{\phi^a}(\bM)\es:=\es\{\,x\in\bM:\es\phi_\bM(x)=x\,\}\,.
\]
\end{definition}

It is a vector space over $\Koh^0_{\phi^a}\bigl(\bM(0,1)\bigr)=\BQ_{p^a}:=W(\BF_{p^a})[p^{-1}]$; use \cite[Proposition 3.3.4]{Kedlaya}.
We have $\Hom_{\phi^a}(\bM,\bM')=\Koh^0_{\phi^a}(\bM\dual\otimes\bM')$ for $\phi^a$-modules $\bM$ and $\bM'$ over $\wt\bB^\dagger_\rig$.

\begin{proposition}\label{Prop1.7}
If $a\ge c>0$ then the $\phi^a$-module $\bM(-c,1)$ over $\wt\bB^\dagger_\rig$ satisfies
\[
\Koh^0_{\phi^a}\bigl(\bM(-c,1)\bigr)\es=\es\Bigl\{\, \sum_{\nu\in\BZ}p^{c\nu}\sum_{j=0}^{c-1}p^j\phi^{-a\nu}([x_j]):\es x_0,\ldots,x_{c-1}\in\wt\bE\,,v_\bE(x_j)>0 \,\Bigr\}\,.
\]
\end{proposition}

\begin{proof}
One easily verifies that the series on the right converges (even in $\wt\bB^{]0,r]}$ for all $r>0$) and is a $\phi^a$-invariant of $\bM(-c,1)$.

Conversely let $x\in\wt\bB^\dagger_\rig$ satisfy $x=p^{-c}\phi^a(x)$. Then $x\in \wt\bB^+_\rig$ by \cite[Proposition I.4.1]{Berger04a} and hence possesses an expansion $x=\sum_{j=-\infty}^\infty p^j[x_j]$ with uniquely determined $x_j\in\wt\bE$ satisfying $v_\bE(x_j)>0$ by \cite[Lemme 3.9.17 and Corollaire 3.9.9]{Fourquaux}. The uniqueness implies $\phi^a([x_j])=[x_{j-c}]$ and the proposition follows. 
\end{proof}

Let us make a remark on radii of convergence. Let $\bM^{]0,r]}$ be a free $\wt\bB^{]0,r]}$-module and for $0<s\le r$ let $\bM^{]0,s]}\;:=\;\bM^{]0,r]}\otimes_{\wt\bB^{]0,r]},\iota}\wt\bB^{]0,s]}$ be obtained by base change via the natural inclusion $\iota:\wt\bB^{]0,r]}\hookrightarrow\wt\bB^{]0,s]}$. Let further
\[
\phi_\bM^{]0,rp^{-a}]}:\es \bM^{]0,r]}\otimes_{\wt\bB^{]0,r]},\phi^a}\wt\bB^{]0,rp^{-a}]}\es\isoto\es\bM^{]0,rp^{-a}]}
\]
be an isomorphism of $\wt\bB^{]0,rp^{-a}]}$-modules. We say that a $\phi^a$-module $(\bM,\phi_\bM)$ over $\wt\bB^\dagger_\rig$ is \emph{represented} by the pair $(\bM^{]0,r]},\phi_\bM^{]0,rp^{-a}]})$ if $(\bM,\phi_\bM)=(\bM^{]0,r]}\otimes_{\wt\bB^{]0,r]}}\wt\bB^\dagger_\rig,\phi_\bM^{]0,rp^{-a}]}\otimes\id)$.

\begin{proposition}\label{Prop1.9b}
If $\bM$ is represented by $(\bM^{]0,r]},\phi_\bM^{]0,rp^{-a}]})$ then
\[
\Koh^0_{\phi^a}(\bM)\es=\es\bigl\{\,x\in\bM^{]0,r]}:\es\phi_\bM^{]0,rp^{-a}]}(x\otimes_{\phi^a}1)=x\otimes_\iota 1\,\bigr\}\,.
\]
\end{proposition}

\begin{proof}
The inclusion ``$\supset$'' follows from the inclusion $\bM^{]0,r]}\subset\bM$.

To prove the opposite inclusion ``$\subset$'' let $x\in\Koh^0_{\phi^a}(\bM)$. Since $\wt\bB^\dagger_\rig=\bigcup_{s>0}\wt\bB^{]0,s]}$ there exists a $0<s\le r$ with $x\in\bM^{]0,s]}$. Choose a $\wt\bB^{]0,r]}$-basis of $\bM^{]0,r]}$ and write $x$ with respect to this basis as a vector $v=(v_1,\ldots,v_n)^T\in\bigl(\wt\bB^{]0,s]}\bigr)^{\oplus n}$. Let $A\in\GL_n\bigl(\wt\bB^{]0,rp^{-a}]}\bigr)$ be the matrix by which the isomorphism $\phi_\bM^{]0,rp^{-a}]}$ acts on this basis. The equation $x=\phi_\bM(x)$ translates into $A^{-1}v=\bigl(\phi^a(v_1),\ldots,\phi^a(v_n)\bigr)^T$. Thus if $s\le rp^{-a}$ we find $\phi^a(v_i)\in\wt\bB^{]0,s]}$, whence $v_i\in\wt\bB^{]0,sp^a]}$. Continuing in this way we see that in fact $v_i\in\wt\bB^{]0,r]}$, that is $x\in\bM^{]0,r]}$.
\end{proof}

\begin{proposition}\label{Prop1.6a}
Let $c,d>0$ and let the $\phi$-module $\bM=\bM(-c,d)$ over $\wt\bB^\dagger_\rig(C)$ be represented by $\bM^{]0,1]}= \bigoplus_{i=1}^d \wt\bB^{]0,1]}\be_i$ and $\phi^{]0,p^{-1}]}$ as in \eqref{EqStandardBasis}. Then $\theta:\wt\bB^{]0,1]}(C)\to C$ induces the following exact sequence of $\BQ_p$-vector spaces
\[
\xymatrix @R=0pc @C=1.5pc { 0\ar[r] & \Koh^0_\phi\bigl(\bM(d-c,d)\bigr) \ar[r]^{\TS T_{-c,d}} & \Koh^0_\phi\bigl(\bM(-c,d)\bigr) \ar[r]^{\TS\theta_\bM} & \bigoplus_{i=1}^d C\,\be_i \ar[r] & 0\,,\\
& \DS\sum_{i=1}^d x_i\be_i \ar@{|->}[r] & \DS\sum_{i=1}^d p^{i-1}t\,x_i\be_i\;,\quad\sum_{i=1}^d y_i\be_i\ar@{|->}[r] & \sum_{i=1}^d \theta(y_i)\be_i\,.\hspace{-1cm}
}
\]
In particular $\theta_\bM$ is an isomorphism for $0<c<d$.
\end{proposition}

\begin{proof}
The elements of $\Koh^0_\phi\bigl(\bM(-c,d)\bigr)$ are of the form $\sum_{i=1}^d\phi^{i-1}(y)\be_i$ for $y=p^{-c}\phi^d(y)$ and $y\in\wt\bB^+_\rig\subset\wt\bB^+_{\rm max}\subset\wt\bB^\dagger_\rig$ by \cite[Proposition I.4.1]{Berger04a}. The map $\theta_\bM$ sends this element to $\sum_{i=1}^d\theta(\phi^{i-1}(y))\be_i$. As $\{\phi^{i-1}:1\le i\le d\}=\Gal(\BQ_{p^d}/\BQ_p)$, the sequence is exact on the right by \cite[Lemme 8.16]{Colmez1}. Since $\phi(t)=pt$ and $\wt\bB^\dagger_\rig$ is an integral domain the sequence is well defined and exact on the left. Finally, the exactness in the middle follows from $t\cdot\wt\bB^+_{\rm max}=\{\,y\in\wt\bB^+_{\rm max}:\theta(\phi^n(y))=0\text{ for all }n\in\bN_0\,\}$; see \cite[Proposition 8.10]{Colmez1}. If $0<c<d$ then $\Koh^0_\phi\bigl(\bM(d-c,d)\bigr)=(0)$ by \cite[Proposition 4.1.3]{Kedlaya} and $\theta_\bM$ is an isomorphism.
\end{proof}

\begin{corollary}\label{Cor1.6b}
If $A:\bM\to\bM(c,d)$ is a homomorphism of $\phi$-modules over $\wt\bB^\dagger_\rig$ represented by $A^{]0,1]}:\bM^{]0,1]}\to(\wt\bB^{]0,1]})^{\oplus d}$ with $\theta_{\bM(c,d)}\circ A^{]0,1]}(\bM^{]0,1]})=(0)$ then $A$ factors as $A=T\circ B$ for a morphism $B:\bM\to\bM(c+d,d)$ where $T=T_{c,d}:\bM(c+d,d)\to\bM(c,d)$ is induced by multiplication with $t$ as in Proposition~\ref{Prop1.6a}.
\end{corollary}

\begin{proof}
By Theorem~\ref{Thm1.4} and Lemma~\ref{Lemma1.3a} we write $\bM\cong\sum_{i=1}^r\bM(c_i,d_i)$ and 
\[
\Hom_\phi\bigl(\bM,\,\bM(c,d)\bigr)\;\cong\;\bigoplus_{i=1}^r\Koh^0_\phi\bigl(\bM(cd_i-c_id,dd_i)\bigr)\,.
\]
Then the assertion follows from Proposition~\ref{Prop1.6a}.
\end{proof}

\begin{proposition}\label{Prop1.6}
The $\phi^a$-modules $\bM(c,d)$ over $\wt\bB^\dagger_\rig$ satisfy
\begin{enumerate}
\item \label{Prop1.6_a}
$\dim_{\BQ_{p^a}}\Koh^0_{\phi^a}\bigl(\bM(c,d)\bigr)=\left\{\begin{array}{cl}
0&\quad\text{if }c>0\\
1&\quad\text{if }c=0,d=1\\
\infty&\quad\text{if }c<0
\end{array}\right.$
\item \label{Prop1.6_b}
$\dim_{\BQ_{p^a}}\Hom_{\phi^a}\bigl(\bM(c,d),\bM(c',d')\bigr)=\left\{\begin{array}{cl}
0&\quad\text{if }c/d<c'/d'\,,\\
d^2&\quad\text{if }c=c',d=d'\,,\\
\infty&\quad\text{if }c/d>c'/d'\,.
\end{array}\right.$
\item \label{Prop1.6_c}
$\End_{\phi^a}\bigl(\bM(c,d)\bigr)$ is a central division algebra over $\BQ_{p^a}$ of
dimension $d^2$ and Hasse invariant $\frac{c}{d}$ if $(c,d)=1$.
\end{enumerate}
\end{proposition}

\begin{proof}
\ref{Prop1.6_a} The assertions for $c>0$ and $c=0,d=1$ were proved by Kedlaya \cite[Propositions 4.1.3 and 3.3.4]{Kedlaya}. The case $c<0$ follows from Proposition~\ref{Prop1.6a} since $\dim_{\BQ_{p^a}}C=\infty$. Note that Proposition~\ref{Prop1.6a} was proved only for $a=1$ but the surjectivity of $\theta_\bM$, which only relies on \cite[Lemme 8.16]{Colmez1}, holds also for $a>1$.

\smallskip\noindent
\ref{Prop1.6_b} follows from \ref{Prop1.6_a} and 
\[
\Hom_{\phi^a}\bigl(\bM(c,d)\,,\,\bM(c',d')\bigr)\;=\;\Koh^0_{\phi^a}\bigl(\bM(-c,d)\otimes\bM(c',d')\bigr)\;\cong\;\Koh^0_{\phi^a}\bigl(\bM(c'd-cd',dd')\bigr)\,.
\]

\smallskip\noindent
\ref{Prop1.6_c} was proved in \cite[Remark 4.1.5]{Kedlaya}. The Hasse invariant can be computed as in \cite[Proposition 8.6]{HP}.
\end{proof}

%
%

\section{Constructing \texorpdfstring{$\phi$}{phi}-Modules from Filtered Isocrystals} \label{SectMinuscule}
\setcounter{equation}{0}

We keep the notation from Section~\ref{SectFontaine'sRings}. In the situation where the minimal and maximal Hodge-Tate weights differ by at most $1$ we will associate with any analytic point $\mu\in\breve\CF^\an$ a $\phi$-module $\bM_\mu$ over $\wt\bB^\dagger_\rig(C)$ where $C$ is the completion of an algebraic closure of $\CH(\mu)$. In case $\CH(\mu)/\breve E$ finite this construction parallels Berger's construction~\cite[\S II]{Berger08} that works for arbitrary Hodge-Tate weights and even produces an ``$\CH(\mu)$-rational'' version of $\bM_\mu$. We begin with the following

\begin{proposition}\label{Prop1.12}
Let $\bN$ be a $\phi$-module over $\wt\bB^\dagger_\rig$ represented by a free $\wt\bB^{]0,1]}$-module $\bN^{]0,1]}$ and an isomorphism $\phi_\bN^{]0,p^{-1}]}:\bN^{]0,1]}\otimes_{\wt\bB^{]0,1]},\phi}\wt\bB^{]0,p^{-1}]}\isoto\bN^{]0,p^{-1}]}$. Consider the morphism $\theta_\bN=\id_\bN\otimes\theta:\bN^{]0,1]}\to \bN^{]0,1]}\otimes_{\wt\bB^{]0,1]},\theta}C=:W_\bN$ and let $W_\bM$ be a $C$-subspace of $W_\bN$.
Then there exists a uniquely determined $\phi$-submodule $\bM\subset\bN$ over $\wt\bB^\dagger_\rig$ with $t\bN\subset\bM$ which is represented by a $\wt\bB^{]0,1]}$-submodule $t\bN^{]0,1]}\subset\bM^{]0,1]}\subset\bN^{]0,1]}$ such that $\theta_\bN(\bM^{]0,1]})=W_\bM$.
\end{proposition}

\begin{proof}
If $I$ is a closed subinterval of $(0,1]$ set $\bN^I:=\bN^{]0,1]}\otimes_{\wt\bB^{]0,1]}}\wt\bB^I$. Also if $I=[s,r]$ we let $pI:=[ps,pr]$. The isomorphism $\phi_\bN^{]0,p^{-1}]}$ induces an isomorphism $\phi_\bN^I:\bN^{pI}\otimes_{\wt\bB^{pI},\phi}\wt\bB^I\isoto\bN^I$.
Let $I_n:=[p^{-n-1},p^{-n}]$ for $n\in\BN_0$. We define $\bM'{}^{I_0}$ as the preimage of $W_\bM$ in $\bN^{I_0}$ under the canonical morphism $\bN^{I_0}\to\bN^{I_0}\otimes_{\wt\bB^{I_0},\theta}C=W_\bN$ and we let $\bM^{I_0}$ be the intersection of $\bM'{}^{I_0}$ with $\phi_\bN^{[p^{-1},p^{-1}]}\bigl(\bM'{}^{I_0}\otimes_{\wt\bB^{I_0},\phi}\wt\bB^{[p^{-1},p^{-1}]}\bigr)$ inside $\bN^{[p^{-1},p^{-1}]}$, that is $\bM^{I_0}$ equals
\[
\ker\Bigl(\;\bM'{}^{I_0}\es\into\es\bN^{[p^{-1},p^{-1}]}\es\onto\es\bN^{[p^{-1},p^{-1}]}/\phi_\bN^{[p^{-1},p^{-1}]}\bigl(\bM'{}^{I_0}\otimes_{\wt\bB^{I_0},\phi}\wt\bB^{[p^{-1},p^{-1}]}\bigr)\;\Bigr)\,.
\]
Since $\wt\bB^{I_0}$ is a principal ideal domain by \cite[Proposition 2.6.8]{Kedlaya} we see that $\bM^{I_0}$ is a free $\wt\bB^{I_0}$-submodule of $\bN^{I_0}$ of full rank.
For $n\ge1$ we define $\bM^{I_n}$ as the image of $\bM^{I_0}\otimes_{\wt\bB^{I_0},\phi^n}\wt\bB^{I_n}$ under the isomorphism
\[
\phi_\bN^{I_n}\circ\ldots\circ\phi_\bN^{I_1}:\bN^{I_0}\otimes_{\wt\bB^{I_0},\phi^n}\wt\bB^{I_n}\isoto\bN^{I_n}.
\]
In the terminology of Kedlaya~\cite[\S 2.8]{Kedlaya} the collection $\bM^{I_n}$ for $n\in\BN_0$ defines a \emph{vector bundle over $\wt\bB^{]0,1]}$} which by \cite[Theorem 2.8.4]{Kedlaya} corresponds to a free $\wt\bB^{]0,1]}$-submodule $\bM^{]0,1]}$ of $\bN^{]0,1]}$. By construction it satisfies $t\bN^{]0,1]}\subset\bM^{]0,1]}$ and $\phi_\bN^{]0,p^{-1}]}$ restricts to an isomorphism on $\bM^{]0,1]}$. This makes $\bM:=\bM^{]0,1]}\otimes_{\wt\bB^{]0,1]}}\wt\bB^\dagger_\rig$ into a $\phi$-module with the desired properties. 
Clearly $\bM^{]0,1]}$ is uniquely determined by the subspace $W_\bM\subset W_\bN$ and by the requirements that $t\bN^{]0,1]}\subset\bM^{]0,1]}\subset\bN^{]0,1]}$ and that $\phi_\bN^{I_n}\bigl(\bM^{I_{n-1}}\otimes_{\wt\bB^{I_{n-1}},\phi}\wt\bB^{I_n}\bigr)=\bM^{I_n}$.
\end{proof}

Now assume that $\LL$ is a (not necessarily finite) extension of $K_0$ and let $\hhh$ be an integer. Let $(D,\phi_D)$ be an $F$-isocrystal over $\BF_p^{\,\alg}$ and let $Fil^\hhh D_\LL$ be an $\LL$-subspace of $D_\LL=D\otimes_{K_0}\LL$. Let $Fil^{\hhh-1}D_\LL=D_\LL$ and $Fil^{\hhh+1}D_\LL=(0)$. Then $\ulD=(D,\phi_D,Fil^\bullet D_\LL)$ is a filtered isocrystal over $\LL$ with Hodge-Tate weights $-\hhh$ and $-\hhh+1$.  

\begin{definition}\label{Def1.12b}
We set $\bD^{]0,1]}:=D\otimes_{K_0}\wt\bB^{]0,1]}$ and $\bD(\ulD):=(\bD,\phi_\bD):=(D,\phi_D)\otimes_{K_0}\wt\bB^\dagger_\rig$. Multiplication with $t^\hhh $ defines natural isomorphisms $\id_D\otimes t^\hhh :t^{-\hhh }\bD^{]0,1]}\isoto D\otimes_{K_0}\wt\bB^{]0,1]}$ and $t^{-\hhh }\bD(\ulD)\isoto(D,p^{-\hhh }\phi_D)\otimes_{K_0}\wt\bB^\dagger_\rig$. We let $\bM(\ulD)$ be the $\phi$-submodule of $t^{-\hhh }\bD$ over $\wt\bB^\dagger_\rig$ represented by the $\wt\bB^{]0,1]}$-submodule $\bM^{]0,1]}$ of $t^{-\hhh }\bD^{]0,1]}$ with $\theta_{t^{-\hhh }\bD}(\bM^{]0,1]})=(Fil^\hhh  D_\LL)\otimes_\LL C$ under the map $\theta_{t^{-\hhh }\bD}:t^{-\hhh }\bD^{]0,1]}\to t^{-\hhh }\bD^{]0,1]}\otimes_{\wt\bB^{]0,1]},\theta}C=D_C$, whose existence was established in Proposition~\ref{Prop1.12}. 
\end{definition}

The following lemma is immediate from the definition.

\begin{lemma}\label{Lemma4.2'}
For all integers $m$ consider the Tate object 
\[
\BOne(m):=\bigl(D=K_0\,,\,\phi_D=p^m\,,\,K_0=Fil^m\supset Fil^{m+1}=(0)\bigr). 
\]
Then the canonical isomorphism $\id_D\otimes t^m:\bD\bigl(\ulD\otimes\BOne(m)\bigr)\isoto t^m\bD(\ulD)$ induces an isomorphism $\bM\bigl(\ulD\otimes\BOne(m)\bigr)\cong\bM(\ulD)$.\qed
\end{lemma}

The construction is compatible with dualizing in the following sense.

\begin{proposition} \label{Prop4.10}
There is a canonical isomorphism $\bM(\ulD\dual)\cong\bM(\ulD)\dual$, where the $\phi$-module $\bM(\ulD)\dual:=\CHom_\phi\bigl(\bM(\ulD),\bM(0,1)\bigr)$ is the inner hom of $\phi$-modules into the unit object $\bM(0,1)$. 
\end{proposition}

\begin{proof}
If $\ulD$ is a filtered isocrystal with $Fil^{\hhh -1}D_\LL=D_\LL\supset Fil^\hhh D_\LL\supset Fil^{\hhh +1}D_\LL=(0)$ then $\ulD\dual$ has $Fil^{-\hhh }D_\LL\dual=D_\LL\dual\supset Fil^{-\hhh +1}D_\LL\dual\supset Fil^{-\hhh +2}D_\LL\dual=(0)$ because $Fil^iD_\LL\dual:=\{\lambda\in D_\LL\dual:\lambda(Fil^{1-i}D_\LL)=(0)\}$ for all $i$. Thus the canonical pairing $D\times D\dual\to K_0$ induces a pairing $\psi$ as follows
\[
\xymatrix {
\bM(\ulD)^{]0,1]} \ar@{^{ (}->}[d]&*=<-3em,0em>\objectbox{\times} & \bM(\ulD\dual)^{]0,1]} \ar@{^{ (}->}[d] \\
t^{-\hhh }\bD(\ulD)^{]0,1]} \ar[d]^\cong &*=<-3em,0em>\objectbox{\times}& t^{\hhh -1}\bD(\ulD\dual)^{]0,1]}\ar[d]^\cong \ar[r]^{\TS\psi} & \bM(-1,1)^{]0,1]}\ar[d]^\cong\\
(D,p^{-\hhh }\phi_D)\otimes_{K_0}\wt\bB^{]0,1]} \ar[d]^{\theta_{t^{-\hhh }\bD(\ulD)}} &*=<-3em,0em>\objectbox{\times}& (D\dual,p^{\hhh -1}\phi_{D\dual})\otimes_{K_0}\wt\bB^{]0,1]} \ar[d]^{\theta_{t^{\hhh -1}\bD(\ulD\dual)}} \ar[r] & (K_0,p^{-1}\phi_{\SSC K_0})\otimes_{K_0}\wt\bB^{]0,1]} \ar[d]^{\theta_{\bM(-1,1)}} \\
D_C &*=<-3em,0em>\objectbox{\times}& D_C\dual \ar[r] & C
}
\]
Since $\theta_{t^{-\hhh }\bD(\ulD)}\bigl(\bM(\ulD)^{]0,1]}\bigr)=Fil^\hhh D_C$ and $\theta_{t^{\hhh -1}\bD(\ulD\dual)}\bigl(\bM(\ulD\dual)^{]0,1]}\bigr)=Fil^{-\hhh +1}D_C\dual$ we obtain $\theta_{\bM(-1,1)}\circ\psi\bigl(\bM(\ulD)^{]0,1]}\otimes_{\wt\bB^{]0,1]}}\bM(\ulD\dual)^{]0,1]}\bigr)=(0)$. By Corollary~\ref{Cor1.6b} the pairing $\bM(\ulD)\otimes_{\wt\bB^\dagger_\rig}\bM(\ulD\dual)\to\bM(-1,1)$ factors through $\bM(0,1)$. This implies $\bM(\ulD\dual)\into\bM(\ulD)\dual$. Since both $\phi$-modules have rank equal to $\dim D$ and by Theorem~\ref{Thm1.13} below their degree satisfies
\[
\deg\bM(\ulD\dual)\;=\;t_N(\ulD\dual)-t_H(\ulD\dual)\;=\;-t_N(\ulD)+t_H(\ulD)\;=\;-\deg\bM(\ulD)\;=\;\deg\bM(\ulD)\dual\,,
\]
\cite[Lemma 3.4.2]{Kedlaya} yields $\bM(\ulD\dual)\isoto\bM(\ulD)\dual$. 
\end{proof}

Unfortunately we do not know how to construct $\bM(\ulD)$ for filtered isocrystals $\ulD$ with Hodge-Tate weights other than $-\hhh $ and $-\hhh +1$. Therefore we cannot make $\ulD\mapsto \bM(\ulD)$ into a tensor functor and we do not follow Berger's argument \cite[Th\'eor\`eme IV.2.1]{Berger08} to prove the next theorem.

\begin{theorem}\label{Thm1.13}
\mbox{ } $\deg\bM(\ulD)\; =\; t_N(\ulD) - t_H(\ulD)$.
\end{theorem}

\begin{proof}
By the Dieudonn\'e-Manin classification~\cite{Manin} 
there exists an $F$-isocrystal $(D',\phi_{D'})$ over $\BF_p^{\,\alg}$ of rank one with $\phi_{D'}=p^{t_N(\ulD)}\cdot\phi$ which is isomorphic to $\det (D,\phi_D)$. So by construction $\deg\bD(\ulD)=t_N(\ulD)$. Moreover, 
\begin{eqnarray*}
t_H(\ulD) & = & (\hhh -1)\cdot\dim_\LL(D_\LL/Fil^\hhh  D_\LL) \;+\; \hhh \cdot\dim_\LL Fil^\hhh D_\LL\\[2mm]
& = & \hhh \cdot\dim_\LL D_\LL \;-\;\dim_C(D_\LL/Fil^\hhh D_\LL)\otimes_\LL C\\[2mm]
& = & \hhh \cdot\dim_{K_0}D\;-\; \dim_C(t^{-\hhh }\bD^{]0,1]}/\bM^{]0,1]})\otimes_{\wt\bB^{]0,1]},\theta}C\,.
\end{eqnarray*}
The following lemma implies that 
\begin{eqnarray*}
\deg t^i\bD(\ulD)-\deg t^{i-1}\bD(\ulD)&=&\dim_{K_0}D \quad\text{for all $i$ and}\\[2mm]
\deg\bM(\ulD)-\deg t^{-\hhh }\bD(\ulD)&=&\dim_C(t^{-\hhh }\bD^{]0,1]}/\bM^{]0,1]})\otimes_{\wt\bB^{]0,1]},\theta}C
\end{eqnarray*}
and this proves the theorem.
\end{proof}

\begin{lemma}\label{Lemma2.8b}
Let $\bM_1$ and $\bM_2$ be $\phi$-modules over $\wt\bB^\dagger_\rig$ represented by $\wt\bB^{]0,1]}$-modules $\bM_i^{]0,1]}$. Assume that $\bM_1^{]0,1]}\supset\bM_2^{]0,1]}\supset t\bM_1^{]0,1]}$. Then
\[
\deg\bM_2-\deg\bM_1\es=\es\dim_C (\bM_1^{]0,1]}/\bM_2^{]0,1]})\otimes_{\wt\bB^{]0,1]},\theta}C\,.
\]
\end{lemma}

\begin{proof}
Clearly the equality holds for $\bM_2=t\bM_1\cong\bM_1\otimes\bM(1,1)$ since $\deg t\bM_1-\deg\bM_1=\rk\bM_1$. We claim that it suffices to prove the inequality 
\begin{equation}\label{EqIneqDeg}
\deg\bM_2-\deg\bM_1\es\ge\es\dim_C W
\end{equation}
where we abbreviate $W:=(\bM_1^{]0,1]}/\bM_2^{]0,1]})\otimes_{\wt\bB^{]0,1]},\theta}C$. Indeed we apply the inequality to the two inclusions $\bM_1\supset\bM_2\supset t\bM_1$ and $\bM_2\supset t\bM_1\supset t\bM_2$ and conclude using the exact sequence of $\wt\bB^{]0,1]}$-modules
\[
0\es\longto\es \bM_2^{]0,1]}/t\bM_1^{]0,1]} \es\longto\es \bM_1^{]0,1]}/t\bM_1^{]0,1]} \es\longto\es \bM_1^{]0,1]}/\bM_2^{]0,1]} \es\longto\es 0\,.
\]

To prove the inequality (\ref{EqIneqDeg}) we argue by induction on $\dim_C W$. Let $\dim_C W=1$. Since $\det\bM_1\supset\det\bM_2$ we know that $\deg\bM_2\ge\deg\bM_1$ from Proposition~\ref{Prop1.6}. If we had $\deg\bM_2=\deg\bM_1$ then $\bM_2=\bM_1$ by \cite[Lemma 3.4.2]{Kedlaya}. So $\deg\bM_2-\deg\bM_1\ge1=\dim_C W$ as desired.

Let now $\dim_C W>1$ and choose a $C$-subspace $W'$ of dimension $1$ of $W$.
By Proposition~\ref{Prop1.12} there is a unique $\wt\bB^\dagger_\rig$-submodule $\bM_2\subset\bM_3\subset\bM_1$ corresponding to $W'$.
By induction $\deg\bM_3-\deg\bM_1\ge \dim_C (W/W')$ and $\deg\bM_2-\deg\bM_3\ge\dim_CW'$ proving the lemma.
\end{proof}

\begin{remark}\label{RemPottharst}
 As J.~Pottharst pointed out to us, one could nevertheless follow \cite[Th\'eor\`eme IV.2.1]{Berger08} to prove Theorem~\ref{Thm1.13} by working in the exact tensor category $\ul{\mathrm{MP}}^\phi_C$ of triples $(D,\phi_D,\CL)$ with $(D,\phi_D)$ an $F$-isocrystal over $\BF_p^{\,\alg}$ and $\CL$ a $\bB_\dR^+$-lattice inside $D\otimes_{K_0}\bB_\dR$, where $\bB_\dR^+=\invlim \wt\bB^+/(\ker\theta)^i$, $\bB_\dR=\bB_\dR^+[t^{-1}]$, and $\theta:\bB_\dR^+/t\bB_\dR^+\isoto C$. Note that to $\ul{\mathrm{MP}}^\phi_C$ also the functions $t_N$ and $t_H$ extend via $t_N(D,\phi_D,\CL):=t_N(D,\phi_D)$ and $t_H(D,\phi_D,\CL):=m$ if $\wedge^{\dim D}\CL=t^{-m}\wedge^{\dim D}D\otimes_{K_0}\bB_\dR^+$.
From $\ul{\mathrm{MP}}^\phi_C$ there are exact tensor functors $F_1$ to the category of filtered isocrystals over $C$ by defining 
\[
Fil^i D_C\;=\;\bigl(t^i\CL\cap D\otimes_{K_0}\bB_\dR^+\bigr)\,\bigr/\,\bigl(t^i\CL\cap t\cdot D\otimes_{K_0}\bB_\dR^+\bigr)\;=\;\id_D\otimes\theta\bigl(t^i\CL\cap D\otimes_{K_0}\bB_\dR^+\bigr)\,,
\]
and $F_2$ to the category of $\phi$-modules over $\wt\bB^\dagger_\rig$ similar to our construction in Proposition~\ref{Prop1.12}.
More precisely $F_2$ assigns to $(D,\phi_D,\CL)$ the $\phi$-module
\begin{eqnarray*}
\bM(D,\phi_D,\CL)&:=&\bM(D,\phi_D,\CL)^{]0,1]}\otimes_{\wt\bB^{]0,1]}}\wt\bB^\dagger_\rig\qquad\text{where}\\[2mm]
\bM(D,\phi_D,\CL)^{]0,1]}&:=&\bigl\{\,m\in D\otimes_{K_0}\wt\bB^{]0,1]}[t^{-1}]:\es \phi_D^i(m\otimes_{\phi^i}1)\in \CL\text{ for all }i\le0\,\bigr\}\,.
\end{eqnarray*}
Restricted to the subcategories where the Hodge-Tate weights are $-\hhh $ and $-\hhh +1$, respectively where $t^{-\hhh +1}D\otimes_{K_0}\bB_\dR^+\subset \CL\subset t^{-\hhh }D\otimes_{K_0}\bB_\dR^+$, the functor $F_1$ is an equivalence. The inverse functor $F_1^{-1}$ assigns to the filtered isocrystal $\ulD$ the $\bB^+_\dR$-lattice $\CL$ such that $t^\hhh \CL$ is the preimage of $Fil^\hhh  D_C$ under the map $\id_D\otimes\theta:D\otimes_{K_0}\bB^+_\dR\to D_C$. Our functor $\ulD\mapsto\bM(\ulD)$ from Definition~\ref{Def1.12b} is then the composition $F_2\circ F_1^{-1}$.
\end{remark}

\begin{remark}\label{Rem1.12c}
The construction of $\bM(\ulD)$ is functorial with respect to automorphisms of $C$. Therefore $\bM(\ulD)$ gives rise to a Galois representation
\[
\Gal\bigl(\CH(x)^\alg/\CH(x)\bigr)\es\longto\es\GL\Bigl(\Koh^0_\phi\bigl(\bM(\ulD)\bigr)\Bigr)\,.
\]
Of course this representation is of most interest when $\bM(\ulD)\cong\bM(0,1)^{\oplus\dim D}$ since then $\Koh^0_\phi\bigl(\bM(\ulD)\bigr)\cong\BQ_p^{\,\oplus\dim D}$.
\end{remark}

If $\LL$ is a discretely valued extension of $K_0$ one can check that $\bM(\ulD)$ equals the $\phi$-module over $\wt\bB^\dagger_\rig$ constructed by Berger~\cite[\S II]{Berger08}. One of Berger's main theorems is the following criterion.

\begin{theorem}\label{Thm1.5}
(\cite[\S IV.2]{Berger08}) \es 
Let $\LL$ be a discretely valued extension of $K_0$. Then
$\ulD$ is weakly admissible if and only if $\bM(\ulD)\cong\bM(0,1)^{\oplus\dim D}$. In this case the crystalline Galois representation associated with $\ulD$ is the representation $\Koh^0_\phi\bigl(\bM(\ulD)\bigr)$ from Remark~\ref{Rem1.12c}.
\end{theorem}

The theorem indicates what the local system on $\breve\CF_b^a$ could be. We will discuss this in detail in the next sections. Let us first make explicit what happens if $\bM(\ulD)\not\cong\bM(0,1)^{\oplus\dim D}$.

\begin{proposition}\label{Prop1.8}
Assume that $t_N(\ulD) = t_H(\ulD)$ and that the Hodge-Tate weights of $\ulD$ are $-\hhh $ and $-\hhh +1$. Then $\bM(\ulD)\not\cong\bM(0,1)^{\oplus\dim D}$ if and only if for some (any) integer $e\ge(\dim D)-1$ there exists a non-zero 
$x\in\Koh^0_{\phi^e}\bigl(\bM(1,1)\otimes[e]_\ast(t^{-\hhh }\bD)\bigr)$ with $\theta_{t^{-\hhh }\bD}\bigl(\phi_{t^{-\hhh }\bD}^m(x)\bigr)\in (Fil^\hhh  D_\LL)\otimes_\LL C$ for all $m=0,\ldots,e-1$.
\end{proposition}

Note that $t^{-\hhh }\bD$ can be represented by $(D\otimes_{K_0}\wt\bB^{]0,r]},p^{-\hhh }b\!\cdot\!\phi)$ for arbitrarily large $r>0$. Hence by Proposition~\ref{Prop1.9b} the element $x$ actually belongs to $D\otimes_{K_0}\wt\bB^{]0,r]}$ and the expression $\theta_{t^{-\hhh }\bD}\bigl(\phi_{t^{-\hhh }\bD}^m(x)\bigr)$ makes sense.

\begin{proof}
Let $\bM(\ulD)\not\cong\bM(0,1)^{\oplus\dim D}$ and let $e\ge(\dim D)-1$ be any integer. By Theorems~\ref{Thm1.13} and \ref{Thm1.4} there is a $\phi$-module $\bM(c,d)$ over $\wt\bB^\dagger_\rig$ with $c<0$ and $d<\dim D$ which is a summand of $\bM(\ulD)$. Thus by Proposition~\ref{Prop1.6} there exists a non-zero morphism of $\phi$-modules $f:\bM(-1,e)\to\bM(c,d)\subset\bM(\ulD)\subset t^{-\hhh }\bD$. Let $\be_1,\ldots,\be_e$ be a basis of $\bM(-1,e)$ satisfying (\ref{EqStandardBasis}) and let $x$ be the image of $\be_1$ in $t^{-\hhh }\bD$ under $f$. Then $\phi_{t^{-\hhh }\bD}^e(x)=p^{-1}x$, that is $x\in\Koh^0_{\phi^e}\bigl(\bM(1,1)\otimes[e]_\ast(t^{-\hhh }\bD)\bigr)$. Moreover, $f(\be_{m+1})=\phi_{t^{-\hhh }\bD}^m(x)$ in $t^{-\hhh }\bD$ for $0\le m<e$. Now the fact that the morphism $f$ factors through $\bM(\ulD)$ amounts by Definition~\ref{Def1.12b} to $\theta_{t^{-\hhh }\bD}\bigl(\phi_{t^{-\hhh }\bD}^m(x)\bigr)\in (Fil^ \hhh D_\LL)\otimes_\LL C$.

Conversely assume that for some integer $e\ge(\dim D)-1$ there exists a non-zero element $x$ in $\Koh^0_{\phi^e}\bigl(\bM(1,1)\otimes[e]_\ast(t^{-\hhh }\bD)\bigr)$ with $\theta\bigl(\phi_{t^{-\hhh }\bD}^m(x)\bigr)\in (Fil^\hhh  D_\LL)\otimes_\LL C$ for all $0\le m<e$. Define the non-trivial morphism of $\phi$-modules $f:\bM(-1,e)\to t^{-\hhh }\bD$ by $f(\be_{m+1}):=\phi_{t^{-\hhh }\bD}^m(x)$ for $0\le m<e$. Since $\theta_{t^{-\hhh }\bD}\bigl(\phi_{t^{-\hhh }\bD}^m(x)\bigr)\in (Fil^\hhh  D_\LL)\otimes_\LL C$ the morphism $f$ factors through $\bM(\ulD)$ by Definition~\ref{Def1.12b}. By Proposition~\ref{Prop1.6} this implies $\bM(\ulD)\not\cong\bM(0,1)^{\oplus\dim D}$.
\end{proof}

%
%

\section{The Minuscule Case} \label{SectOpenness}
\setcounter{equation}{0}

We retain the notation from Section~\ref{SectPeriodSpaces}. In particular $G$ is a reductive group, $b\in G(K_0)$, and $\{\mu\}$ is a conjugacy class of cocharacters of $G$. Recall from \cite[1.7]{RZ} that there is a slope quasi-cocharacter $\nu\in\Hom_{K_0}(\BG_m,G)\otimes_\BZ \BQ$ associated with $b$. We make the following assumption on the triple $(G,b,\{\mu\})$:
\begin{eqnarray} \label{Eq2.2}
&\text{$\nu^\#=\mu^\#$ in $(\pi_1(G)_{\Gal(\ol\BQ_p/\BQ_p)})\otimes_\BZ\BQ$, and}\nonumber \\
&\text{there is an $\hhh\in\BZ$ and a faithful $\BQ_p$-rational representation $V$ of $G$,}\\
&\text{ such that all the weights of $\{\mu\}$ on $V$ are $\hhh-1$ and $\hhh$, and} \nonumber\\
&\text{ the isocrystal $\ulD_{b,\mu}(V)$ from (\ref{EqIsocrystal}) has Newton slopes in the interval $[\hhh-1,\hhh]$.} \nonumber 
\end{eqnarray}
The significance of $\nu^\#=\mu^\#$ is explained by the following

\begin{lemma}\label{LemmaKottwitzCond}
Let $\nu$ be the slope quasi-cocharacter in $\Hom_{K_0}(\BG_m,G)\otimes_{\BZ}\BQ$ associated with $b$; see \cite[1.7]{RZ}, and let $\rho:G\to\GL(V)$ be in $\Rep_{\BQ_p}\!\!\!G$. Then under the canonical identification $\Hom(\BG_m,\BG_m)=\BZ$  we have
\[
\TS\det_V\circ\rho\circ\nu\;=\;t_N\bigl(\ulD_{b,\mu}(V)\bigr)\qquad\text{and}\qquad\det_V\circ\rho\circ\mu\;=\;t_H\bigl(\ulD_{b,\mu}(V)\bigr)\,.
\]
In particular, if the images $\nu^\#$ and $\mu^\#$ of $\nu$ and $\mu$ in $(\pi_1(G)_{\Gal(\ol\BQ_p/\BQ_p)})\otimes_\BZ\BQ$ coincide then $t_N\bigl(\ulD_{b,\mu}(V)\bigr)=t_H\bigl(\ulD_{b,\mu}(V)\bigr)$ for all $V\in\Rep_{\BQ_p}\!\!\!G$.
\end{lemma}

\begin{proof}
The statement about $t_N$ follows from the construction of $\nu$ in \cite[\S4.2]{Kottwitz}. The statement about $t_H$ is immediate from the definition of $\ulD_{b,\mu}(V)$. If $\nu^\#=\mu^\#$ holds in $(\pi_1(G)_{\Gal(\ol\BQ_p/\BQ_p)})\otimes_\BZ\BQ$ then $(\rho\circ\nu)^\#=(\rho\circ\mu)^\#$ in $(\pi_1(\GL(V))_{\Gal(\ol\BQ_p/\BQ_p)})\otimes_\BZ\BQ\cong\BQ$. Under the last isomorphism we have $(\rho\circ\nu)^\#=\det\circ\rho\circ\nu$ and likewise for $\mu$. This proves the lemma.
\end{proof}

If $(G,\{\mu\})$ satisfies \eqref{Eq2.2} with $\hhh=1$ then $\mu$ is a minuscule coweight by Serre~\cite[\S3]{Serre}. Furthermore, assume that $G$ does not possess a proper normal subgroup through which $\mu$ factors. Then all weights of $\mu$ in $V$ form one orbit under the Weyl group, the group $G$ has no simply connected factor group, and all irreducible components of $G^\der$ are of type $A_n , B_n , C_n$ or $D_n$. Moreover, in case $V$ is an irreducible representation, $\dim V$ is even, except possibly when all irreducible factors of $G^\der$ are of type $A$.

Assume now that $(G,b,\{\mu\})$ satisfies condition (\ref{Eq2.2}) and let $\rho:G\hookrightarrow\GL(V)=:G'$ be the faithful representation from \eqref{Eq2.2}. This defines a closed embedding $\CF\hookrightarrow\CF':=\CF lag(V)$ into the flag variety from (\ref{Eq2.0}). Here $\CF lag(V)$ is actually the Grassmannian over $K_0$ of $c$-dimensional subspaces of $V_{K_0}=V\otimes_{\BQ_p}K_0$, where $c$ is the multiplicity of the weight $\hhh$ of $\mu$ on $V$. By \cite[1.18]{RZ} we obtain a closed embedding of $\breve E$-analytic period spaces $\breve\CF^{wa}_b\hookrightarrow\breve\CF'{}_{\!\!\!\!\rho(b)}^{wa},\,\mu\mapsto\rho\circ\mu$.

Let $\mu\in\breve\CF^\an$ be an analytic point. Let $\LL=\CH(\mu)$ and let $C$ be the completion of an algebraic closure of $\LL$. Let $\ulD_\mu:=\ulD_{b,\mu}(V):=(V_{K_0},b\!\cdot\!\phi,Fil^\bullet_\mu D_\LL)$ be the filtered isocrystal from (\ref{EqIsocrystal}). In particular $Fil_\mu^{\hhh-1}D_\LL=D_\LL$ and $Fil_\mu^{\hhh+1} D_\LL=(0)$. Now we let $\bM_\mu:=\bM_{b,\mu}(V):=\bM\bigl(\ulD_{b,\mu}(V)\bigr)$ be the $\phi$-module over $\wt\bB^\dagger_\rig(C)$ from Definition~\ref{Def1.12b}. Note that $\deg\bM_\mu=0$ under our assumption~\eqref{Eq2.2} by Theorem~\ref{Thm1.13} and Lemma~\ref{LemmaKottwitzCond}. We define
\begin{eqnarray} \label{EqFAdm}
\breve\CF_b^0&:=&\breve\CF_b^0(V)\es:=\es\bigl\{\,\mu\in\breve\CF^\an\text{ analytic points} : \es\bM_{b,\mu}(V)\cong\bM(0,1)\otimes_{\BQ_p}\!V\,\bigr\}.
\end{eqnarray}

If $b'=g\,b\,\phi(g^{-1})$ for some $g\in G(K_0)$ one easily checks that $\mu\mapsto g\mu g^{-1}=:\mu'$ maps $\breve\CF_b^0$ isomorphically onto $\breve\CF_{b'}^0$, because $g$ induces isomorphisms $\ulD_{b'\!,\mu'}(V)\cong\ulD_{b,\mu}(V)$ and $\bM_{b'\!,\mu'}(V)\cong\bM_{b,\mu}(V)$. In particular $\breve\CF^0_b$ is invariant under the group $J(\BQ_p)$ from (\ref{EqJ}).

\begin{proposition}\label{PropDOR}
The set $\breve\CF_b^0$ is independent of the representation $V$ of $G$ with (\ref{Eq2.2}).
\end{proposition}

\begin{proof}
Let $\rho:G\hookrightarrow\GL(V)$ and $\tilde\rho:G\hookrightarrow\GL(\wt V)$ be two representations of $G$ satisfying (\ref{Eq2.2}) for integers $\hhh$ respectively $\tilde\hhh$. Let $\mu\in\breve\CF^\an$ and assume $\bM_{b,\mu}(V)\cong\bM(0,1)\otimes_{\BQ_p}V$. We must show that also $\bM_{b,\mu}(\wt V)\cong\bM(0,1)\otimes_{\BQ_p}\wt V$. Since $V$ is a faithful representation, $\wt V$ is a direct summand of $\bigoplus_i V^{\otimes n_i}\otimes_{\BQ_p}V\dual{}^{\otimes m_i}$ for suitable positive integers $n_i,m_i$. As in Remark~\ref{RemPottharst} consider the object $(D,\phi_D,\CL):=\bigl(V_{K_0}\,,\,\rho(b)\phi\,,\,t^{-\hhh}(\id_{V_{K_0}}\otimes\theta)^{-1}Fil^\hhh D_C\bigr)\in\ul{\mathrm{MP}}^\phi_C$ associated with $\ulD_{b,\mu}(V)$ under the inverse $F_1^{-1}$ of the functor $F_1$ for Hodge-Tate weights $-\hhh$ and $-\hhh+1$. Set
\[
(\wt D,\phi_{\wt D},\wt \CL)\es:=\es\bigl(\wt V_{K_0}\,,\; \tilde\rho(b)\phi\,,\;\wt V\otimes_{\BQ_p}\bB_\dR\;\cap\;\bigoplus_i \CL^{\otimes n_i}\otimes_{\bB^+_\dR}\CL\dual{}^{\otimes m_i}\bigr)
\]
Since $F_1$ is a tensor functor it assigns to $(\wt D,\phi_{\wt D},\wt \CL)$ the filtered isocrystal $\ulD_{b,\mu}(\wt V)=(\wt D,\phi_{\wt D},Fil^\bullet_\mu \wt D_C)$ over $C$. The latter has Hodge-Tate weights $-\tilde\hhh$ and $-\tilde\hhh+1$, and hence $t^{-\tilde\hhh+1}\wt D\otimes_{K_0}\bB^+_\dR\subset\wt \CL\subset t^{-\tilde\hhh}\wt D\otimes_{K_0}\bB^+_\dR$ and the inverse $F_1^{-1}$ of $F_1$ for Hodge-Tate weights $-\tilde\hhh$ and $-\tilde\hhh+1$ assigns to $\ulD_{b,\mu}(\wt V)$ the object $(\wt D,\phi_{\wt D},\wt \CL)$ of $\ul{\mathrm{MP}}^\phi_C$. Therefore  under the tensor functor $F_2$ we have 
\[
\bM_{b,\mu}(\wt V)\;=\;F_2(\wt D,\phi_{\wt D},\wt \CL)\;=\;\wt V\otimes_{\BQ_p}\wt\bB^\dagger_\rig[t^{-1}]\;\cap\;\bigoplus_i \bM_{b,\mu}(V)^{\otimes n_i}\otimes_{\wt\bB^\dagger_\rig}\bM_{b,\mu}(V)\dual{}^{\otimes m_i}
\]
as $\phi$-modules over $\wt\bB^\dagger_\rig$. By Theorem~\ref{Thm1.4} we can write $\bM_{b,\mu}(\wt V)\cong\bigoplus_j\bM(\tilde c_j,\tilde d_j)$. From $\bM_{b,\mu}(V)\cong\bM(0,1)\otimes_{\BQ_p}V$ it follows that also $\bigoplus_i \bM_{b,\mu}(V)^{\otimes n_i}\otimes_{\wt\bB^\dagger_\rig}\bM_{b,\mu}(V)\dual{}^{\otimes m_i}$ is isomorphic to a direct sum of $\bM(0,1)$. By Proposition~\ref{Prop1.6} also its direct summand $\bM_{b,\mu}(\wt V)$ is isomorphic to $\bM(0,1)\otimes_{\BQ_p}\wt V$ as desired.
\end{proof}

\begin{proposition}\label{Prop2.7}
The set $\breve\CF_b^0$ is invariant under the group $J(\BQ_p)$ from (\ref{EqJ}) and contained in $\breve\CF^{wa}_b$ with $\breve\CF_b^0(\LL)=\breve\CF^{wa}_b(\LL)$ for any finite extension $\LL/\breve E$. 
\end{proposition}

\begin{proof}
The invariance under $J(\BQ_p)$ was remarked before Proposition~\ref{PropDOR}.

Let $\mu\in\breve\CF_b^0$ be an analytic point and set $\LL=\CH(\mu)$. Consider a faithful representation $V$ as in (\ref{Eq2.2}) and let $\ulD_\mu=\ulD_{b,\mu}(V)=(D,\phi_D,Fil^\bullet_\mu D_\LL)$. Let $D'\subset D$ be a $\phi_D$-stable $K_0$-subspace and let $Fil^i_\mu D'_\LL:=D'_\LL\cap Fil^i_\mu D_\LL$. We have to show that the inequality \mbox{$t_H(\ulD')\le t_N(\ulD')$} holds for the subobject $\ulD'=(D',\phi_D|_{D'},Fil^\bullet_\mu D'_\LL)\subset \ulD_\mu$ with equality if $\ulD'=\ulD_\mu$. Consider the $\phi$-submodule $\bM':=\bM(\ulD')\subset\bM(\ulD_\mu)=\bM_\mu$. Then 
\[
\rk\bM'\cdot\weight\bM'\;=\;\deg\bM'\;=\;t_N(\ulD')-t_H(\ulD')
\]
by Theorem~\ref{Thm1.13}. Since $\mu\in\breve\CF_b^0$ we have $\bM_\mu\cong\bM(0,1)\otimes_{\BQ_p}V$ and hence $t_H(\ulD_\mu)=t_N(\ulD_\mu)$. Moreover, $\weight\bM'\ge\weight\bM_\mu=0$ by Proposition~\ref{Prop2.9} proving $t_N(\ulD')\ge t_H(\ulD')$. This shows that $\breve\CF_b^0\subset\breve\CF^{wa}_b$. The equality on $\LL$-valued points follows from Theorem~\ref{Thm1.5}. 
\end{proof}

We expect that the equality on $\LL$-valued points holds for many more $\LL$.

\begin{conjecture}\label{Conj5.3b}
If $\LL/\breve E$ is a complete valued extension such that the value group of $\LL$ is finitely generated then $\breve\CF^0_b(\LL)=\breve\CF^{wa}_b(\LL)$.
\end{conjecture}

We will see in a moment that $\breve\CF^0_b\subset\breve\CF^{wa}_b$ is open. We also make the

\begin{conjecture}\label{Conj5.3c}
$\breve\CF^0_b$ is arcwise connected.
\end{conjecture}

We give evidence for these two conjectures in Section~\ref{SectConjectures}. In particular we show in Proposition~\ref{Prop1.2b} that Conjecture~\ref{Conj5.3c} follows from Conjecture~\ref{Conj5.3b}. It is known in some cases that Conjecture~\ref{Conj5.3c} holds; see Remarks~\ref{RemConnected1}(a) and \ref{Rem7.5}(b).

\begin{theorem}\label{Thm2.1}
The set $\breve\CF_b^0$ is a paracompact open $\breve E$-analytic subspace of $\breve\CF^\an$. If $b$ is decent with the integer $s$; cf.~\eqref{EqDecent}; then $\breve\CF_b^0$ has a natural structure of paracompact open $E_s$-analytic subspace of $(\CF\otimes_E E_s)^\an$ from which it arises by base change to $\breve E$.
\end{theorem}

\begin{proof}
(a) \es Since $\breve\CF^\an$ has a finite covering by affine spaces, all its open $\breve E$-analytic subspaces are paracompact by Lemma~\ref{LemmaParacompact}. To prove that $\breve\CF^0_b\subset\breve\CF^\an$ is open consider a faithful representation $V$ of dimension $h$ as in (\ref{Eq2.2}) and let $G'=\GL(V)$. This defines a closed embedding $\CF\hookrightarrow\CF':=\CF lag(V)\otimes_{\BQ_p}E$ into the flag variety from (\ref{Eq2.0}). Here $\CF lag(V)$ is a Grassmannian isomorphic to $G'/S'$ where $S'=\Stab_{G'}(V_0)$ is the stabilizer of an appropriate subspace $V_0$ of $V$. By definition $\breve\CF_b^0=\breve\CF^\an\cap\breve\CF'_b{}^0$. So it suffices to prove the theorem for $G'$ instead of $G$. Since $G'$ is connected we may assume by \cite{Kottwitz} that $b$ is decent, say with integer $s$. We let $\CF'_s{}^\an:=(\CF'\otimes_E E_s)^\an$ and define the subset $\CF'_b{}^0\subset\CF'_s{}^\an$  by the same condition as in (\ref{EqFAdm}). We will show that it is open. For $\mu\in\CF'_s{}^\an$ we set $\ulD_\mu=\ulD_{b,\mu}(V)$. 
By Lemma~\ref{LemmaKottwitzCond} our assumption \eqref{Eq2.2} implies  $t_N(\ulD_\mu)=t_H(\ulD_\mu)$. We choose an integer $e\ge(\dim V)-1$ which is a multiple of $s$. Since $\ulD_\mu$ has Hodge-Tate weights $-\hhh$ and $-\hhh+1$, the set $\CF_s'{}^\an\setminus\CF'_b{}^0$ is by Proposition~\ref{Prop1.8} equal to
\begin{eqnarray*} 
&&\Bigl\{\,\mu\in\CF_s'{}^\an\text{ analytic points} : \text{ there exists an algebraically closed complete}\nonumber \\[-1mm]
&&\quad \text{extension $C$ of $\CH(\mu)$ and an element $x\in\Koh^0_{\phi^e}\bigl(\bM(1,1)\otimes[e]_\ast(t^{-\hhh}\bD)\bigr)$, }x\ne0\qquad\\[-1mm]
&&\quad \text{with }\es\theta_{t^{-\hhh}\bD}\bigl(\phi_{t^{-\hhh}\bD}^m(x)\bigr)\;\in\;(Fil^\hhh_\mu V_{\CH(\mu)})\otimes_{\CH(\mu)}C\es\text{ for all }m=0,\ldots,e-1\quad\Bigr\}\,. 
\end{eqnarray*}
Here $\bD:=\bD(\ulD_\mu)=(D,\phi_D)\otimes_{K_0}\wt\bB^\dagger_\rig(C)$ is the $\phi$-module over $\wt\bB^\dagger_\rig(C)$ associated with $\ulD_\mu$ in Definition~\ref{Def1.12b} and $[e]_\ast(t^{-1}\bD)=(D,p^{-e}\phi_D^e)\otimes_{K_0}\wt\bB^\dagger_\rig(C)$.
  
\medskip

(b) \es We identify $V\otimes_{\BQ_p}\BA^1_{E_s}$ with affine $h$-space $\BA^h_{E_s}$ over $E_s$. For $\zeta\in E_s$ consider the $E_s$-analytic polydisc with radii $(|\zeta|,\ldots,|\zeta|)$
\[
\BD(\zeta)^{eh} \;=\;\CM\bigl(E_s\langle\tfrac{w_{im}}{\zeta}\,:\;i=1,\ldots,h\,;\,m=0,\ldots,e-1\rangle\bigr)\;\subset\; (\BA^h_{E_s})^e.
\]
We will construct in (c) below a constant $\zeta\in E_s$ and a compact subset $Z$ of $\BD(\zeta)^{eh}$ with the following property:
If $C$ is an algebraically closed complete extension of $E_s$ and $x\in\Koh^0_{\phi^e}\bigl(\bM(1,1)\otimes[e]_\ast(t^{-\hhh}\bD)\bigr)$ with $x\ne0$ then for some integer $N$
\[
\bigl(w_m\bigr)_{m=0}^{e-1}\es:=\es\Bigl(\bigl(w_{1m},\ldots,w_{hm}\bigr)^T\Bigr)_{m=0}^{e-1}\es:=\es\Bigl(p^N\theta_{t^{-\hhh}\bD}\bigl(\phi_{t^{-\hhh}\bD}^m(x)\bigr)\Bigr)_{m=0}^{e-1}
\]
is a $C$ valued point of $Z$ and $Z$ consists precisely of those points.

Now let $G'_s{}^\an$ be the $E_s$-analytic space associated with the group scheme $G'\otimes_{\BQ_p}E_s$ and consider the morphism of $E_s$-analytic spaces
\begin{eqnarray*}
\beta:\es G'_s{}^\an\times_{E_s}\BD(\zeta)^{eh} & \longto & (\BA^h_{E_s})^e\es\cong\es(V\otimes_{\BQ_p}\BA^1_{E_s})^e\\[2mm]
\bigl(\;g\quad,\, (w_m)_{m=0}^{e-1}\,\bigr) & \longmapsto & \bigl(g^{-1}w_m\bigr)_{m=0}^{e-1}\es.
\end{eqnarray*}
Let $Y$ be the closed subset of $G'_s{}^\an\times_{E_s}\BD(\zeta)^{eh}$ defined by the condition that $(w_m)_{m=0}^{e-1}$ belongs to $Z$ and that $\beta(Y)\subset (V_0\otimes_{\BQ_p}\BA^1_{E_s})^e$. Furthermore consider the projection map
\[
pr_1:\es G'_s{}^\an\times_{E_s}\BD(\zeta)^{eh}\es\longto\es G'_s{}^\an
\]
onto the first factor and the canonical map $\gamma:G'_s{}^\an\to\CF'_s{}^\an$ coming from the isomorphism $\CF'_s\cong G'_s/\Stab_{G'_s}(V_0)$. Then $\mu\in\CF'_s{}^\an$ does not belong to $\CF'_b{}^0$ if and only if $\mu\in\gamma\circ pr_1(Y)$. Since $\BD(\zeta)^{eh}$ is quasi-compact the projection $pr_1$ is a proper map of topological Hausdorff spaces by \cite[Proposition 3.3.2]{Berkovich1}. Thus in particular it is closed and $pr_1(Y)$ is closed. Note that $\CF'_s{}^\an$ carries the quotient topology under $\gamma$ since $\gamma$ is a smooth morphism of schemes, hence open by \cite[Proposition 3.5.8 and Corollary 3.7.4]{Berkovich2}. Since by construction $pr_1(Y)=\gamma^{-1}\bigl(\gamma\circ pr_1(Y)\bigr)$ we conclude that $\CF'_b{}^0=\CF'_s{}^\an\setminus\gamma\circ pr_1(Y)$ is open in $\CF'_s{}^\an$ as desired.

\medskip

(c) \es It remains to construct the compact set $Z$. Since $b$ is decent, the $\phi^e$-module $\bM(1,1)\otimes[e]_\ast(t^{-\hhh}\bD)$ is isomorphic to $\bigoplus_{i=1}^h \bM(-c_i,1)$ for suitable integers $c_i$. We assume that the identification of $V\otimes_{\BQ_p}\BA^1_{E_s}$ with $\BA^h_{E_s}$ in (b) was chosen compatibly with this direct sum decomposition. Let 
\[
c_1,\ldots,c_k\;>\;0\,=\,c_{k+1}\,=\,\ldots=c_\ell\;>\;c_{\ell+1},\ldots,c_h\,.
\]
Since all Newton slopes of $\bD$, respectively $t^{-\hhh}\bD$, lie in $[\hhh-1,\hhh]$, respectively $[-1,0]$, we have $e-1\ge c_i\ge-1$ for all $i$. Then by Propositions~\ref{Prop1.7} and \ref{Prop1.6}
\begin{eqnarray*}
\Koh^0_{\phi^e}\bigl(\bM(1,1)\otimes[e]_\ast(t^{-\hhh}\bD)\bigr)&\cong&
\bigoplus_{i=1}^k\Bigl\{\,\sum_{\nu\in\BZ}p^{c_i\nu}\sum_{j=0}^{c_i-1}p^j\phi^{-e\nu}([u_{ij}]):\es u_{ij}\in\wt\bE,\,v_\bE(u_{ij})>0\,\Bigr\}\\
&\oplus&\bigoplus_{i=k+1}^\ell W(\BF_{p^e}){\TS[\frac{1}{p}]}\\
&\oplus&\bigoplus_{i=\ell+1}^h (0) \quad.
\end{eqnarray*}
For $1\le i\le k, 0\le j\le c_i-1$ and all integers $\nnn\ge0$ consider the compact $E_s$-analytic spaces which we just view as compact sets 
\begin{eqnarray*}
U^{(0)}_{ij}&:=&\CM\bigl(E_s\langle{\TS\frac{u^{(0)}_{ij}}{p}}\rangle\bigr)\es=\es\bigl\{\,|u^{(0)}_{ij}|\le|p|\,\bigr\} \qquad\text{and}\\[2mm]
U^{(\nnn)}_{ij}&:=&\CM\bigl(E_s\langle u^{(\nnn)}_{ij}\rangle\bigr)\es=\es\bigl\{\,|u^{(\nnn)}_{ij}|\le 1\,\bigr\}\qquad\text{for }\nnn>0\,.
\end{eqnarray*}
Then the sets
\begin{eqnarray*}
U_{ij}&:=&\Bigl\{\,(u^{(\nnn)}_{ij})_{\nnn\in\BN_0}\in\prod_{\nnn\in\BN_0}U^{(\nnn)}_{ij}:\es(u^{(\nnn+1)}_{ij})^p=u^{(\nnn)}_{ij} \text{ for all }\nnn\ge0\,\Bigr\}\quad\text{and} \\[2mm]
U&:=& \prod_{i=1}^k\prod_{j=0}^{c_i-1}U_{ij}\times\prod_{i=k+1}^\ell W(\BF_{p^e})\times\prod_{i=\ell+1}^h\{0\}
\end{eqnarray*}
are compact by Tychonoff's theorem.
For an arbitrary algebraically closed complete extension $C$ of $E_s$ consider a $C$-valued point $u$ of $U$ given by
\[
\Bigl( \bigl((u^{(\nnn)}_{ij})_{\nnn\in\BN_0}\bigr)_{i=1,\ldots,k\,;\;j= 0,\ldots, c_i-1}\,,\,(a_i)_{i=k+1,\ldots,\ell}\,,\,(0)_{i=\ell+1,\ldots,h}\Bigr)
\]
with $u^{(\nnn)}_{ij}\in C$ and $a_i\in W(\BF_{p^e})$. We assign to $u$ the $C$-valued point $y$ of $\BA_{E_s}^{eh}$ with $(w_m)_{m=0}^{e-1}=\,\bigl(\theta_{t^{-\hhh}\bD}(\phi_{t^{-\hhh}\bD}^m(x))\bigr)_{m=0}^{e-1}$ where $x$ is the element of $\Koh^0_{\phi^e}\bigl(\bM(1,1)\otimes[e]_\ast(t^{-\hhh}\bD)\bigr)$ associated with the $u_{ij}^{(\nnn)}$ and $a_i$. In concrete terms this means
\begin{eqnarray*}
w_m
&:=& b\cdot\phi(b)\cdots\phi^{m-1}(b)\cdot\left(\begin{array}{c}
\DS\biggl(\;\sum_{j=0}^{c_i-1}p^j \Bigl(\sum_{\nu>0} p^{c_i\nu}u^{(e\nu-m)}_{ij}\;+\;\sum_{\nu\le 0}p^{c_i\nu}\bigl(u^{(0)}_{ij}\bigr)^{p^{m-\nu e}}\Bigr)\biggr)_{i=1}^k \\[5mm]
\DS\bigl(\phi^m(a_i)\bigr)_{i=k+1}^\ell\\[2mm]
\DS(0)_{i=\ell+1}^h
\end{array}\right)
\end{eqnarray*}
for $m=0,\ldots,e-1$. This defines a map $\alpha:U\to\BA_{E_s}^{eh},u\mapsto y$ of topological Hausdorff spaces. We claim that $\alpha$ is continuous. By definition of the topology on $\BA_{E_s}^{eh}$ (see Appendix~\ref{AppBerkovichSpaces}) the map $\alpha$ is continuous if and only if for any polynomial $f\in E_s[w_{im}:1\le i\le h,\,0\le m\le e-1]$ and any open interval $I\subset\BR$ the preimage under $\alpha$ of the open set
\[
W\es=\es\bigl\{\,y\in\BA_{E_s}^{eh}:\es|f|_y\in I\,\bigr\}
\]
is open in $U$. Since $|u^{(0)}_{ij}|_u\le|p|$ there is a constant $\zeta\in E_s$ such that $|w_{im}|_{\alpha(u)}\le|\zeta|$ for all $i$ and $m$ and all $u\in U$. Write $f=\sum_{\ul n}b_{\ul n}\,\ul w^{\ul n}$ where $b_{\ul n}\in E_s$ for every multi index $\ul n\in\BN_0^{eh}$. If we set $w_{im}=w'_{im}+w''_{im}$ we obtain from the Taylor expansion of $f$ in powers of $\ul w''$ a bound $\delta$ such that for all $y\in\BA^{eh}_{E_s}$ the condition $|w''_{im}|_y\le\delta$ for all $i,m$ implies
\[
|f(\ul w)|_y \in I \text{ and } |\ul w|_y\le|\zeta|\qquad\Longleftrightarrow\qquad|f(\ul w')|_y \in I \text{ and } |\ul w'|_y\le|\zeta|\,.
\]
Now we fix a positive integer $r$ and set for all $m=0,\ldots,e-1$
\begin{eqnarray*}
\bigl(w''_{im}\bigr)_{i=1}^h &:=&
b\cdot\phi(b)\cdots\phi^{m-1}(b)\cdot
\left(\begin{array}{c}
\DS\Bigl(\sum_{j=0}^{c_i-1}p^j \sum_{\nu> r} p^{c_i\nu}u^{(e\nu-m)}_{ij}\Bigr)_{i=1}^k \\[5mm]
\DS(0)_{i=k+1}^\ell\\[2mm]
\DS(0)_{i=\ell+1}^h
\end{array}\right) \qquad\text{and} \\[8mm]
\bigl(w'_{im}\bigr)_{i=1}^h &:=&
b\cdot\phi(b)\cdots\phi^{m-1}(b)\cdot
\left(\begin{array}{c}
\DS\Bigl(\sum_{j=0}^{c_i-1}p^j \sum_{\nu\le r}p^{c_i\nu}\bigl(u^{(er)}_{ij}\bigr)^{p^{m+(r-\nu) e}}\Bigr)_{i=1}^k \\[5mm]
\DS\bigl(\phi^m(a_i)\bigr)_{i=k+1}^\ell\\[2mm]
\DS(0)_{i=\ell+1}^h
\end{array}\right)\es.
\end{eqnarray*}
Since $|u^{(\nnn)}_{ij}|_u<1$ for all $\nnn,i,j$ and for all points $u\in U$ we can find a large enough integer $r$ such that $|w''_{im}|_{\alpha(u)}\le\delta$ for all $i,m$ and for all $u\in U$. This shows that
\[
\alpha^{-1}(W)\es=\es W'\es:=\es\bigl\{\,u\in U:\es|f(w'_{im})|_{\alpha(u)}\in I\,\bigr\}\,.
\]
To prove that $W'$ is open in $U$ observe that the map
\begin{eqnarray*}
\alpha'_m:\es U'\es:=\es\prod_{i=1}^k\prod_{j=0}^{c_i-1}U^{(er)}_{ij}\times\prod_{i=k+1}^\ell \BA^1_{E_s}&\longto& \BA_{E_s}^{h} \\[2mm]
\Bigl(\;(u^{(er)}_{ij})_{i,j}\;,\es\bigl(\phi^m(a_i)\bigr)_i\;\Bigr)&\longmapsto& (w'_{im})_{i=1}^h
\end{eqnarray*}
is a morphism of $E_s$-analytic spaces, and in particular continuous. Furthermore, the projection maps $\alpha_{ij}:U_{ij}\to U^{(er)}_{ij}$ and the inclusion 
\[
\alpha_i:\es W(\BF_{p^s})\es\longto\es\BA^e_{E_s}\es=\es\CM\bigl(E_s[T_0,\ldots,T_{e-1}]\bigr)\,,\quad a\es\longmapsto\es \bigl(T_m\mapsto \phi^m(a)\bigr)
\]
for $i=k+1,\ldots,\ell$ are continuous. Since 
\[
\alpha'\es:=\es(\alpha'_0\times\ldots\times\alpha'_{e-1})\circ\Bigl(\,\prod_{i=1}^k\prod_{j=0}^{c_i-1}\alpha_{ij}\times\prod_{i=k+1}^\ell \alpha_i\Bigr):\es U\es\longto\es\BA^{eh}_{E_s} 
\]
maps $u$ to $\alpha'(u)=(w'_m)_{m=0}^{e-1}$
the set $W'=(\alpha')^{-1}(W)$ is open in $U$ and this proves that $\alpha$ is continuous. 
 
Now multiplying $(w_{im})_{i,m}$ with $p$ amounts to replacing $u^{(\nnn)}_{ij}$ by $u^{(\nnn)}_{i,j-1}$ for $j=1,\ldots,c_i-1$, and $u^{(\nnn)}_{i,0}$ by $(u^{(\nnn)}_{i,c_i-1})^{p^e}$, and $a_i$ by $pa_i$. Thus we may take
\[
Z\es:=\es\alpha\Bigl(U\setminus\bigl\{\,u\in U:\es|u^{(0)}_{ij}|_u<|p|^{p^e}\,,\,a_i\in pW(\BF_{p^e})\es\text{for all $i$ and $j$}\,\bigr\}\Bigr)\,.
\]
Then $Z$ is the continuous image of a compact set and satisfies the property required in (b). This proves the theorem.
\end{proof}

We give an example showing that the inclusion $\breve\CF_b^0\subset\breve\CF^{wa}_b$ may be strict. Similar examples were independently obtained by A.~Genestier and V.~Lafforgue. In Section~\ref{SectWAImpliesAdm} we will determine all $b\in\GL_n(K_0)$ for which $\breve\CF_b^0=\breve\CF_b^{wa}$.

\begin{example}\label{Example6.4}
Let $G=\GL_5$ and consider the conjugacy class $\{\mu\}$ of cocharacters containing $\mu:\BG_m\to G\,,\,z\mapsto\diag(z,z,1,1,1)$. We have $E=\BQ_p$ and $\CF\cong\Grass(2,5)$ the Grassmannian. Let
\begin{equation}\label{EqExample6.4b}
b\es=\es \left(
\begin{array}{ccccc@{\es}} \,0\, & \,0\, & \,0\, & \,0\, & \!\!p^{-2}\!\! \\ p & 0 & 0 & 0 & 0 \\ 0 & p & 0 & 0 & 0 \\ 0 & 0 & p & 0 & 0 \\0 & 0 & 0 & p & 0 \end{array}\right)
\es\in\es G(K_0)\,.
\end{equation}
The element $b$ is decent with integer $s=5$ and $(G,b,\{\mu\})$ satisfies \eqref{Eq2.2} for $\hhh=1$. Since the isocrystal \mbox{$({K_0}^{\oplus 5},b\!\cdot\!\phi)$} is simple every cocharacter in $\{\mu\}$ is weakly admissible, that is $\breve\CF^{wa}_b=\breve\CF^\an$.

Let $\mu\in\breve\CF^\an$ and let $C$ be the completion of an algebraic closure of $\CH(\mu)$. The $\phi$-module $\bM_\mu=\bM({K_0}^{\oplus 5},b\!\cdot\!\phi,Fil^\bullet_\mu {\CH(\mu)}^{\oplus 5})$ constructed from $b$ and $\mu$ satisfies $\deg\bM_\mu=0$ and 
\[
t^{-1}\bD=\es\bD\otimes\bM(-1,1)\es=\es\bM(-3,5)\es\supset\es\bM_\mu\es\supset\es\bD\es\es\cong\es\bM(2,5)
\]
where $\bD=\bD({K_0}^{\oplus 5},b\!\cdot\!\phi)=({K_0}^{\oplus 5},b\!\cdot\!\phi)\otimes_{K_0}\wt\bB^\dagger_\rig(C)$ and $t^{-1}\bD=\bD({K_0}^{\oplus 5},p^{-1}b\!\cdot\!\phi)$.
Since by Proposition~\ref{Prop1.6} the weight of every summand of $\bM_\mu$ lies between $-3/5$ and $2/5$, either $\bM_\mu\cong\bM(0,1)^{\oplus 5}$ or $\bM_\mu\cong\bM(-1,2)\oplus\bM(1,3)$. (The first entries must sum to $0=\deg\bM_\mu$ and the second entries must sum to $5=\rk\bM_\mu$.) Now the bad situation $\bM_\mu\cong\bM(-1,2)\oplus\bM(1,3)$ occurs if and only if
\begin{eqnarray*}
(0)&\ne&\Hom_{\phi}\bigl(\bM(-1,2)\,,\,\bM_\mu\bigr)\es=\\[2mm]
&=&\Bigl\{\,f\in\Hom_\phi\bigl(\bM(-1,2)\,,\,t^{-1}\bD\bigr)\colon\;\theta_{t^{-1}\bD}(\im f)\subset Fil^1_\mu {\CH(\mu)}^{\oplus 5}\,\Bigr\}\,.
\end{eqnarray*}
We compute $\Hom_\phi\bigl(\bM(-1,2)\,,\,\bM(-3,5)\bigr)=$
\[
= \es \Bigl\{\,A\;=\;(a_{i,j})_{i=1\ldots5,\,j=1,2}\in\wt\bB^\dagger_\rig(C)^{5\times2}\colon p^{-1}b\cdot\phi(A)=A\cdot\left(  \begin{array}{cc} 0 & p^{-1} \\ 1 & 0 \end{array}\right)\Bigr\} 
\]
with $b$ as in \eqref{EqExample6.4b}. This implies $p^{-3}\phi(a_{5,1})=a_{1,2},\;p^{-3}\phi(a_{5,2})=p^{-1}a_{1,1}$ and $\phi(a_{i-1,1})=a_{i,2},\;\phi(a_{i-1,2})=p^{-1}a_{i,1}$ for $2\le i\le5$. So $a_{1,2}$ must satisfy $p^{-1}\cdot\phi^{10}(a_{1,2})=a_{1,2}$, that is, $x:=a_{1,2}\in\Koh^0_{\phi^{10}}\bigl(\bM(-1,1)\bigr)$. By Proposition~\ref{Prop1.7} it follows that $x=\sum_{\nu\in\BZ}p^\nu\phi^{-10\nu}([u])$ for some $u\in\wt\bE(C)$ with $v_\bE(u)>0$. We conclude that $\Hom_\phi\bigl(\bM(-1,2)\,,\,\bM(-3,5)\bigr)=$
\[
= \es \Biggl\{\,A\;=\;\left(  \begin{array}{cc}
\phi^5(x) & x \\ \phi^{11}(x) & \phi^6(x) \\ \phi^{17}(x) & \phi^{12}(x) \\ \phi^{23}(x) & \phi^{18}(x) \\ \phi^{29}(x) & \phi^{24}(x)   \end{array}\right)
\begin{array}{l}
\\[2mm]
\DS:\es x\;=\;\sum_{\nu\in\BZ}p^\nu\phi^{-10\nu}([u])\;\in\;\Koh^0_{\phi^{10}}\bigl(\bM(-1,1)\bigr)\,, \\[2mm]
\qquad\quad\qquad u\;\in\;\wt\bE(C),\,v_\bE(u)>0
\end{array}
\Biggr\}\,.
\]
Thus the points $\mu\in\breve\CF^\an$ for which the columns of $\theta_{t^{-1}\bD}(A)$ are contained in $Fil^1_\mu {\CH(\mu)}^{\oplus 5}$ do not belong to $\breve\CF_b^0$. We may therefore take any nonzero $u\in\wt\bE(C)$ with $v_\bE(u)>0$, define the matrix $A$ correspondingly, and let $\mu\in\breve\CF^\an(C)$ be a point for which $Fil^1_\mu {\CH(\mu)}^{\oplus 5}$ contains the columns of $\theta_{t^{-1}\bD}(A)$. All these points lie in the complement of $\breve\CF_b^0$ in $\breve\CF^{wa}_b$.

Note that $\rk\theta_{t^{-1}\bD}(A)=2$ and $\mu$ is uniquely determined by $A$ and hence by $u$. Indeed, assume the contrary and consider $\wedge^2A\in\Hom_\phi\bigl(\wedge^2\bM(-1,2),\wedge^2\bM(-3,5)\bigr)$ which is non-zero because $A$ is an inclusion $\bM(-1,2)\into\bM(-3,5)$. Here $\wedge^2\bM(-1,2)=\bM(-1,1)$ and $\wedge^2\bM(-3,5)\cong\bM(-6,5)^{\oplus2}$ by \cite[Proposition 4.6.3]{Kedlaya}. Then $\theta_{\wedge^2t^{-1}\bD}(\wedge^2A)=\wedge^2\theta_{t^{-1}\bD}(A)=0$ implies that $\wedge^2A=T\circ B$ for $0\ne B\in\Hom_\phi\bigl(\bM(-1,1)\,,\,\bM(-1,5)^{\oplus2}\bigr)$ by Corollary~\ref{Cor1.6b}. But since $-1<-\frac{1}{5}$, the later space is zero by Proposition~\ref{Prop1.6} which gives the contradiction to our assumption.

In particular the map $\alpha:\{\,u\in\wt\bE:0<|u|<1\,\}\to\breve\CF^\an,\,u\mapsto\mu$ parametrizes the complement $\breve\CF^\an\setminus\breve\CF^0_b$. The nature of the map $\alpha$ is somewhat mysterious. It is continuous by similar arguments as in part (c) of the proof of Theorem~\ref{Thm2.1} but its image does not meet the dense subset of points $x\in\breve\CF^\an$ with $\CH(x)$ finite over $K_0$.
\end{example}

\medskip

We conclude this section by stating the following theorem which we will prove at the end of Section~\ref{SectPeriodMorph}.

\begin{theorem}\label{ThmDOR}
For any representation $V$ of $G$ as in (\ref{Eq2.2}) there is a local system $\CV$ of $\BQ_p$-vector spaces on $\breve\CF^0_b$ such that the fiber of $\CV$ at any point $\mu\in\breve\CF^0_b(\LL)$ for a finite extension $\LL/\breve E$ induces via \eqref{EqGalRep} the crystalline Galois representation $V_\cris\bigl(\ulD_{b,\mu}(V)\bigr)$ associated with the (weakly) admissible filtered isocrystal $\ulD_{b,\mu}(V)=(V_{K_0},\rho(b)\!\cdot\!\phi,Fil^\bullet_\mu V_\LL)$.
\end{theorem}

\begin{remark}
Note that this theorem is weaker than Conjecture~\ref{ConjRZ} in two respects. Firstly we do not know in general whether $V\mapsto\CV$ is a tensor functor $\Rep_{\BQ_p}G\to\PLoc_{\breve\CF^0_b}$ as required in Conjecture~\ref{ConjRZ}. However this is true if $G=\GL_n$ (see Theorem~\ref{ThmC}) or if $G$ is the group associated to an EL or PEL situation (see Theorem~\ref{ThmPEL}).

Secondly we do neither know whether $\breve\CF^0_b$ is connected nor whether it is the largest open subspace of $\breve\CF^{wa}_b$ with this property. However, we conjecture that this is true.
\end{remark}

\begin{conjecture}\label{Conj5.10}
$\breve\CF^0_b$ is the largest open $\breve E$-analytic subspace of $\breve\CF^{wa}_b$ carrying a local system $\CV$ that has the properties stated in the theorem.
\end{conjecture}

 We give some evidence for this in Section~\ref{SectConjectures}.

%
%

\section{Relation with Period Morphisms} \label{SectPeriodMorph}
\setcounter{equation}{0}

For $G=\GL(V)$, Rapoport and Zink study a period morphism $\breve\pi:\breve\CG^\an\to\breve\CF_b^{wa}$. We keep the notation of Section~\ref{SectOpenness}. In particular let $\tilde b\in\GL(V)(K_0)$, let $\{\tilde\mu\}$ be a conjugacy class of cocharacters $\BG_m\to\GL(V)$. Then $\breve E=K_0$. We assume that $V$ satisfies (\ref{Eq2.2}) for an integer $\tilde\hhh$. We set $b:=p^{1-\tilde\hhh}\!\cdot\!\tilde b\in G(K_0)$ and define $\mu:\BG_m\to G$ by $\mu(z)=z^{1-\tilde n}\!\cdot\!\tilde\mu(z)$. Then $\ulD_{b,\mu}(V)=\ulD_{\tilde b,\tilde\mu}(V)\otimes\BOne(1-\tilde\hhh)$ and $\bM_{b,\mu}(V)=\bM_{\tilde b,\tilde\mu}(V)$ by Lemma~\ref{Lemma4.2'}. This implies that the assignment $\tilde\mu\mapsto\mu$ is an isomorphism $\breve\CF_{\tilde b}^0\isoto\breve\CF_b^0$. 

Now $(G,b,\{\mu\})$ and $V$ satisfy \eqref{Eq2.2} for $\hhh=1$. Therefore $(V_{K_0},b\!\cdot\!\phi)$ is the covariant Dieudonn\'e-module $(D,\phi_D)$ of a $p$-divisible group $\bX=\bX_b(V)$ over $\BF_p^{\,\alg}$; see Messing~\cite{Messing}. This means $D:=\BD(\bX)$ is the tensor product with $K_0$ of the Lie algebra of the universal vector extension of some (any) lift of $\bX$ to a $p$-divisible group over $W(\BF_p^{\,\alg})$. On $\bX$ there is the Verschiebung morphism $Ver_\bX:\Frob_p^\ast\bX\to\bX$. By the crystalline nature of $\BD(\bX)$ it induces the $\phi$-semilinear isomorphism $\phi_D:=\BD(Ver_\bX)\circ\phi:D\to D$. Moreover, \eqref{Eq2.2} and Lemma~\ref{LemmaKottwitzCond} imply 
\[
\dim\bX=\dim_{K_0}D-t_N\bigl(\ulD_{b,\mu}(V)\bigr)=\dim_{K_0}D-t_H\bigl(\ulD_{b,\mu}(V)\bigr)=\dim V-\dim_\LL Fil^1 D_\LL=:d
\]
and $\breve\CF^\an$ is the Grassmannian of $(\dim V-d)$-dimensional subspaces of $V$. 

The period morphism is constructed as follows. We denote by $\Nilp_{W(\BF_p^{\,\alg})}$ the category of $W(\BF_p^{\,\alg})$-schemes on which $p$ is locally nilpotent. For $S\in\Nilp_{W(\BF_p^{\,\alg})}$ denote by $\bar S$ the closed subscheme defined by the ideal $p\CO_S$. We set $\bX_{\bar S}:=\bX\times_{\BF_p^\alg}\bar S$.

\begin{proposition}[{\cite[Theorem 2.16]{RZ}}]\label{Prop6.1}
The contravariant functor $\Nilp_{W(\BF_p^{\,\alg})} \longto \Sets$
\begin{eqnarray*}
S & \longmapsto & \Bigl\{\,\text{Isomorphism classes of }(X,\eta:\bX_{\bar S}\to X_{\bar S}) \text{ where} \\
& & \qquad X \text{ is a $p$-divisible group over $S$ and}\\
& & \qquad \eta \text{ is a quasi-isogeny over $\bar S$ to $X_{\bar S}=X\times_S \bar S$}\quad\Bigr\}
\end{eqnarray*}
is representable by a formal scheme $\breve\CG$ locally formally of finite type over $\Spf W(\BF_p^{\,\alg})$.
\end{proposition}

The group $J(\BQ_p)$ from \eqref{EqJ} equals the quasi-isogeny group of $\bX$. It acts on $\breve\CG$ by $\gamma\colon(X,\eta)\mapsto(X,\eta\circ\gamma^{-1})$ for $\gamma\in J(\BQ_p)$.

Rapoport and Zink also study formal moduli schemes corresponding to additional data on $(X,\eta)$ of type EL and PEL; see Section~\ref{SectPELPeriodM}. Their period morphisms are derived from the following prototype.

Let $(X,\eta:\bX_{\ol{\breve\CG}}\to X_{\ol{\breve\CG}})$ be the universal $p$-divisible group with quasi-isogeny over $\breve\CG$. Let $\BD(X)_{\breve\CG}$ be the Lie algebra of the universal vector extension of $X$ over $\breve\CG$. It sits in the exact sequence
\[
0\es\longto\es(\Lie X\dual)\dual \es\longto\es\BD(X)_{\breve\CG}\es\longto\es\Lie X\es\longto\es0\,.
\]
Let $\breve\CG^\an$ be the $K_0$-analytic space associated with $\breve\CG$; see Theorem~\ref{ThmFormalRigBerkovich}. The quasi-isogeny $\eta$ induces an isomorphism 
\[
\BD(\eta):\es D\otimes_{K_0}\CO_{\breve\CG^\an}\es=\es\BD(\bX)_{K_0}\otimes_{K_0}\CO_{\breve\CG^\an}\es\isoto\es\BD(X)_{\breve\CG^\an}
\]
see \cite[Proposition 5.15]{RZ}. The kernel of the morphism 
\[
V_{K_0}\otimes_{K_0}\CO_{\breve\CG^\an}\es\isoto\es\BD(X)_{\breve\CG^\an}\es\longto\es(\Lie X)_{\breve\CG^\an}
\]
defines the $\breve\CG^\an$-valued point $\BD(\eta)^{-1}(\Lie X\dual)\dual_{\breve\CG^\an}$ of the Grassmannian $\breve\CF^\an$. This is the desired \emph{period morphism} $\breve\pi:\breve\CG^\an\to\breve\CF^\an$. It is equivariant with respect to the $J(\BQ_p)$-action on $\breve\CF^\an$ under which $\gamma\in J(\BQ_p)\subset G(K_0)$ sends the universal subspace $Fil^1$ on $\breve\CF^\an$ to $\gamma(Fil^1)$. Rapoport and Zink show that $\breve\pi$ factors through $\breve\CF_b^{wa}$. They noticed that $\breve\pi:\breve\CG^\an\to\breve\CF_b^{wa}$ is in general not quasi-compact \cite[5.53]{RZ}. We will give an explanation for this fact in Remark~\ref{Rem7.5}(d) below.

\begin{proposition}\label{Prop2.4}
The period morphism factors through $\breve\CF_b^0$ and is surjective on the level of rigid analytic points.
\end{proposition}

\begin{proof}
We fix the tautological representation $V$ of $G$, set $D=V_{K_0}$, and consider $\breve\CF^\an$ as the Grassmannian of $(\dim V-d)$-dimensional subspaces of $V$.
Let $x\in\breve\CG^\an$ be an analytic point and let $\mu=\breve\pi(x)\in\breve\CF^\an$. Let $\LL=\CH(x)$ and let $C$ be the completion of an algebraic closure of $\LL$. Let $X_x$ be the fiber of the universal $p$-divisible group $X$ at $x$ and consider the Tate module $T_p X_x$ of $X_x$. An element $\lambda\in T_pX_x$ corresponds to a morphism of $p$-divisible groups $\lambda:\BQ_p/\BZ_p\to X_{\CO_C}$ over $\CO_C$. By functoriality of the universal vector extension this yields the following diagram of $C$-vector spaces
\begin{equation}\label{Eq2.5}
\xymatrix  {
0 \ar[r] & \BD(\BQ_p/\BZ_p)_C \ar@{=}[r] \ar[d] & \BD(\BQ_p/\BZ_p)_C \ar[r]\ar[d]^{\BD(\lambda)_C} & 0 \\
0 \ar[r] & (Fil_\mu^1D_\LL)\otimes_\LL C \ar[r] & \BD(X)_{\CO_\LL}\otimes_{\CO_\LL} C \ar[r] & (\Lie X_x)_C \ar[r] & 0\,.
}
\end{equation}
Note that $\Lie\BQ_p/\BZ_p=(0)$ and that $Fil_\mu^1D_\LL$ is the $\LL$-subspace of $D_\LL$ associated with $X_x$ via the period morphism.
We also obtain a morphism of crystals $\BD(\lambda):\BD(\BQ_p/\BZ_p)\to\BD(X_{\CO_C})$ which we evaluate on the pd-thickening $\bB^+_\cris(C)$ of $\CO_C$. Since $\BD(\BQ_p/\BZ_p)_{\bB^+_\cris(C)}\;=\;\bB^+_\cris(C)$ (because the universal vector extension of $\BQ_p/\BZ_p$ over $\bB^+_\cris(C)$ is obtained from the sequence $0\to\BZ_p\to\BQ_p\to\BQ_p/\BZ_p\to0$ by pushout via $\BZ_p\to\bB^+_\cris(C)$) we obtain a morphism
\[
\xymatrix @R=.7pc {
T_pX_x\otimes_{\BZ_p}\bB^+_\cris(C) \es \ar[r] & \es\BD(X_x)_{\bB^+_\cris(C)} \es & \es D\otimes_{K_0} \bB^+_\cris(C)\ar[l]_{\es\sim} \\
\lambda\es\otimes\es a\qquad  \ar@{|->}[r] & \es\BD(\lambda)(a)\,,\;
}
\]
which is compatible with the Frobenius on both sides. Here the isomorphism on the right arises from the quasi-isogeny $\eta_x$ by the same reasoning as above. By Faltings's~\cite[Theorem 7]{Faltings2} the morphism on the left is injective. Since the elements of $T_pX_x$ are $\phi$-invariant inside $T_pX_x\otimes_{\BZ_p}\bB^+_\cris(C)$ and $\wt\bB^+_\rig(C)$ equals $\bigcap_{n\in\BN_0}\phi^n \bB^+_\cris(C)$ we get a monomorphism
\[
\xymatrix{
T_pX_x\otimes_{\BZ_p}\wt\bB^+_\rig(C) \es\ar@{^{ (}->}[r] & \es D\otimes_{K_0} \wt\bB^+_\rig(C)\,.
}
\]
It gives rise to a monomorphism
\[
\xymatrix @R=.7pc {
T_pX_x \es\ar@{^{ (}->}[r] & \es T_pX_x\otimes_{\BZ_p}\wt\bB^{]0,1]}(C) \es\ar@{^{ (}->}[r] & \es D\otimes_{K_0} \wt\bB^{]0,1]}(C) \\
\lambda \quad\ar@{|->}[r] & \quad\lambda\es\otimes\es 1\quad\qquad \ar@{|->}[r] & \quad\BD(\lambda)(1)
}
\]
since $\wt\bB^{]0,1]}(C)$ is a flat $\wt\bB^+_\rig(C)$-algebra. Consider the morphism $\theta:D\otimes_{K_0}\wt\bB^{]0,1]}\to D\otimes_{K_0}C$. Diagram (\ref{Eq2.5}) shows that $\theta(T_pX_x)\subset(Fil_\mu^1D_\LL)\otimes_\LL C$ and so by Proposition~\ref{Prop1.12}, the $\phi$-module $T_pX_x\otimes_{\BZ_p}\wt\bB^\dagger_\rig(C)=\bM(0,1)\otimes_{\BZ_p}T_pX_x$ is in fact contained in the $\phi$-module $\bM_\mu:=\bM(D,\phi_D,Fil_\mu^1D_\LL)$. But then  $\bM_\mu\cong\bM(0,1)\otimes_{\BQ_p}V$ by \cite[Lemma 3.4.2]{Kedlaya} because both $\phi$-modules have the same rank and $\deg\bM_\mu=0$. This proves that $\breve\pi$ factors through $\breve\CF_b^0$.

It remains to show that
$\breve\CG^\an\to\breve\CF_b^0\hookrightarrow\breve\CF^{wa}_b$ is surjective on the level of rigid analytic points. This follows from the Colmez-Fontaine Theorem~\cite{CF} and Kisin's proof~\cite[Theorem 0.3]{Kisin} of Fontaine's conjecture that every crystalline representation with Hodge-Tate weights $0$ and $-1$ arises from a $p$-divisible group. See also Breuil~\cite[Theorem 1.4]{Breuil00}.
\end{proof}

Using a result of Faltings~\cite{Faltings08}, or alternatively a result of Scholze and Weinstein~\cite{ScholzeWeinstein}, we show in Theorem~\ref{ThmC} below that $\breve\CF_b^0$ actually equals the image of the period morphism $\breve\pi$. In order to formulate part (c) of the theorem we recall from Remark~\ref{Rem1.8}(a) that with each tensor functor $\ul\CV:\Rep_{\BQ_p}\!\!\!G\to\PLoc_{\breve\CF^0_b}$ one associates a tower of \'etale covering spaces $\{\HeckeTower_\wt{K}\}_{\wt{K}\subset\wt G(\BQ_p)}$ of $\breve\CF^0_b$. In the situation we treat in this section, $G$ is the algebraic group $\GL(V)$ over $\BQ_p$ and in the notation introduced before Corollary~\ref{Cor1.7} we use Hilbert 90 (see \cite[Proposition~III.4.9]{Milne}) to identify $\wt\omega=\omega_0$, $\wt V=\wt\omega(V)=\omega_0(V)=V$, and $\wt G=G=\GL(\wt V)$. The situation will be more general in the next section. Similarly to the covering spaces $\HeckeTower_{\wt K}$, Rapoport and Zink \cite[5.32--5.39]{RZ} constructed a tower of coverings $\breve\CG^\an_\wt{K}$ of $\breve\CG^\an$: For any compact open subgroup $\wt{K}\subset \wt G(\BQ_p)$ let $\breve\CG^\an_\wt{K}$ be the $\breve E$-analytic space representing $\wt{K}$-level structures on the universal $p$-divisible group $X_{\breve\CG^\an}$, that is classes modulo $\wt{K}$ of trivializations $\alpha:\wt V\isoto(V_p X_{\breve\CG^\an})_{\bar\mu}$ of the rational Tate module (at some geometric base point $\bar\mu$ of $\breve\CF^0_b$), such that $\alpha\mod \wt{K}$ is invariant under the \'etale fundamental group of $\breve\CG^\an_\wt{K}$. The $\breve\CG^\an_\wt{K}$ do not depend on the choice of $\bar\mu$. After fixing a $\BZ_p$-lattice $\wt\Lambda$ in $\wt V$ and a preimage $\bar x\in\breve\CG^\an$ of $\bar\mu$ we can identify $\breve\CG^\an_{\GL(\wt\Lambda)}$ with $\breve\CG^\an$ by considering on $\breve\CG^\an$ the unique residue class of trivializations $(\alpha:\wt\Lambda\isoto (T_p X_{\breve\CG^\an})_{\bar x})\mod \GL(\wt\Lambda)$, where $\alpha$ is any isomorphism. For $\wt{K}\subset\GL(\wt\Lambda)$ the space $\breve\CG^\an_\wt{K}$ is then a finite \'etale covering space of $\breve\CG^\an$. For any $\wt{K}$ and any $g\in \wt G(\BQ_p)$ there are $\breve\CF^0_b$-isomorphisms $i_\wt{K}(g):\breve\CG^\an_\wt{K}\isoto\breve\CG^\an_{g^{-1}\wt{K}g}$ given by $(\alpha\mod \wt{K})\mapsto(\alpha\,g\mod g^{-1}\wt{K}g)$. This proves that the $\breve\CG^\an_\wt{K}$ are $\breve E$-analytic spaces and Rapoport and Zink \cite[5.39]{RZ} show that they are independent of the choice of $\wt\Lambda$ and $\bar x$. On the tower $(\breve\CG^\an_\wt{K})_{\wt{K}\subset \wt G(\BQ_p)}$ the group $\wt G(\BQ_p)$ acts through Hecke correspondences; see \cite[5.34]{RZ}. Also the action of $J(\BQ_p)$ on $\breve\CG^\an$ lifts to an action on each $\breve\CG^\an_{\wt K}$ by $\gamma\colon(X,\eta,\alpha\wt K)\mapsto (X,\eta\circ\gamma^{-1},\alpha\wt K)$ for $\gamma\in J(\BQ_p)$, which commutes with the Hecke-action. Rapoport and Zink conjectured that the $\breve\CG^\an_{\wt K}$ are \'etale covering spaces of the image $\breve\pi(\breve\CG^\an)$; compare \cite[1.37]{RZ}. We prove this in part \ref{ThmC_3} of the following

\begin{theorem}\label{ThmC}
Let $G$ be the algebraic group $\GL(V)$ over $\BQ_p$ and set $\wt G=G$. 
\begin{enumerate}
\item \label{ThmC_1}
The image $\breve\pi(\breve\CG^\an)$ of the period morphism $\breve\pi$ for $\GL(V)$ equals $\breve\CF_b^0$.
\item  \label{ThmC_2}
The rational Tate module $V_p X_{\breve\CG^\an}$ of the universal $p$-divisible group $X_{\breve\CG^\an}$ over $\breve\CG^\an$ descends to a local system $\CV$ of $\BQ_p$-vector spaces on $\breve\CF_b^0$. This induces a tensor functor $\ul\CV$ from $\Rep_{\BQ_p}\!\!\!G$ to the category $\PLoc_{\breve\CF_b^0}$ of local systems of $\BQ_p$-vector spaces on $\breve\CF_b^0$ satisfying (\ref{EqRZ}) of Conjecture~\ref{ConjRZ} and $\ul\CV(V)\cong\CV$.
\item \label{ThmC_3}
The tower of $\breve E$-analytic spaces $(\breve\CG^\an_\wt{K})_{\wt{K}\subset \wt G(\BQ_p)}$ is canonically isomorphic over $\breve\CF^0_b$ in a Hecke equivariant way to the tower of \'etale covering spaces $(\HeckeTower_\wt{K})_{\wt{K}\subset \wt G(\BQ_p)}$ of $\breve\CF^0_b$ associated with the tensor functor $\ul\CV$. In particular, $\breve\CG^\an$ is (non-canonically) isomorphic to the space of $\BZ_p$-lattices inside the local system $\CV$.
\item \label{ThmC_4}
The tensor functor $\ul\CV$ carries a canonical $J(\BQ_p)$-linearization which by Remark~\ref{Rem1.8}(b) induces an action of $J(\BQ_p)$ on the tower $(\HeckeTower_\wt{K})_{\wt{K}\subset \wt G(\BQ_p)}$. The isomorphism from \ref{ThmC_3} is equivariant for the $J(\BQ_p)$-action on both towers.
\end{enumerate}
\end{theorem}

\noindent
A previously known example for this theorem is the Lubin-Tate situation where $\breve\CF_b^0=\breve\CF^{wa}_b=\breve\CF^\an=\BP^{h-1}_{K_0}$. In this case the theorem was proved by de Jong~\cite[Proposition 7.2]{dJ} .

\begin{proof}[Proof of Theorem~\ref{ThmC}]
\ref{ThmC_1}  The inclusion $\breve\pi(\breve\CG^\an)\subset\breve\CF_b^0$ was proved in Proposition~\ref{Prop2.4}. So let now $\mu\in\breve\CF_b^0$ and let $C$ be the completion of an algebraic closure of $\LL:=\CH(\mu)$. Consider the twisted dual 
\[
\ulTD\;:=\;\bigl((D,\phi_D,Fil^\bullet_\mu D_\LL)\otimes\BOne(-1)\bigr)\dual\;=\;\CHom\bigl((D,\phi_D,Fil^\bullet_\mu D_\LL),\BOne(1)\bigr).
\]
It satisfies $Fil^0\ulTD_C=\ulTD_C$ and $Fil^2\ulTD_C=(0)$. We abbreviate $\wt\bD:=\bD(\ulTD)$ and $\wt\bM:=\bM(\ulTD)$. Then $\wt\bM\cong\bM(0,1)\otimes_{\BQ_p}V^{\SSC\lor}$ by Lemma~\ref{Lemma4.2'} and Proposition~\ref{Prop4.10}. Consider the  morphisms 
\[
\alpha:(\wt\bB^\dagger_\rig)\otimes_{\BQ_p}V^{\SSC\lor}=\bM(0,1)\otimes_{\BQ_p}V^{\SSC\lor}\cong\wt\bM\hookrightarrow t^{-1}\wt\bD \quad\text{and}\quad\beta:\wt\bD\hookrightarrow\wt\bM\cong(\wt\bB^\dagger_\rig)\otimes_{\BQ_p}V^{\SSC\lor}
\]
and note that $t^{-1}\wt\bD\cong\bD\bigl(\ulTD\otimes\BOne(-1)\bigr)$.
Let $\nn=\dim V^{\SSC\lor}$ and fix a basis of $V^{\SSC\lor}$. This induces a basis of $\wt\bD=V^{\SSC\lor}\otimes_{\BQ_p}\wt\bB^\dagger_\rig$. With respect to this basis, $\alpha$ and $\beta$ are represented by matrices $A,B\in M_\nn(\wt\bB^\dagger_\rig)$ satisfying $AB=t\Id_\nn$. Then $A,B\in M_\nn(\wt\bB^+_\rig)\subset M_\nn(\bB^+_\cris)$ by \cite[Proposition I.4.1]{Berger04a}. So $A$ defines an isomorphism $\alpha_\cris:\bB_\cris^{\oplus \nn}\isoto \bigl(\ulTD\otimes\BOne(-1)\bigr)\otimes_{K_0}\bB_\cris$ compatible with the Frobenius on both sides. Moreover, it maps $(\bB_\cris^+)^{\oplus \nn}$ onto the preimage of $Fil^1_\mu \ulTD_C$ under the map $\id\otimes\theta:\bigl(\ulTD\otimes\BOne(-1)\bigr)\otimes_{K_0}\bB_\cris^+\to \ulTD\otimes_{K_0}C$. We denote this preimage by $Fil^1(\ulTD\otimes_{K_0}\bB_\cris)$. Tensoring $\alpha_\cris$ with $\BOne(1)\otimes_{K_0}\bB_\cris$ and precomposing with 
$\bB_\cris^{\oplus h}\isoto\BOne(1)\otimes_{K_0}\bB_\cris^{\oplus h},\,x\mapsto t^{-1}x$ yields an isomorphism $\bB_\cris^{\oplus h}\isoto\ulTD\otimes_{K_0}\bB_\cris$ compatible with Frobenius which maps $Fil^1\bB_\cris^{\oplus h}:=(t\cdot\bB_\cris^+)^{\oplus h}$ onto $Fil^1(\ulTD\otimes_{K_0}\bB_\cris)$. This means that $\ulTD$ is admissible in the sense of Faltings \cite[Definition 1]{Faltings08}. 

By \cite[Theorems 9 and 14]{Faltings08}, the filtered isocrystal $\ulTD$ comes from a $p$-divisible group $X$ over $\CO_\LL$ and a quasi-isogeny $\eta:\bX_{\CO_\LL/(p)}\to X_{\CO_\LL/(p)}$. Note that Faltings uses contravariant Dieudonn\'e modules and that his assertions for the covariant Dieudonn\'e functor $\BD$ mean $Fil^1_\mu D_\LL=\BD(\eta^{-1})(\Lie X\dual)_\LL\dual$. Therefore $\mu$ lies in the image of $\breve\pi^\an$ as desired.

\medskip\noindent
\ref{ThmC_2} 
This is proved by de Jong~\cite{dJ}. Indeed by \cite[Proposition 6.2]{dJ} the rational Tate module $V_pX_{\breve\CG^\an}$ is a local system of $\BQ_p$-vector spaces on $\breve\CG^\an$. In order that it defines a local system $\CV$ on $\breve\CF^0_b$ it suffices by \cite[Definition 4.1]{dJ} to show that
\begin{enumerate}
\item[(i)] $\breve\pi:\breve\CG^\an\to\breve\CF^0_b$ is a covering for the \'etale topology.
\item[(ii)] There is a descent datum $\wt\psi:pr_1^\ast V_pX_{\breve\CG^\an}\isoto pr_2^\ast V_pX_{\breve\CG^\an}$ over $\breve\CG^\an\times_{\breve\CF^0_b}\breve\CG^\an$ where $pr_i:\breve\CG^\an\times_{\breve\CF^0_b}\breve\CG^\an\to\breve\CG^\an$ is the projection onto the $i$-th factor, such that $\wt\psi$ satisfies the cocycle condition on $\breve\CG^\an\times_{\breve\CF^0_b}\breve\CG^\an\times_{\breve\CF^0_b}\breve\CG^\an$.
\end{enumerate}

Statement (i) was proved by Fargues~\cite[Lemma 2.3.24]{Fargues} and also follows from Faltings's result \cite[Theorem 14]{Faltings08} which says that $\breve\CF^0_b$ has an open covering $\breve\CF^0_b=\bigcup_i U_i$ such that the morphism $\coprod_i U_i\to \breve\CF^0_b$ factors through $\breve\pi:\breve\CG^\an\to\breve\CF^0_b$.

For (ii) consider the universal filtered isocrystal $\ul\CalD=(V_{K_0},b\phi,Fil^\bullet)$ over $\breve\CF^0_b$ and its canonical descent datum $\id:pr_1^\ast\breve\pi^\ast\ul\CalD=pr_2^\ast\breve\pi^\ast\ul\CalD$ over $\breve\CG^\an\times_{\breve\CF^0_b}\breve\CG^\an$. De Jong considers $pr_i^\ast X_{\breve\CG^\an}$ as an object of the stack $BT^\rig_{\BQ,\breve\CG^\an\times_{\breve\CF^0_b}\breve\CG^\an}$ of $p$-divisible groups up to isogeny over $\breve\CG^\an\times_{\breve\CF^0_b}\breve\CG^\an$. The functor from this stack to filtered isocrystals sends $pr_i^\ast X_{\breve\CG^\an}$ to $pr_i^\ast\breve\pi^\ast\ul\CalD$. Since this functor is fully faithful by \cite[Proposition 6.6]{dJ} one obtains a descent datum $\psi:pr_1^\ast X_{\breve\CG^\an}\isoto pr_2^\ast X_{\breve\CG^\an}$ in $BT^\rig_{\BQ,\breve\CG^\an\times_{\breve\CF^0_b}\breve\CG^\an}$. One now applies the Tate-module functor $V_p$ and takes $\wt\psi=V_p\psi:pr_1^\ast V_pX_{\breve\CG^\an}\isoto pr_2^\ast V_pX_{\breve\CG^\an}$. This yields the local system $\CV$ on $\breve\CF^0_b$.

The local system $\CV$ induces the desired tensor functor $\ul\CV$ because $V$ is a tensor generator of $\Rep_{\BQ_p}\!\!\!\GL(V)$. Let us give some more details. We do not know whether $\breve\CF^0_b$ is connected, although we expect this to hold; see Conjecture~\ref{Conj5.3c}. Therefore we consider each connected component $Y$ of it separately and fix a geometric base point $\bar\mu$ of $Y$. Then $\CV$ induces a representation 
\[
\pi_1^\et(Y,\bar\mu)\;\longto\;\GL(\CV_{\bar\mu})(\BQ_p)
\]
by Proposition~\ref{Prop1.2c}. We use the notation introduced before Corollary~\ref{Cor1.7}. Since $\dim V=\dim\CV_{\bar\mu}$ and $\Koh^1\bigl(\BQ_p,\GL(V)\bigr)=(0)$ by Hilbert 90, we can take $\wt\omega=\omega_0$ and $\wt G=\GL(V)\cong\GL(\CV_{\bar\mu})$. This defines a continuous group homomorphism $\pi_1^\et(Y,\bar\mu)\to\wt G(\BQ_p)$ and via Corollary~\ref{Cor1.7} a tensor functor $\ul\CV:\Rep_{\BQ_p}\!\!\!\GL(V)\to\PLoc_Y$ satisfying (\ref{EqRZb}). Tracing through the proof of Corollary~\ref{Cor1.7} we see that its composition with the forgetful fiber functor $\omega_{\bar\mu}$ satisfies $\omega_{\bar\mu}\circ\ul\CV\cong\wt\omega=\omega_0$. Hence $\omega_{\bar\mu}\circ\ul\CV(V)\cong V\cong\CV_{\bar\mu}$ and $\ul\CV(V)\cong\CV$. Moreover, for every $\mu\in\breve\CF^0_b(\LL)$ with $\LL/K_0$ finite the $\Gal(\LL^\alg/\LL)$-representations satisfy
\[
\ul\CV(V)_\mu\;\cong\; V_pX_{\CO_\LL}\;\cong\;V_\cris\bigl(\BD(X_{\CO_\LL}),(\Lie X_{\CO_\LL}\dual)\dual\bigr)\;\cong\;V_\cris\bigl(\ulD_{b,\mu}(V)\bigr)\,.
\]
This holds on all connected components $Y$ and therefore gives the tensor functor $\ul\CV$ satisfying (\ref{EqRZ}) of Conjecture~\ref{ConjRZ} on all of $\breve\CF^0_b$ as claimed.

\medskip\noindent
\ref{ThmC_3} 
We fix a geometric base point $\bar\mu$ of $\breve\CF^0_b$ and consider the canonical family of morphisms $f_\wt{K}:\breve\CG^\an_\wt{K}\to\HeckeTower_\wt{K}$ which sends $\alpha \wt{K}$ for $\alpha:\wt V\isoto(V_p X_{\breve\CG^\an})_{\bar\mu}$ to the $\wt{K}$-residue class $\beta \wt{K}$ of tensor isomorphisms $\beta:\wt\omega\isoto\omega_{\bar\mu}\circ\ul\CV$ induced by $\beta_V:=\alpha:\wt\omega(V)=\wt V\isoto(V_pX_{\breve\CG^\an})_{\bar\mu}=\omega_{\bar\mu}\circ\ul\CV(V)$ via the fact that $V$ is a tensor generator of $\Rep_{\BQ_p}\!\!\!G$. Again $f_\wt{K}$ is independent of the choice of $\bar\mu$ by \cite[Theorem 2.9]{dJ}. By construction the family $(f_\wt{K})_\wt{K}$ is equivariant for the Hecke action of $\wt G(\BQ_p)$ on both towers. For any algebraically closed complete extension $C$ of $\breve E$ the morphism $f_\wt{K}$ is bijective on $C$-valued points because the fibers of $\breve\CG^\an_\wt{K}(C)$ and $\HeckeTower_\wt{K}(C)$ over a fixed $C$-valued point of $\breve\CF^0_b$ are both isomorphic to the coset $\wt G(\BQ_p)/\wt{K}$ by \cite[Proposition 5.37]{RZ} and Remark~\ref{Rem1.8}(a). Hence $f_\wt{K}$ is quasi-finite by \cite[Proposition 3.1.4]{Berkovich2}. Since $\breve\CG^\an_\wt{K}$ and $\HeckeTower_\wt{K}$ are \'etale over $\breve\CF^0_b$ the morphisms $f_\wt{K}$ are \'etale by \cite[Corollary 3.3.9]{Berkovich2} and hence isomorphisms by Proposition~\ref{PropA.3}.

The statement about $\breve\CG^\an$ comes from the identification $\breve\CG^\an\cong\breve\CG^\an_{\GL(\wt\Lambda)}$ described before the theorem and the fact that $\HeckeTower_{\GL(\wt\Lambda)}$ is the space of $\BZ_p$-lattices inside $\CV$%
.

\medskip\noindent
\ref{ThmC_4} Recall that the action $\gamma:\breve\CG\to\breve\CG$ of $\gamma\in J(\BQ_p)$ is defined using the universal property of $\breve\CG$ by requiring that $\gamma^\ast(X_{\breve\CG},\eta)\cong(X_{\breve\CG},\eta\circ\gamma^{-1})$, i.e., there is an isomorphism $\Phi_\gamma\colon\gamma^*X_{\breve\CG}\isoto X_{\breve\CG}$ with $(\Phi_\gamma)_{\ol{\breve\CG}}\circ\gamma^\ast\eta=\eta\circ\gamma^{-1}$
. By rigidity of quasi-isogenies the isomorophism $\Phi_\gamma$ is uniquely determined and a straight forward calculation shows that $\Phi_\gamma\circ\gamma^\ast\Phi_\delta=\Phi_{\delta\gamma}$. By definition of the local system $\CV$ as $V_pX_{\breve\CG^\an}$ this yields a canonical $J(\BQ_p)$-linearization $V_p\Phi_\gamma:\gamma^*\ul\CV\isoto\ul\CV$ on $\ul\CV$.
Now the equivariance of the isomorphism $f_{\wt K}\colon\breve\CG^\an_{\wt K}\isoto\HeckeTower_{\wt K}\,,\,(X,\eta,\alpha\wt K)\mapsto \beta\wt K$ from \ref{ThmC_3} follows from the explicit description of the $J(\BQ_p)$-actions given in \eqref{EqGammaAction} and before Theorem~\ref{ThmC}.
\end{proof}

\begin{remark}\label{RemConnected1}
(a) If the Newton polygon and the Hodge polygon of $b$ have no point in common except for the two endpoints, then Chen~\cite[Lemme 5.1.1.1 and Th\'eor\`eme 5.1.1.1]{ChenThesis} proved that the image of the period morphism $\breve\CF_b^0$ is (arcwise) connected.

\medskip\noindent
(b)  In the case when $J$ is an inner form of $G$, that is, when the isocrystal $(V_{K_0},b\cdot\phi)$ is isoclinic, Kottwitz has formulated a conjecture on how the compactly supported $\ell$-adic \'etale cohomology of the towers in Theorem~\ref{ThmC}\ref{ThmC_3} decomposes as a representation of $\wt G(\BQ_p)\times J(\BQ_p)\times W_E$, where $W_E$ is the Weil group of $E$; see \cite[Conjecture 5.1]{Rapoport95}.

\medskip\noindent \label{RemFaltings}
(c) Let $\ulD:=(V_{K_0},b\cdot\phi,Fil^\bullet_\mu)$ with $\mu\in\breve\CF^\an$ be a filtered isocrystal satisfying \eqref{Eq2.2} for $n=1$. One can show that Faltings's notion of admissibility \cite[Definition~1]{Faltings08} for $\ulD$ is equivalent to our condition $\bM(\ulD)\cong\bM(0,1)\otimes_{\BQ_p}V$ from \eqref{EqFAdm}. In fact, our proof of Theorem~\ref{ThmC}\ref{ThmC_1} shows that our condition implies Faltings's admissibility for $\ulTD:=\CHom(\ulD,\BOne(1))$. Note also that $\ulTD$ is admissible if and only if $\ulD$ is by \cite[Remarks after Definition~1]{Faltings08}, and that $\bM(\ulTD)\cong\bM(\ulD)\dual$ by Lemma~\ref{Lemma4.2'} and Proposition~\ref{Prop4.10}.

To prove the converse we use the notation of the proof of Theorem~\ref{ThmC}\ref{ThmC_1} and assume that $\ulTD$ is admissible in the sense of Faltings via an isomorphism $\bB_\cris^{\oplus h}\isoto\ulTD\otimes_{K_0}\bB_\cris$ compatible with Frobenius and the filtrations $Fil^1$. Tensoring this isomorphism with $\BOne(-1)\otimes_{K_0}\bB_\cris$ and precomposing with 
$\bB_\cris^{\oplus h}\isoto\BOne(-1)\otimes_{K_0}\bB_\cris^{\oplus h},\,x\mapsto tx$ yields an isomorphism $\alpha_\cris\colon\bB_\cris^{\oplus h}\isoto\bigl(\ulTD\otimes\BOne(-1)\bigr)\otimes_{K_0}\bB_\cris$ compatible with Frobenius, which maps $(\bB_\cris^+)^{\oplus h}$ onto $Fil^1(\ulTD\otimes_{K_0}\bB_\cris)\subset(\ulTD\otimes\BOne(-1))\otimes_{K_0}\bB^+_\cris$. Since $\wt\bB^+_\rig\,=\,\bigcap_{n\in\BN_0}\phi^n \bB^+_\cris$ the isomorphism $\alpha_\cris$ comes from a monomorphism $(\wt\bB^+_\rig)^{\oplus h}\into\bigl(\ulTD\otimes\BOne(-1)\bigr)\otimes_{K_0}\wt\bB_\rig^+$, which by flatness of $\wt\bB^\dagger_\rig$ over $\wt\bB^+_\rig$ induces a monomorphism 
\[
\bM(0,1)^{\oplus h}\;\into\;\bigl(\ulTD\otimes\BOne(-1)\bigr)\otimes_{K_0}\wt\bB_\rig^\dagger\;\cong\; t^{-1}\bD(\ulTD).
\]
By construction of $\bM(\ulTD)$ in \ref{Def1.12b} and the fact that $\alpha_\cris(\bB_\cris^+)^{\oplus h}\;=\;Fil^1(\ulTD\otimes_{K_0}\bB_\cris)$, this monomorphism maps $\bM(0,1)^{\oplus h}$ isomorphically onto $\bM(\ulTD)$. This proves that 
\[
\bM(\ulD)\cong\bM(\ulTD)\dual\cong(\bM(0,1)^{\oplus h})\dual\cong\bM(0,1)^{\oplus h}.
\]
Therefore Faltings's notion of admissibility is equivalent to ours. Using this equivalence, the openness of $\breve\CF_b^0$ in Theorem~\ref{Thm2.1} also follows from \cite[Corollary~10]{Faltings08}.

\medskip\noindent
(d) There is an alternative proof by Scholze and Weinstein~\cite{ScholzeWeinstein}, not relying on Faltings's \cite{Faltings08} results, which shows that $\breve\CF_b^0$ equals the image of the period morphism $\breve\pi$. By using this alternative proof and our Theorem~\ref{Thm2.1}, one could make the whole treatment in this article independent of \cite{Faltings08}.
\end{remark}

\begin{corollary}\label{Cor6.4}
Let $G=\GL(V)$, let $\tilde b\in G(K_0)$, and let $\{\tilde\mu\}$ be a conjugacy class of cocharacters $\BG_m\to G$ such that $V$ satisfies (\ref{Eq2.2}) for some integer $\tilde\hhh$. Then there is a tensor functor $\ul{\wt\CV}$ from $\Rep_{\BQ_p}\!\!\!G$ to the category $\PLoc_{\breve\CF_{\tilde b}^0}$ of local systems of $\BQ_p$-vector spaces on $\breve\CF_{\tilde b}^0$ satisfying (\ref{EqRZ}) of Conjecture~\ref{ConjRZ}.
\end{corollary}

\begin{proof}
Like at the beginning of this section we set $b:=p^{1-\tilde\hhh}\!\cdot\!\tilde b\in G(K_0)$ and define the cocharacter $\mu\colon\BG_m\to G$ by $\mu(z)=z^{1-\tilde n}\!\cdot\!\tilde\mu(z)$ to obtain identifications $\ulD_{b,\mu}(V)=\ulD_{\tilde b,\tilde\mu}(V)\otimes\BOne(1-\tilde\hhh)$ and $\breve\CF_{\tilde b}^0\isoto\breve\CF_b^0\,,\,\tilde\mu\mapsto\mu$. Then Theorem~\ref{ThmC}\ref{ThmC_2} yields a tensor functor $\ul\CV$, which for every $\LL$-valued point $\tilde\mu\in\breve\CF_{\tilde b}^0$ for $\LL/K_0$ finite satisfies
\begin{eqnarray*}
\ul\CV(V)_\mu\;\cong\;V_\cris\bigl(\ulD_{b,\mu}(V)\bigr)\;=\;V_\cris\bigl(\ulD_{\tilde b,\tilde\mu}(V)\bigr)\otimes\BQ_p(\tilde\hhh-1)
\end{eqnarray*}
as $\Gal(\LL^\alg/\LL)$-representations because $V_\cris\bigl(\BOne(1)\bigr)=\BQ_p(-1)$. Thus it remains to twist $\ul\CV$ with the cyclotomic character. In concrete terms this means the following. On each connected component $Y$ of $\breve\CF_{\tilde b}^0$ with some geometric base point $\bar\mu$ the tensor functor $\ul\CV$ corresponds to a representation $\rho':\pi_1^\et(Y,\bar\mu)\to\wt G(\BQ_p)\cong\GL_\nn(\BQ_p)$. We now consider the representation $\tilde\rho':\pi_1^\et(Y,\bar\mu)\to\GL_\nn(\BQ_p)$ defined by $\tilde\rho'(\gamma):=\chi_{\rm cyc}(\gamma)^{1-\tilde\hhh}\cdot\rho'(\gamma)$ where $\chi_{\rm cyc}:\pi_1^\et(Y,\bar\mu)\to\Gal(\breve E^\alg/\breve E)\xrightarrow{\;\chi_{_{\rm cyc}}\,}\BZ_p\mal$ is the cyclotomic character of the base field $\breve E$. Let $\ul{\wt\CV}$ be the tensor functor from $\Rep_{\BQ_p}\!\!\!G$ to $\PLoc_{\breve\CF_{\tilde b}^0}$ associated with $\tilde\rho'$ by Corollary~\ref{Cor1.7}. It is independent of the base point $\bar\mu$ and hence exists on all of $\breve\CF_{\tilde b}^0$. Then there is an isomorphism of $\Gal(\LL^\alg/\LL)$-representations $\ul{\wt\CV}(V)_{\tilde\mu}\cong V_\cris\bigl(\ulD_{\tilde b,\tilde\mu}(V)\bigr)$ for all $\tilde\mu$ and this proves the corollary.
\end{proof}

Finally we can give the 

\begin{proof}[Proof of Theorem~\ref{ThmDOR}]
Let $\rho:G\hookrightarrow\GL(V)=:G'$ be a representation satisfying (\ref{Eq2.2}) and consider the induced embedding $\breve\CF\hookrightarrow\breve\CF':=\CF lag(V)$. By Proposition~\ref{PropDOR} we have $\breve\CF^0_b=\breve\CF^\an\cap\breve\CF'_b{}^0$. Now the local system $\ul{\wt\CV}(V)$ of Corollary~\ref{Cor6.4} on $\breve\CF'_b{}^0$ pulls back to a local system of $\BQ_p$-vector spaces on $\breve\CF^0_b$ which has the required property.
\end{proof}

\begin{conjecture}\label{Conj6.5}
The open subspace $\breve\CF^0_{\tilde b}$ from Corollary~\ref{Cor6.4} is the largest open $\breve E$-analytic subspace $\breve\CF^a_{\tilde b}$ of $\breve\CF^{wa}_{\tilde b}$ on which a tensor functor $\ul\CV:\Rep_{\BQ_p}\!\!\!G\to\PLoc_{\breve\CF_{\tilde b}^a}$ exists which satisfies (\ref{EqRZ}) of Conjecture~\ref{ConjRZ}.
\end{conjecture}

In particular, we expect that $\breve\CF^0_{\tilde b}$ and $\ul\CV$ from Corollary~\ref{Cor6.4} solve the problem posed in Conjecture~\ref{ConjRZ} for $G=\GL(V)$. Clearly the conjecture follows from Conjecture~\ref{Conj5.10} because $\ul\CV(V)$ is a local system as in Conjecture~\ref{Conj5.10}. We give some evidence for both conjectures in Section~\ref{SectConjectures}.

%
%

\section{PEL period morphisms} \label{SectPELPeriodM}
\setcounter{equation}{0}

Rapoport and Zink also consider period morphisms in the PEL situation with parahoric level structure. They fix data as follows:  
\begin{description}
\item[Case EL:] \mbox{ }\\
Let $B$ be a finite semi-simple algebra over $\BQ_p$.\\
Let $\CO_B$ be a maximal order in $B$.\\
Let $V$ be a finite dimensional $B$-module.\\
Let $G=\GL_B(V)$ as an algebraic group over $\BQ_p$.
\item[Case PEL:]\mbox{ }\\
Let $B,\CO_B,V$ be as in case EL.\\
Let $(\,.\;,\,.\,)$ be a non-degenerate alternating $\BQ_p$-bilinear form on $V$.\\
Let $\ast:a\mapsto a^\ast$ be an involution of $B$ which satisfies
\[
(av,w)\;=\;(v,a^\ast w)\qquad\text{for all }v,w\in V.
\]
Let $G$ be the algebraic group over $\BQ_p$ whose points with values in a $\BQ_p$-algebra $R$ are\\[2mm]
$G(R)\;=\;\bigl\{\,g\in \GL_{B\otimes R}(V_R):\;\exists\;c(g)\in R\text{ with }(gv,gw)=c(g)(v,w)\es\forall\,v,w\in V_R\,\bigr\}$,\\[2mm]
where $V_R:=V\otimes_{\BQ_p}R$.
\end{description}
In addition they fix an element $b\in G(K_0)$, a conjugacy class of cocharacters $\mu:\BG_m\to G$ and a multichain $\CL$ of $\CO_B$-lattices in $V$ \cite[Definition 3.4]{RZ}. They assume that the representation $V$ of $G$ satisfies our condition (\ref{Eq2.2}) with $\hhh=1$. In case PEL they assume in addition that the character $c:G\to\BG_m$ pairs with the slope quasi-cocharacter $\nu\in\Hom_{K_0}(\BG_m,G)\otimes_{\BZ}\BQ$ to $1$, and that $\CL$ is self-dual \cite[Definition 3.13]{RZ}.

Under these assumptions there is a $p$-divisible group $\bX$ over $\BF_p^\alg$ whose covariant Dieudonn\'e module is $(V_{K_0}\,,\,b\!\cdot\!\phi)$. In particular the group $\bX$ is equipped with an action $i_\bX:B\to\End(\bX)\otimes_{\BZ_p}\BQ_p$ and, in case PEL a polarization $\lambda:\bX\to\bX\dual$ with $\lambda\dual=-\lambda$, which is uniquely determined up to multiplication by an element of $\BQ_p\mal$; see \cite[3.20]{RZ}. Let $E$ be the Shimura field, that is the field of definition of the conjugacy class of $\mu$, and let $\breve E=E\cdot K_0$. Rapoport and Zink consider the following moduli problem.

\begin{definition}[{\cite[Definition 3.21]{RZ}}]
Let $\breve\CG:\Nilp_{\CO_{\!\breve E}}\to\Sets$ be the functor which assigns to $S\in\Nilp_{\CO_{\!\breve E}}$ the following data up to isomorphism.
\begin{enumerate}
\item For each lattice $\Lambda\in\CL$ a $p$-divisible group $X_\Lambda$ over $S$, with an $\CO_B$-action, $i_\Lambda:\CO_B\to\End X_\Lambda$, and
\item for each lattice $\Lambda\in\CL$ a quasi-isogeny $\eta_\Lambda:\bX_{\bar S}\to X_{\Lambda,\bar S}$ which commutes with the action of $\CO_B$, (here again $\bar S$ is the closed subscheme of $S$ defined by the ideal $p\CO_S$)
\end{enumerate}
such that certain conditions (i)--(v) are satisfied for which we refer to \cite[Definition 3.21 and the discussion on pp.~82--88]{RZ}.
\end{definition}

\begin{proposition}[{\cite[Theorem 3.25]{RZ}}]\label{Prop7.1}
The functor $\breve\CG:\Nilp_{\CO_{\!\breve E}}\to\Sets$ is representable by a formal scheme locally formally of finite type over $\CO_{\!\breve E}$.
\end{proposition}

The group $J(\BQ_p)$ from \eqref{EqJ} equals the group of quasi-isogenies of $\bX$ that respect the $B$-action and the polarization. This group acts on $\breve\CG$ by $\gamma\colon(X_\Lambda,\eta_\Lambda)_\Lambda\mapsto(X_\Lambda,\eta_\Lambda\circ\gamma^{-1})_\Lambda$ for $\gamma\in J(\BQ_p)$.

The period morphism $\breve\CG^\an\to\breve\CF_b^{wa}$ is defined as follows. Fix a $\Lambda\in\CL$, let $X_\Lambda$ be the universal $p$-divisible group over $\breve\CG$, and consider the exact sequence
\[
0\es\longto\es(\Lie X_\Lambda\dual)\dual \es\longto\es\BD(X_\Lambda)_{\breve\CG}\es\longto\es\Lie X_\Lambda\es\longto\es0\,.
\]
Then the quasi-isogeny $\eta_\Lambda$ induces an isomorphism
\[
\BD(\eta_\Lambda)_{\breve\CG^\an}:\es V\otimes_{\BQ_p}\CO_{\breve\CG^\an}\;=\;\BD(\bX)_{\breve\CG^\an}\;\isoto\;\BD(X_\Lambda)_{\breve\CG^\an}
\]
that makes $\BD(\eta_\Lambda)_{\breve\CG^\an}^{-1}(\Lie
X\dual)\dual_{\breve\CG^\an}$ into a family over $\breve\CG^\an$ of $B$-invariant subspaces of $V$, which in case PEL are totally isotropic. These subspaces are classified by $\breve\CF^\an$ and this defines the \emph{period morphism} $\breve\pi:\breve\CG^\an\to\breve\CF^\an$. The whole construction does not dependent on the choice of $\Lambda\in\CL$; see \cite[5.16]{RZ}. The period morphism factors through $\breve\CF_b^{wa}$ and by the same argument as in Proposition~\ref{Prop2.4} it even factors through $\breve\CF^0_b$. We will show in Theorem~\ref{ThmPEL} below, that its image equals $\breve\CF^0_b$. In Remark~\ref{Rem7.5}(d) we will give an explanation for the fact noticed by Rapoport and Zink \cite[5.53]{RZ}, that $\breve\pi:\breve\CG^\an\to\breve\CF_b^{wa}$ is in general not quasi-compact. The period morphism is also equivariant with respect to the $J(\BQ_p)$-action on $\breve\CF^\an$ under which $\gamma\in J(\BQ_p)\subset G(K_0)$ sends $\mu\in\breve\CF^\an$ to $\gamma\mu\gamma^{-1}$.

\bigskip

Recall the  cohomology class
\[
cl(b,\mu)\;=\;\kappa(b)-\mu^\#\;\in\;\Koh^1(\BQ_p,G)
\]
and the associated fiber functor $\wt\omega:\Rep_{\BQ_p}\!\!\!G\to(\BQ_p\text{-vector spaces})$ and inner form $\wt G$ of $G$ described before Corollary~\ref{Cor1.7}. Let $\wt V:=\wt\omega(V)$. Since $B$ acts as endomorphisms of the representation $V$ of $G$, also $\wt V$ is a $B$-module. In the PEL case $\wt V$ is equipped with a non-degenerate alternating $\BQ_p$-bilinear form $(\,.\,,\,.\,)^{^\sim}$ such that 
\begin{equation}\label{Eq7.1}
(a\tilde v,\tilde w)^{^\sim}\;=\;(\tilde v,a^\ast\tilde w)^{^\sim}\qquad\text{for all }\tilde v,\tilde w\in\wt V \text{ and all }a\in B\,.
\end{equation}
Namely, the original form $(\,.\,,\,.\,):V\otimes_{\BQ_p}V\to\BQ_p$ is a morphism from the $G$-representation $V\otimes_{\BQ_p}V$ to the $G$-representation $\BQ_p$ on which $G$ acts through the character $c$. We let $(\,.\,,\,.\,)^{^\sim}$ be the image of $(\,.\,,\,.\,)$ under $\wt\omega$. Formula (\ref{Eq7.1}), being an equality of morphisms $\wt V\otimes_{\BQ_p}\wt V\to\BQ_p$, is the image under $\wt\omega$ of the corresponding equality for $V$. By construction of $\wt G$ we have $\wt G\cong\Aut^\otimes(\wt\omega)$ and $\Rep_{\BQ_p}\!\!\!G\cong\Rep_{\BQ_p}\!\!\!\wt G$; see \cite[Proposition 2.8 and Theorem 2.11]{DM}. Applying $\wt\omega$ to the character $c:G\to\BG_m$ we obtain another character which we view as a $\wt G$-representation $\tilde c:\wt G\to\BG_m$. In this way we get identifications of algebraic groups over $\BQ_p$
\begin{align}\label{WTG}
\wt G&\;=\;\GL_B(\wt V) &&\text{in case EL,}\\
\wt G&\;=\;\bigl\{\,\tilde g\in\GL_B(\wt V):\es(\tilde g\tilde v,\tilde g\tilde w)^{^\sim}=\,\tilde c(\tilde g)(\tilde v,\tilde w)^{^\sim}\quad\forall \,\tilde v,\tilde w\in\wt V\,\bigr\}&&\text{in case PEL.}\nonumber
\end{align}

\bigskip

Consider a decomposition of $B$ into a product $\prod_{i=1}^r B_i$ of simple algebras $B_i\cong M_{n_i}(D_i)$, which are matrix algebras over division algebras $D_i$, in such a way that $\CO_B\cong M_{n_i}(\CO_{D_i})$. Then any $\CO_B$-module $\Lambda$ decomposes accordingly $\Lambda\cong\bigoplus_{i=1}^r\Lambda_i$. 

\begin{definition}\label{Def7.2}
In the EL case a \emph{multichain} of type $(\CL)$ of $\CO_B$-lattices in $\wt V$ is a functor $\wt{\es}:\Lambda\mapsto\wt\Lambda$ from $\CL$ to the category of $\CO_B$-lattices in $\wt V$ such that
\begin{enumerate}
\item[($\wt{\rm ii}$)] If $\Lambda\subset\Lambda'$ then $\wt\Lambda\subset\wt{\Lambda'\,}\!$ and $[\wt{\Lambda'_i}:\wt\Lambda_i]=[\Lambda'_i:\Lambda_i]$ for $i=1,\ldots,r$.
\item[($\wt{\rm iii}$)] For any $a\in B$ with $a^{-1}\CO_Ba=\CO_B$ we have $\wt{a\Lambda}=a\cdot\wt\Lambda$\,.
\end{enumerate}
In the PEL case we speak of a \emph{polarized multichain} of type $(\CL)$ and we require in addition that there is a unit $\tilde\ell\in\BQ_p\mal$ with
\begin{enumerate}
\item[($\wt{\rm v}$)] $\tilde\ell\cdot\wt{(\Lambda\dual)}=(\wt\Lambda)\dual:=\{\,\tilde v\in\wt V:\;(\tilde v,\tilde w)^{^\sim}\in\BZ_p\text{ for all }\tilde w\in\wt\Lambda\,\}$ for every $\Lambda\in\CL$.
\end{enumerate}
The unusual numbering of our conditions is chosen because the items of our definition correspond to the respective items of \cite[Definition 3.21]{RZ}.
\end{definition}

For the next theorem again recall from Remark~\ref{Rem1.8}(a) that with each tensor functor $\ul\CV:\Rep_{\BQ_p}\!\!\!G\to\PLoc_{\breve\CF^0_b}$ one associates a tower of \'etale covering spaces $\HeckeTower_{\wt K}$ of $\breve\CF^0_b$ for $\wt K\subset\wt G(\BQ_p)$ compact open. Moreover, we again consider the tower of \'etale coverings $\breve\CG^\an_{\wt K}$ of $\breve\CG^\an$ with its commuting actions of $J(\BQ_p)$ by $\gamma\colon(X_\Lambda,\eta_\Lambda,\alpha\wt K)\mapsto (X_\Lambda,\eta_\Lambda\circ\gamma^{-1},\alpha\wt K)$ for $\gamma\in J(\BQ_p)$, and of $\wt G(\BQ_p)$ by Hecke-correspondences; see Rapoport and Zink \cite[5.32--5.39]{RZ} and our discussion before Theorem~\ref{ThmC}.

\begin{theorem}\label{ThmPEL} 
Assume that $\breve\CG^\an\ne\emptyset$.
\begin{enumerate}
\item \label{ThmPEL_1}
Then the image $\breve\pi(\breve\CG^\an)$ of the EL or PEL period morphism $\breve\pi$ equals the open $\breve E$-analytic subspace $\breve\CF_b^0\subset\breve\CF^\an$ from Theorem~\ref{Thm2.1}.
\item  \label{ThmPEL_2} Fix any lattice $\Lambda\in\CL$. Then the rational Tate module $V_p X_{\Lambda,\breve\CG^\an}$ of the universal $p$-divisible group $X_{\Lambda,\breve\CG^\an}$ over $\breve\CG^\an$ descends to a local system $\CV$ of $\BQ_p$-vector spaces on $\breve\CF_b^0$ which is independent of $\Lambda$ up to canonical isomorphism. This induces a tensor functor $\ul\CV$ from $\Rep_{\BQ_p}\!\!\!G$ to the category $\PLoc_{\breve\CF_b^0}$ of local systems of $\BQ_p$-vector spaces on $\breve\CF_b^0$ satisfying (\ref{EqRZ}) of Conjecture~\ref{ConjRZ} and $\ul\CV(V)\cong\CV$.
\item \label{ThmPEL_3}
The tower of $\breve E$-analytic spaces $(\breve\CG^\an_{\wt K})_{\wt K\subset \wt G(\BQ_p)}$ is canonically isomorphic over $\breve\CF^0_b$ in a Hecke equivariant way to the tower of \'etale covering spaces $(\HeckeTower_{\wt K})_{\wt K\subset \wt G(\BQ_p)}$ of $\breve\CF^0_b$ associated with the tensor functor $\ul\CV$. In particular, $\breve\CG^\an$ is (non-canonically) isomorphic to the space of (polarized in case PEL) multichains of type $(\CL)$ inside the local system $\CV$ (which are defined as in Definition~\ref{Def7.2} using Proposition~\ref{Prop1.2c}).
\item \label{ThmPEL_4}
The tensor functor $\ul\CV$ carries a canonical $J(\BQ_p)$-linearization which by Remark~\ref{Rem1.8}(b) induces an action of $J(\BQ_p)$ on the tower $(\HeckeTower_\wt{K})_{\wt{K}\subset \wt G(\BQ_p)}$. The isomorphism from \ref{ThmPEL_3} is equivariant for the $J(\BQ_p)$-action on both towers.
\end{enumerate}
\end{theorem}

\begin{remark}\label{Rem7.5}
(a)  For each connected component $Y$ of $\breve\CF^0_b$ and each geometric base point $\bar\mu$ of $Y$ the tensor functor $\ul\CV$ corresponds to a representation $\rho':\pi^\et_1(Y,\bar\mu)\to\wt G(\BQ_p)$ by Corollary~\ref{Cor1.7}. In the course of the proof we show that the composition $\tilde c\circ\rho'$ with the multiplier $\tilde c$ of $\wt G$ equals the cyclotomic character $\pi^\et_1(Y,\bar\mu)\to\Gal(\breve E^\alg/\breve E)\xrightarrow{\;\chi_{_{\rm cyc}}\,}\BZ_p\mal$ of the base field $\breve E$. 

\medskip\noindent
(b) If $F/\BQ_p$ is unramified and the Newton polygon and the Hodge polygon of $b$ have no point in common except for the two endpoints, then Chen, Kisin, and Viehmann~\cite{CKV13} proved that the image of the period morphism $\breve\CF_b^0$ is (arcwise) connected; see also \cite[Lemme 5.1.1.1]{ChenThesis}.

\medskip\noindent
(c) Theorem~\ref{ThmPEL}\ref{ThmPEL_3} implies that $\breve\CG^\an_{\wt K}\times_{\breve\CF_b^0}Y$ is the \'etale covering space of $Y$ corresponding to the $\pi_1^\et(Y,\bar\mu)$-set $\Isom^\otimes(\wt\omega,\omega_{\bar\mu}\circ\ul\CV)/\wt K$, which can be identified non-canonically with $\wt G(\BQ_p)/{\wt K}$; see Remark~\ref{Rem1.8}. This had been conjectured by Rapoport and Zink \cite[1.37]{RZ}. 

\medskip\noindent
(d) The theorem gives a natural explanation for the fact that in general $\breve\pi$ is not quasi-compact, which was discovered in \cite[5.53]{RZ}. Firstly the fibers of $\breve\CG^\an$ over $\breve\CF_b^0$ are in bijection with a coset of $\wt G(\BQ_p)$ modulo a compact open subgroup; see the proof of Theorem~\ref{ThmPEL}\ref{ThmPEL_3}. In particular the fibers are infinite in general. Secondly, if the inclusion $\breve\CF_b^0\subset\breve\CF_b^{wa}$ were quasi-compact then it would be an isomorphism. But Example~\ref{Example6.4} and Theorem~\ref{Thm8.2} show that this occurs only in a few low-dimensional cases.

\medskip\noindent
(e) In the case when $b$ is \emph{basic}, that is, when $J$ is an inner form of $G$, Kottwitz has formulated a conjecture on how the compactly supported $\ell$-adic \'etale cohomology of the towers in Theorem~\ref{ThmPEL}\ref{ThmPEL_3} decomposes as a representation of $\wt G(\BQ_p)\times J(\BQ_p)\times W_E$, where $W_E$ is the Weil group of $E$; see \cite[Conjecture 5.1]{Rapoport95}.
\end{remark}

Before proving the theorem let us mention the following

\begin{conjecture}\label{Conj7.6}
The open subspace $\breve\CF^0_b$ is the largest open $\breve E$-analytic subspace $\breve\CF^a_b$ of $\breve\CF^{wa}_b$ on which a tensor functor $\Rep_{\BQ_p}\!\!\!G\to\PLoc_{\breve\CF_b^a}$ exists which satisfies (\ref{EqRZ}) of Conjecture~\ref{ConjRZ}.
\end{conjecture}

In particular, we expect that $\breve\CF^0_b$ and $\ul\CV$ from Theorem~\ref{ThmPEL}\ref{ThmPEL_2} solve the problem posed in Conjecture~\ref{ConjRZ}. Again the conjecture follows from Conjecture~\ref{Conj5.10} because $\ul\CV(V)$ is a local system as in Conjecture~\ref{Conj5.10}. We give some evidence for both conjectures in Section~\ref{SectConjectures}. We now turn to the 

\begin{proof}[Proof of Theorem~\ref{ThmPEL}]
\ref{ThmPEL_1} \es 1. By our assumption $\breve\CG^\an\ne\emptyset$, there is a finite field extension $L/\breve E$ and a point $x_0\in\breve\CG(\CO_L)$ corresponding to an $\CL$-set of $p$-divisible groups with $\CO_B$-action $\{X_\Lambda,i_\Lambda\}_{\Lambda\in\CL}$ over $\CO_L$ and $\CO_B$-equivariant quasi-isogenies $\eta_\Lambda:\bX_{\CO_L/(p)}\to X_{\Lambda,\CO_L/(p)}$ satisfying conditions (i)--(v) of \cite[Definition 3.21]{RZ}. In particular $\eta_{\Lambda'}\eta_\Lambda^{-1}$ lifts to a quasi-isogeny $\eta_{\Lambda',\Lambda}:X_\Lambda\to X_{\Lambda'}$ over $\CO_L$ for all $\Lambda,\Lambda'\in\CL$. Let $\mu_0=\breve\pi(x_0)\in\breve\CF^0_b$ be the image of $x_0$. By Wintenberger's result~\cite[Corollary to Proposition 4.5.3]{Wintenberger97} the rational Tate modules $V_p X_{\Lambda,L^\alg}$ are all $B$-isomorphic to $\wt V$ via $B$-isomorphisms $\delta_\Lambda:V_p X_\Lambda\isoto \wt V$ satisfying $\delta_\Lambda=\delta_{\Lambda'}\circ V_p(\eta_{\Lambda',\Lambda})$. We claim that $\delta_\Lambda(T_p X_\Lambda)=:\wt\Lambda\subset\wt V$ forms a (polarized) multichain of type $(\CL)$ of $\CO_B$-lattices in $\wt V$. Indeed the decomposition $\CO_B\cong\prod_{i=1}^r M_{n_i}(\CO_{D_i})$ induces a decomposition $X_\Lambda\cong\prod_{i=1}^r X_{\Lambda_i}$ of $p$-divisible groups. If $\Lambda\subset\Lambda'$ then \cite[Definition 3.21(ii) and (ii bis) from p.~88]{RZ} implies that $\eta_{\Lambda',\Lambda}$ is an isogeny whose $i$-th composition factor $\eta_{\Lambda'_i,\Lambda_i}:X_{\Lambda_i}\to X_{\Lambda'_i}$ has height $\log_p[\Lambda'_i:\Lambda_i]$. Therefore we obtain $\wt\Lambda\subset\wt{\Lambda'\,}\!$ and $[\wt{\Lambda'_i}:\wt\Lambda_i]=p^{{\rm ht}(\eta_{\Lambda'_i,\Lambda_i})}=[\Lambda'_i:\Lambda_i]$ proving ($\wt{\rm ii}$).

If $a\in B\mal$ satisfies $a^{-1}\CO_B a=\CO_B$ then \cite[Definition 3.21(iii)]{RZ} says that $\eta_{a\Lambda,\Lambda}=\theta_{a,\Lambda}\circ i_\Lambda(a)^{-1}$ for an isomorphism $\theta_{a,\Lambda}:X_\Lambda\to X_{a\Lambda}$ of $p$-divisible groups. Therefore condition ($\wt{\rm iii}$) follows from
\[
\wt{a\Lambda}\;=\;\delta_\Lambda V_p(\eta_{\Lambda,a\Lambda})(T_p X_{a\Lambda})\;=\;\delta_\Lambda V_p\bigl(i_\Lambda(a)\,\theta_{a,\Lambda}^{-1}\bigr)(T_p X_{a\Lambda})\;=\;a\cdot\delta_\Lambda(T_p X_\Lambda)\;=\;a\cdot\wt\Lambda\,.
\] 

In case PEL there exists by \cite[Definition 3.21(v)]{RZ} a constant $\ell\in\BQ_p\mal$ such that for all $\Lambda\in\CL$ the quasi-isogeny $(\eta_{\Lambda\dual}\dual)^{-1}\circ\ell\lambda\circ\eta_\Lambda^{-1}$ 
\[
\xymatrix {
X_{\Lambda,\CO_L/(p)} \ar[r] & (X_{\Lambda\dual,\CO_L/(p)})\dual \ar[d]^{(\eta_{\Lambda\dual})\dual}\\
\bX_{\CO_L/(p)}\ar[u]^{\eta_\Lambda} \ar[r]^{\ell\lambda} & \bX\dual_{\CO_L/(p)}
}
\]
lifts to an isomorphism $p_\Lambda:X_\Lambda\isoto(X_{\Lambda\dual})\dual$. Here we obtain the polarization $\ell\lambda:\bX\to\bX\dual$ from the alternating form $(\,.\,,\,.\,)$ on $V$. That is $(\,.\,,\,.\,)$ on $V\otimes_{\BQ_p}K_0=\BD(\bX)_{K_0}$ coincides up to multiplication with a unit in $\BQ_p\mal$ with the pairing
\[
\BD(\bX)_{K_0}\otimes_{K_0}\BD(\bX)_{K_0} \;\xrightarrow{\es\BD(\ell\lambda)\otimes\id\;}\;\BD(\bX\dual)_{K_0}\otimes_{K_0}\BD(\bX)_{K_0}\;\xrightarrow{\es\BD(\psi_\bX)_{K_0}\;}\;\BD(\BG_m)_{K_0}\;=\;K_0
\]
induced from the pairing $\psi_\bX:\bX\dual\times\bX\to\BG_m$.
Now under the identification $\delta_\Lambda:V_p X_\Lambda\isoto\wt V$ the isomorphism $p_\Lambda$ induces an alternating form 
\begin{equation}\label{Eq7.2}
\wt V\otimes_{\BQ_p}\wt V\;\xrightarrow{\es V_p(p_{\Lambda\dual}\circ\eta_{\Lambda\dual,\Lambda})\delta_\Lambda^{-1}\otimes\delta_\Lambda^{-1}\;}\;V_p X_\Lambda\dual\otimes_{\BQ_p}V_p X_\Lambda \;\xrightarrow{\es V_p\psi_\Lambda\;}\; V_p\BG_{m,L^\alg}\;=\;\BQ_p
\end{equation}
which is independent of $\Lambda$. Up to multiplication with a unit $\tilde\ell\in\BQ_p\mal$ this alternating form equals $(\,.\,,\,.\,)^{^\sim}$, the image of $(\,.\,,\,.\,)$ under the fiber functor $\tilde\omega$ by Wintenberger~\cite[Corollary to Proposition 4.5.3]{Wintenberger97}. Hence ($\wt{\rm v}$) follows from
\begin{eqnarray*}
\wt{\Lambda\dual}& :=& \delta_\Lambda\circ V_p(\eta_{\Lambda,\Lambda\dual})(T_pX_{\Lambda\dual})\\[2mm]
& =& \delta_\Lambda\circ V_p(\eta_{\Lambda,\Lambda\dual}\circ p_{\Lambda\dual}^{-1})(T_pX_\Lambda\dual)\\[2mm]
& =& \bigl\{\,\tilde v\in\wt V:\;\tilde\ell\cdot(\tilde v,\tilde w)^{^\sim}\in\BZ_p\text{ for all }\tilde w\in\delta_\Lambda(T_p X_\Lambda)=\wt\Lambda\,\bigr\}\\[2mm]
& =& \tilde\ell^{-1}(\wt\Lambda)\dual\,.
\end{eqnarray*}

\medskip
\noindent
2. Now in addition to $\mu_0=\breve\pi(x_0)$ let $\mu\in\breve\CF^0_b$ be any point. Consider the canonical representation $\rho:G\into\GL(V)=:G'$ and the induced embedding $\breve\CF^0_b\into\breve\CF'^0_{\rho(b)},\,\mu\mapsto\rho\circ\mu$ into the corresponding open subspace $\breve\CF'^0_{\rho(b)}$ for $G'$. The properties of the latter subspace were discussed in the previous section. Moreover, $\breve\CF^0_b=\breve\CF^\an\cap\breve\CF'^0_{\rho(b)}$ by Proposition~\ref{PropDOR}. In addition we denote the formal scheme from Proposition~\ref{Prop6.1} in this section by $\breve\CG'$. 
By Theorem~\ref{ThmC} and Faltings's \cite[Theorem 14]{Faltings08} there is an admissible affine formal scheme $S'$ over $\CO_{\!\breve E}$, a $p$-divisible group $X$ over $S'$, and a quasi-isogeny $\eta:\bX_\ol{S'}\to X_\ol{S'}$ for $\ol{S'}=\Var(p)\subset S'$, such that the $S'$-valued point $(X,\eta)\in\breve\CG'(S')$ induces under the period morphism $\breve\pi':\breve\CG'{}^\an\to\breve\CF'^0_{\rho(b)}$ an affinoid subdomain $S'{}^\an\subset\breve\CF'^0_{\rho(b)}$ containing $\mu$. By changing $\mu_0$ we may assume that $S'{}^\an$ also contains $\mu_0$. Let $S\subset S'$ be the Zariski closed formal subscheme with $S^\an=\breve\CF^0_b\times_{\breve\CF'^0_{\rho(b)}}S'{}^\an$. We denote by $V_p X$ the $\pi_1^\et(S^\an,\bar\mu_0)$-module $V_p X_{\bar\mu_0}$ at some fixed geometric base point $\bar\mu_0$ of $S^\an$ above $\mu_0$. By Wintenberger~\cite[Corollary to Proposition 4.5.3]{Wintenberger97} for $(G',\rho(b))$ there is an isomorphism $\delta:V_p X\isoto\wt V$. 

Consider the (polarized) multichain of type $(\CL)$ in $\wt V$ constructed in part 1 above from an $\CL$-set of $p$-divisible groups with $\CO_B$-action above $x_0\in\breve\CG(\CO_L)$, which we now call $\{\ol X_\Lambda\}_{\Lambda\in\CL}$. 
We claim that there is a projective morphism $T\to S$ of admissible formal $\CO_{\!\breve E}$-schemes with $T^\an\to S^\an$ a finite \'etale covering, and for each $\Lambda\in\CL$ a $p$-divisible group $X_\Lambda$ over $T$ together with a quasi-isogeny $\eta_{X,\Lambda}:X_\Lambda\to X$ with $\delta_\Lambda(T_p X_\Lambda)=\wt\Lambda\subset\wt V$ for $\delta_\Lambda:=\delta\circ V_p(\eta_{X,\Lambda}):V_p X_\Lambda\isoto\wt V$.

Indeed, by using the periodicity $\wt{p\Lambda}=p\cdot\wt\Lambda$ which allows to set $\eta_{\,p^n\!\Lambda,\Lambda}:=p^{-n}:X_\Lambda\to X_\Lambda=:X_{p^n\!\Lambda}$, we only need to check this for finitely many lattices \mbox{$\Lambda_\nu\in\CL$} with \mbox{$\wt\Lambda_\nu\supset\delta(T_p X)$}. Let $T^\an\to S^\an$ be the finite \'etale covering corresponding to the $\pi_1^\et(S^\an,\bar\mu_0)$-set $\prod_\nu\delta^{-1}\wt\Lambda_\nu\big/T_p X$. That is, $\bar\mu_0$ lifts to a point of $T^\an$ and $\pi_1^\et(T^\an,\bar\mu_0)$ acts trivially on $\prod_\nu\delta^{-1}\wt\Lambda_\nu\big/T_p X$. This implies the existence of a $p$-divisible group $X_{\Lambda_\nu}$ over $T^\an$ and a quasi-isogeny $\eta_{\Lambda_\nu,X}:X\to X_{\Lambda_\nu}$ with $\ker\eta_{\Lambda_\nu,X}=\bigl(\delta^{-1}\wt\Lambda_\nu\big/T_p X\bigr)_{T^\an}$. The family of finite flat subgroup schemes $\{\,\ker\eta_{\Lambda_\nu,X}\subset X\,\}_\nu$ over $T^\an$ corresponds to a morphism $f$ from $T^\an$ to a moduli space of families of finite flat subgroup schemes of $X$. This moduli space is projective over $S$. Let $T$ be the scheme theoretic closure of the graph of this morphism $f$. Then $T\to S$ is projective, $\ker\eta_{\Lambda_\nu,X}$ extends to a finite flat subgroup scheme of $X$ over $T$, and we set $\eta_{\Lambda_\nu,X}:X\to X_{\Lambda_\nu}:=X/\ker\eta_{\lambda_\nu,X}$ and $\eta_{X,\Lambda_\nu}:=\eta_{\Lambda_\nu,X}^{-1}$. This proves the claim. 

We may replace $S$ by $T$ and $(X,\eta)$ by $(X_\Lambda,\eta_{X,\Lambda}^{-1}\circ\eta)$ for any $\Lambda$. The construction yields that the fiber over $\bar\mu_0$ of $X_\Lambda$ is isomorphic to $\ol X_\Lambda$. This implies that $V_p X_\Lambda$ is a $B$-module and $\delta_\Lambda$ is a $B$-isomorphism.

\medskip\noindent
3. The quasi-isogeny $\eta_\Lambda:=\eta_{X,\Lambda}^{-1}\circ\eta:\bX_{\ol{S}}\to X_{\Lambda,\ol{S}}$ induces an action $i_\Lambda(a)=\eta_\Lambda\circ i_\bX(a)\circ\eta_\Lambda^{-1}$ of $a\in B$ as quasi-isogenies on $X_{\Lambda,\ol{S}}$. For each $a\in B$ there is an integer $n$ such that $i_\Lambda(p^n a)$ is even an isogeny. We denote by $\BD(X_\Lambda)_{S}$ the Lie algebra of the universal vector extension of $X_\Lambda$ over $S$. It only depends on the reduction $X_{\Lambda,\ol{S}}$. We set $\BD(X_\Lambda)_{S^\an}=\BD(X_\Lambda)_{S}\otimes_{\CO_S}\CO_{S^\an}$ and denote by $\BD(\eta_\Lambda)_{S^\an}:V\otimes_{\BQ_p}\CO_{S^\an}=\BD(\bX)_{K_0}\otimes_{K_0}\CO_{S^\an}\isoto\BD(X_\Lambda)_{S^\an}$ the isomorphism induced from $\eta_\Lambda$ by the crystalline nature of $\BD(X_\Lambda)_S$. By definition of the period morphism $\breve\pi'$ 
\[
Fil^1 X_\Lambda\;=\;\BD(X_{\Lambda})_S\,\cap\,(Fil^1X_\Lambda)\otimes_{\CO_S}\CO_{S^\an}\;=\;\BD(X_{\Lambda})_S\,\cap\,\BD(\eta_\Lambda)_{S^\an}\bigl(Fil^1(V\otimes_{\BQ_p}\CO_{S^\an})\bigr)
\]
for the universal filtration $Fil^1(V\otimes_{\BQ_p}\CO_{S^\an})$. Indeed, the first equation holds by definition and for the second equation the inclusion ``$\subset$'' is clear. On the other hand, if $m$ belongs to the right hand side then $p^nm\in Fil^1X_\Lambda$ for $n\gg0$. Thus the image $\bar m$ of $m$ inside $\Lie X_\Lambda=\BD(X_\Lambda)_S/Fil^1X_\Lambda$ satisfies $p^n\bar m=0$. But $\Lie X_\Lambda$, as a locally free $\CO_S$-module, does not have $p$-torsion, and so $m\in Fil^1X_\Lambda$.

Since $Fil^1(V\otimes_{\BQ_p}\CO_{S^\an})$ is $B$-invariant, $i_\Lambda(p^n a)$ lifts to an isogeny of $X_\Lambda$ over $S$. Hence $i_\Lambda(a):=i_\Lambda(p^n a)\otimes p^{-n}\in\End(X_\Lambda)\otimes_{\BZ_p}\BQ_p$. By construction the $B$-module structure on $V_pX_\Lambda$ induced from $V_p(i_\Lambda(a))$ coincides with the structure considered at the end of part 2 above. For $a\in\CO_B$ then $a\cdot\wt\Lambda\subset\wt\Lambda$ implies $V_p(i_\Lambda(a))(T_p X_\Lambda)\subset T_p X_\Lambda$ and therefore $i_\Lambda(a)\in\End(X_\Lambda)$. In order that $(X_\Lambda,i_\Lambda,\eta_\Lambda)_{\Lambda\in\CL}$ is an $S$-valued point of $\breve\CG$ we have to verify conditions (i)--(v) of \cite[Definition 3.21]{RZ}.

By construction $\BD(\eta_\Lambda)_{S^\an}$ induces an isomorphism of exact sequences of $B$-modules
\[
\xymatrix @C=1pc {
0 \ar[r] & (Fil^1 X_\Lambda)\otimes_{\CO_S}{\CO_{S^\an}} \ar[r] & \BD(X_\Lambda)_{S^\an} \ar[r] & \Lie X_\Lambda\otimes_{\CO_S}{\CO_{S^\an}}\ar[r] & 0\\
0 \ar[r] & Fil^1(V\otimes_{\BQ_p}{\CO_{S^\an}}) \ar[r] \ar[u]^\sim & \BD(\bX)_{K_0}\otimes_{K_0}{\CO_{S^\an}} \ar[r] \ar[u]^\sim_{\BD(\eta_\Lambda)_{S^\an}} & gr^0(V\otimes_{\BQ_p}{\CO_{S^\an}}) \ar[u]^\sim \ar[r] & 0\\
0 \ar[r] & (V\otimes_{\BQ_p}\CO_{S^\an})^{\mu^{\rm univ}_{S^\an}(p)=p} \ar[r] \ar@{=}[u] & V\otimes_{\BQ_p}{\CO_{S^\an}} \ar[r] \ar@{=}[u] & (V\otimes_{\BQ_p}\CO_{S^\an})^{\mu^{\rm univ}_{S^\an}(p)=1} \ar[u]^\sim \ar[r] & 0\,,
}
\]
where $\mu^{\rm univ}_{S^\an}$ is the universal cocharacter on $S^\an$.
Hence $\det_{S}(a;\Lie X_\Lambda)=\det_{S^\an}\bigl(a;(V\otimes_{\BQ_p}\CO_{S^\an})^{\mu^{\rm univ}_{S^\an}(p)=1}\bigr)$ for every $a\in\CO_B$. This implies condition (iv). Note that (iv) further implies (i) by \cite[3.23(c)]{RZ}.

Next, for any $\Lambda\subset\Lambda'$ in $\CL$ our condition 
$\wt\Lambda\subset\wt{\Lambda'\,}\!$ implies that $V_p(\eta_{\Lambda',\Lambda})(T_p X_\Lambda)\subset T_p X_{\Lambda'}$ and so $\eta_{\Lambda',\Lambda}:X_\Lambda\to X_{\Lambda'}$ is an isogeny. Under the decompositions $\CO_B\cong\prod_{i=1}^r M_{n_i}(\CO_{D_i})$,\, $\Lambda\cong\bigoplus_i\Lambda_i$ and $X_\Lambda\cong\prod_{i=1}^r X_{\Lambda_i}$ it induces an isogeny $X_{\Lambda_i}\to X_{\Lambda'_i}$ of height equal to $\log_p[T_pX_{\Lambda'_i}:V_p(\eta_{\Lambda',\Lambda})(T_p X_{\Lambda_i})]=\log_p[\wt{\Lambda'_i}:\wt\Lambda_i]=\log_p[\Lambda'_i:\Lambda_i]$. Condition (ii) follows from this by \cite[3.23(d)]{RZ}.

If $a\in B\mal$ satisfies $a^{-1}\CO_B a=\CO_B$ then the quasi-isogeny $\theta_{a,\Lambda}=\eta_{a\Lambda,\Lambda}\circ i_\Lambda(a):X_\Lambda\to X_{a\Lambda}$ satisfies 
\[
V_p(\theta_{a,\Lambda})(T_p X_\Lambda)\;=\;V_p(\theta_{a,\Lambda})\circ\delta_\Lambda^{-1}(\wt\Lambda)\;=\;V_p(\eta_{a\Lambda,\Lambda})\circ\delta_\Lambda^{-1}(a\cdot\wt\Lambda)\;=\;\delta_{a\Lambda}^{-1}(\wt{a\Lambda})\;=\;T_p X_{a\Lambda}\,.
\]
Hence $\theta_{a,\Lambda}:X_\Lambda\to X_{a\Lambda}$ is an isomorphism and this proves condition (iii).

In case PEL the quasi-isogeny $\lambda:\bX\to\bX\dual$ induces a quasi-isogeny 
\[
\ol p_{\Lambda\dual}\;:=\;(\eta_\Lambda\dual)^{-1}\circ\lambda\circ\eta_{\Lambda\dual}^{-1}:\es X_{\Lambda\dual,\ol{S}}\;\longto\;(X_{\Lambda,\ol{S}})\dual\,.
\]
By the crystalline nature of the functor $\BD(\,.\,)_{S^\an}$ we obtain an isomorphism $\BD(\ol p_{\Lambda\dual})_{S^\an}:\BD(X_{\Lambda\dual})_{S^\an}\isoto\BD(X_\Lambda\dual)_{S^\an}\cong\BD(X_\Lambda)\dual_{S^\an}$ where the last isomorphism is induced by the pairing $\psi_\Lambda$. Since the alternating form $(\,.\,,\,.\,)$ on $V\otimes_{\BQ_p}K_0=\BD(\bX)_{K_0}$ coincides with 
\[
\BD(\bX)_{K_0}\otimes_{K_0}\BD(\bX)_{K_0} \;\xrightarrow{\es\BD(\lambda)\otimes\id\;}\;\BD(\bX\dual)_{K_0}\otimes_{K_0}\BD(\bX)_{K_0}\;\xrightarrow{\es\BD(\psi_\bX)_{K_0}\;}\;\BD(\BG_m)_{K_0}\;=\;K_0
\]
up to multiplication with a unit in $\BQ_p\mal$ and since $Fil^1(V\otimes_{\BQ_p}\CO_{S^\an})$ is totally isotropic for $(\,.\,,\,.\,)$ we see that
\begin{eqnarray*}
\lefteqn{\BD(\ol p_{\Lambda\dual})_{S^\an}(Fil^1 X_{\Lambda\dual})\otimes_{\CO_S}\CO_{S^\an} \es = \es \BD(\eta_\Lambda^{-1})\dual\BD(\lambda)\bigl(Fil^1(V\otimes_{\BQ_p}\CO_{S^\an})\bigr)}\hspace{3.5cm}\\[2mm]
& \subset & \BD(\eta_{\Lambda}^{-1})\dual\bigl(\{\,h\in V^{\SSC\lor}\otimes_{\BQ_p}\CO_{S^\an}:\;h\bigl(Fil^1(V\otimes_{\BQ_p}\CO_{S^\an})\bigr)=0\,\}\bigr)\\[2mm]
& = & \bigl\{\,h'\in\BD(X_\Lambda)\dual_{S^\an}:\;h'(Fil^1 X_\Lambda)=0\,\bigr\}\\[2mm]
& = & (Fil^1 X_\Lambda\dual)\otimes_{S}\CO_{S^\an}\,.
\end{eqnarray*}
This shows that $\ol p_{\Lambda\dual}$ lifts to a quasi-isogeny $p_{\Lambda\dual}:X_{\Lambda\dual}\to X_\Lambda\dual$ over $S$. By part 1 above the quasi-isogeny $p_{\Lambda\dual}\circ\eta_{\Lambda\dual,\Lambda}:X_{\Lambda}\to X_\Lambda\dual$ induces an alternating form on $\wt V\cong V_p X_\Lambda$ as in \eqref{Eq7.2} which equals $\tilde\ell\cdot(\,.\,,\,.\,)^{^\sim}$ for a constant $\tilde\ell\in\BQ_p\mal$. This $\tilde\ell$ is also the constant from Definition~\ref{Def7.2}($\wt{\rm v}$) of the multichain $\{\wt\Lambda\}_{\Lambda\in\CL}$. Then
\begin{eqnarray*}
\TS V_p(p_{\Lambda\dual})(T_pX_{\Lambda\dual}) & = & \TS V_p(\frac{1}{\tilde\ell}\,p_{\Lambda\dual}\circ\eta_{\Lambda\dual,\Lambda})\delta_\Lambda^{-1}(\tilde\ell\cdot\wt{\Lambda\dual})\\[2mm]
& = & \frac{1}{\tilde\ell} \,V_p(p_{\Lambda\dual}\circ\eta_{\Lambda\dual,\Lambda})\delta_\Lambda^{-1}(\wt\Lambda)\dual\\[2mm]
& = & \frac{1}{\tilde\ell}\, V_p(p_{\Lambda\dual}\circ\eta_{\Lambda\dual,\Lambda})\delta_\Lambda^{-1}\bigl(\{\,\tilde v\in\wt V:\;(\tilde v,\tilde w)^{^\sim}\in\BZ_p\text{ for all }\tilde w\in\wt\Lambda\,\}\bigr)\\[2mm]
& = & T_pX_\Lambda\dual\,.
\end{eqnarray*}
Therefore $p_{\Lambda\dual}:X_{\Lambda\dual}\to X_\Lambda\dual$ is an isomorphism and this establishes \cite[Definition 3.21(v)]{RZ}. Hence $(X_\Lambda,i_\Lambda,\eta_\Lambda)$ is a point in $\breve\CG(S)$ such that the induced morphism $S^\an\to\breve\CG^\an\xrightarrow{\;\breve\pi\,}\breve\CF^0_b$ coincides with the morphism constructed in part 2 above. This shows that $\mu\in\breve\CF^0_b$ lies in the image of $\breve\pi$ and proves \ref{ThmPEL_1}.

\bigskip
\noindent
\ref{ThmPEL_2} As in the proof of Theorem~\ref{ThmC}\ref{ThmC_2} the rational Tate module $V_p X_\Lambda$ descends to a local system $\CV$ of $\BQ_p$-vector spaces on the space $\breve\CF^0_b$. Consider a connected component $Y$ of $\breve\CF^0_b$ and a geometric base point $\bar\mu\in Y$ which lies over an $L$-valued point $\mu\in Y(L)$ for a finite field extension $L/\breve E$. By Proposition~\ref{Prop1.2c} the local system $\CV$ on $Y$ corresponds to a representation $\rho':\pi_1^\et(Y,\bar\mu)\to\GL(\CV_{\bar\mu})\cong\GL(\wt V)$ where we identify $V_p X_\Lambda\cong\wt V$ as in part 1 of (a) above. There is even an isomorphism of $\Gal(L^\alg/L)$-representations
\begin{equation}\label{Eq7.7}
\CV_\mu\;=\; V_pX_{\Lambda,\CO_L}\;\cong\;V_\cris\bigl(\BD(X_{\Lambda,\CO_L}),Fil^1 X_{\Lambda,\CO_L}\bigr)\;\cong\;V_\cris\bigl(\ulD_{b,\mu}(V)\bigr)\,.
\end{equation}
In order to prove that $\CV$ induces a tensor functor and \eqref{Eq7.7} is a tensor isomorphism it suffices by Corollary~\ref{Cor1.7} to show that $\rho'$ factors through $\wt G(\BQ_p)\subset\GL(\wt V)$, which in turn follows from the explicit description of $\wt G$ given in \eqref{WTG}. Indeed, since the universal $X_\Lambda$ over $\breve\CG$ carries an action of $\CO_B$ the image of $\pi_1^\et(Y,\bar\mu)$ in $\GL(\wt V)$ commutes with $B$ and this already proves the claim in case EL. 
In case PEL, diagram \eqref{Eq7.2} shows that the action of $\gamma\in\pi_1^\et(Y,\bar\mu)$ on $\wt V$ satisfies $(\gamma\cdot\tilde v,\gamma\cdot\tilde w)^{^\sim}=\chi_{\rm cyc}(\gamma)\cdot(\tilde v,\tilde w)^{^\sim}$ for all $\tilde v,\tilde w\in\wt V$, where $\chi_{\rm cyc}:\pi_1^\et(Y,\bar\mu)\to\Gal(\breve E^\alg/\breve E)\xrightarrow{\;\chi_{_{\rm cyc}}\,}\BZ_p\mal$ is the cyclotomic character of the base field $\breve E$. In particular, $\rho'$ factors through a morphism $\rho':\pi_1^\et(Y,\bar\mu)\to\wt G(\BQ_p)$ with $\tilde c\circ\rho'=\chi_{\rm cyc}$,
and by Corollary~\ref{Cor1.7} this induces the desired tensor functor $\ul\CV:\Rep_{\BQ_p}\!\!\!G \to\PLoc_{Y}$ with $\ul\CV(V)=\CV$ satisfying \eqref{EqRZb}. Moreover, \eqref{Eq7.7} is compatible with the action of $B$ on both sides and, in the PEL case, the alternating pairings $V_p\psi_\Lambda$ and $V_\cris\bigl(\BD(\psi_\bX)_{K_0}\bigr)$; see \eqref{Eq7.2}. Hence \eqref{Eq7.7} induces an isomorphism of tensor functors on $\Rep_{\BQ_p}\!\!\!G$ and also (\ref{EqRZ}) of Conjecture~\ref{ConjRZ} holds at the fixed $\mu$. But the construction is independent of $\mu$. So it holds on all connected components $Y$ of $\breve\CF^0_b$ and this establishes \ref{ThmPEL_2} and Remark~\ref{Rem7.5}(a).

\medskip\noindent
\ref{ThmPEL_3} The argument of Theorem~\ref{ThmC}\ref{ThmC_3} applies literally to prove the first part of (c). For the statement about $\breve\CG^\an$ let $\wt K\subset\wt G(\BQ_p)$ be the stabilizer of the (polarized) multichain $(\wt\Lambda)_{\Lambda\in\CL}$ of type $(\CL)$ in $\wt V$ constructed in part 1 above. The construction of parts 2 and 3 yields an isomorphism $\breve\CG^\an\isoto\breve\CG^\an_{\wt K}\cong\breve\CE_{\wt K}$ where $\breve\CE_{\wt K}$ is the space of (polarized) multichains of type $(\CL)$ in $\CV$.

\medskip\noindent
\ref{ThmC_4} is proved by the same argument as Theorem~\ref{ThmC}\ref{ThmC_4}.
\end{proof}

%
%

\section{When does weakly admissible imply admissible ?} \label{SectWAImpliesAdm}
\setcounter{equation}{0}

In this section we approach the question for which groups $G$ and which elements $b\in G(K_0)$ as in \eqref{Eq2.2} ``weakly admissible implies admissible'', that is $\breve\CF_b^{wa}=\breve\CF_b^0$. Since we do not yet have a good group theoretic understanding of the phenomenon we only treat the case $G=\GL_n$ comprehensively. We begin with the following observations

\begin{lemma}\label{Lemma8.1}
Let $G=G_1\times G_2$ be a product of two reductive groups over $\BQ_p$ and $b=(b_1,b_2)$ such that each $b_i\in G_i(K_0)$ satisfies \eqref{Eq2.2} for the same integer $n$. If $\breve\CF_{G_1,b_1}^{wa}\neq\breve\CF_{G_1,b_1}^0$ and $\breve\CF_{G_2,b_2}^{wa}\ne\emptyset$ then $\breve\CF_{G,b}^{wa}\ne\breve\CF_{G,b}^0$.
\end{lemma}

\begin{proof}
Let $\mu_1\in\breve\CF_{G_1,b_1}^{wa}\setminus\breve\CF_{G_1,b_1}^0$ be an $\LL_1$-valued point and $\mu_2\in\breve\CF_{G_2,b_2}^{wa}(\LL_2)$ for complete field extensions $\LL_1/\breve E_1$ and $\LL_2/\breve E_2$ with $[\LL_2:\breve E_2]<\infty$. Let $\LL$ be some compositum of $\LL_1$ and $\LL_2$ in some algebraic closure of $\LL_1$ and set $\mu=(\mu_1,\mu_2):\BG_m\to G$. Let $\rho_i:G_i\into\GL(V_i)$ be representations satisfying \eqref{Eq2.2} for $n$. Then $\rho:G\into\GL(V_1\oplus V_2)$ likewise satisfies \eqref{Eq2.2} for $n$. Consider the filtered isocrystals $\ulD_i=\ulD_{b_i,\mu_i}(V_i)$ and 
\[
\ulD\;=\;\ulD_1\oplus\ulD_2\;=\;\ulD_{b,\mu}(V_1\oplus V_2)
\]
over $\LL$ which all three are weakly admissible. Then
\[
\bM(\ulD)\;=\;\bM(\ulD_1)\oplus\bM(\ulD_2)\;\not\cong\;\bM(0,1)\otimes_{\BQ_p}(V_1\oplus V_2)
\]
because $\bM(\ulD_1)\not\cong\bM(0,1)\otimes_{\BQ_p}V_1$. This shows that $\mu\in\breve\CF_{G,b}^{wa}\setminus\breve\CF_{G,b}^0$ as desired.
\end{proof}

The next lemma shows that the question is invariant under dualization.

\begin{lemma}\label{Lemma8.1'}
Let $G_i=\GL_n$ and $b_i\in\GL_n(K_0)$ for $i=1,2$ such that the isocrystals $\ulD_i=(K_0^n\,,\,b_i\!\cdot\!\phi)$ satisfy \eqref{Eq2.2}. Assume that
\begin{enumerate}
\item 
$\ulD_2\cong\ulD_1\otimes\BOne(r)$ or
\item 
$\ulD_2\cong\Hom\bigl(\ulD_1,\BOne(r)\bigr)$ 
\end{enumerate}
holds for some $r\in\BZ$ (see Lemma~\ref{Lemma4.2'}). Then $\breve\CF_{G_1,b_1}^{wa}=\breve\CF_{G_1,b_1}^0$ if and only if $\breve\CF_{G_2,b_2}^{wa}=\breve\CF_{G_2,b_2}^0$.
\end{lemma}

\begin{proof}
By symmetry it suffices to prove one direction. So assume that $\breve\CF_{G_2,b_2}^{wa}=\breve\CF_{G_2,b_2}^0$ and let $\mu_1\in\breve\CF_{G_1,b_1}^{wa}$. Since the filtered isocrystal $\ulD_1=(K_0^n,b_1\!\cdot\!\phi,Fil^\bullet_{\mu_1})$ is weakly admissible, also the filtered isocrystal $\ulD_2:=\ulD_1\otimes\BOne(r)$ in case (a), respectively $\ulD_2:=\Hom\bigl(\ulD_1,\BOne(r)\bigr)$ in case (b) is weakly admissible and corresponds to a point $\mu_2\in\breve\CF_{G_2,b_2}^{wa}$. Then our assumption implies $\bM(\ulD_2)\cong\bM(0,1)^n$. Lemma~\ref{Lemma4.2'} and Proposition~\ref{Prop4.10} yield in case 
\begin{enumerate}
\item 
$\bM(\ulD_1)\cong\bM(\ulD_2)\cong\bM(0,1)^n$, respectively 
\item 
$\bM(\ulD_1)\cong\bM\bigl(\ulD_2\dual\otimes\BOne(r)\bigr)\cong\bM(\ulD_2\dual)\cong\bigr(\bM(0,1)^n\bigl)\dual\cong\bM(0,1)^n$,
\end{enumerate}
and therefore $\mu_1\in\breve\CF_{G_1,b_1}^0$ as desired.
\end{proof}

Next we state the main theorem of this section.

\begin{theorem}\label{Thm8.2}
Let $G=\GL_n$, let $b\in\GL_n(K_0)$ be such that the isocrystal $(K_0^n,\,b\!\cdot\!\phi)$ has Newton slopes in the interval $[0,1]$, and let $\{\mu_0\}$ be the conjugacy class for which the standard representation $\BQ_p^n$ of $G$ satisfies \eqref{Eq2.2}. Then the equality $\breve\CF_b^{wa}=\breve\CF_b^0$ holds if and only if the Newton slopes of $b$ are
\begin{enumerate}
\item \label{Thm8.2_1}
$(1^{(h_1)},\frac{1}{h}^{(h)},0^{(h_0)})$ for $h_1,h,h_0\in\BN_0$, or
\item \label{Thm8.2_2}
$(1^{(h_1)},\frac{h-1}{h}^{(h)},0^{(h_0)})$ for $h_1,h,h_0\in\BN_0$, or
\item \label{Thm8.2_3}
$(1^{(h_1)},\frac{1}{2}^{(4)},0^{(h_0)})$ for $h_1,h_0\in\BN_0$.
\end{enumerate}
In particular $\breve\CF_b^{wa}=\breve\CF_b^0$ for $n\le 4$. (Here $s^{(m)}$ means that the slope $s\in\BQ$ occurs with multiplicity $m\in\BN_0$. If $m=0$ this slope does not occur. The sum of all multiplicities is $n$.)
\end{theorem}

Before proving the theorem we derive the following

\begin{corollary}\label{Cor8.3}
Let $(G,b,\{\mu\})$ be arbitrary such that there is a faithful representation $\rho:G\into\GL(V)$ as in \eqref{Eq2.2} with $\dim_{\BQ_p}V\le 4$. Then $\breve\CF_b^{wa}=\breve\CF_b^0$.
\end{corollary}

\begin{proof}
This follows from the theorem and the fact that $\breve\CF_{G,b}^0=\breve\CF_{G,b}^{wa}\cap\breve\CF_{\GL(V),\rho(b)}^0$, see \eqref{EqFAdm} and Proposition~\ref{PropDOR}.
\end{proof}

\begin{remark} 
1. The case of Theorem~\ref{Thm8.2} where $h=0$ and $h_1=h_0$ was discussed by Rapoport and Zink \cite[5.51]{RZ}. In this case the flag variety is the Grassmannian $\Grass(h_0,2h_0)$, and $\breve\CF_b^{wa}$ is the big cell of subspaces transversal to the sub-isocrystal of slope zero. It is isomorphic to affine space $\BA^{h_0^2}$. The period morphism was constructed by Dwork (compare \cite{KatzDwork}) and, according to \cite[Proposition 5.52]{RZ}, is given component-wise as the $p$-adic logarithm. By the surjectivity of the latter, the image $\breve\CF_b^0$ of the period morphism equals $\breve\CF_b^{wa}$.

\medskip\noindent
2. Let us now assume that $h_1=h_0=0$ in Theorem~\ref{Thm8.2}.

Then case \ref{Thm8.2_2} (and its dual case \ref{Thm8.2_1}) correspond to the Lubin-Tate deformation space \cite{LT} and were studied by Gross and Hopkins \cite{HG1,HG2}. In these cases the flag variety $\breve\CF$ equals $\BP^{h-1}$, the space of hyperplanes (respectively lines) in $\BQ_p^{\oplus h}$, the Rapoport-Zink space $\breve\CG$ is $\coprod_\BZ\Spf W(\BF_p^\alg)\dbl t_1,\ldots,t_{h-1}\dbr$, and the period morphism has image $\breve\CF^0_b=\breve\CF_b^{wa}=(\BP^{h-1})^\an$ by \cite{HG1,HG2}; see also \cite[5.50]{RZ}. In this case our Theorem~\ref{ThmC} was proved by \cite[\S\S1,7]{dJ} who also determines the image of the representation $\rho:\pi_1^\et(\breve\CF_b^1,\bar\mu)\to\GL_h(\BQ_p)$ corresponding to the tensor functor $\ul\CV$.

The third alternative, case \ref{Thm8.2_3} in Theorem~\ref{Thm8.2} comes as a surprise and the author has no interpretation for it so far.

\medskip\noindent
3. Among examples for Corollary~\ref{Cor8.3} is the EL case where $B$ is the quaternion algebra over $\BQ_p$, where $V$ is a free $B$-module of rank one, and $G=(B^{opp})\mal$; see Section~\ref{SectPELPeriodM}. We present $B$ as
\[
B=\BQ_{p^2}[\Pi]\quad\text{with}\quad \Pi^2=p,\text{ and }\Pi x=\phi(x)\Pi
\]
and let $b=p^{-1}\Pi\in G(K_0)$. This case was studied by Drinfeld~\cite{Drinfeld76}. Here the flag variety is $\BP^1$, and $\breve\CF_b^{wa}$ is the base change $\Omega\times_{\BQ_p}K_0$ of Drinfeld's upper half-plane $\Omega=\BP^1\setminus\BP^1(\BQ_p)$; see \cite[1.44--1.46]{RZ}. The Rapoport-Zink space $\breve\CG$ is the disjoint union indexed by $\BZ$ of Deligne's integral model of $\Omega\times_{\BQ_p}K_0$. In particular each connected component of $\breve\CG^\an$ maps isomorphically onto $\breve\CF_b^{wa}$ and hence $\breve\CF_b^0=\breve\CF_b^{wa}$; see \cite[3.45--3.77 and 5.48--5.49]{RZ}. This case directly generalizes to the situation where $B$ is a central division algebra of dimension $e$ and Hasse invariant $\tfrac{1}{e}$.
\end{remark}

The rest of this section is devoted to the

\begin{proof}[Proof of Theorem~\ref{Thm8.2}]
We first prove that the conditions (a), (b), and (c) are sufficient. Let $\mu\in\breve\CF^\an\setminus\breve\CF^0_b$. We must show that $\mu\notin\breve\CF_b^{wa}$. We set $\bM_\mu=\bM(\ulD_{b,\mu}(\BQ_p^n))$ and $\bD=\bD(\ulD_{b,\mu}(\BQ_p^n))$ (Definition~\ref{Def1.12b}), and consider the decompositions of $\bD$ and $\bM_\mu$ from Theorem~\ref{Thm1.4}. We write $\bM_\mu=\bM_+\oplus\bM_0\oplus\bM_-$ such that the weights of all components of $\bM_+$ (respectively $\bM_0$, respectively $\bM_-$) are positive (respectively zero, respectively negative). Note that $\bM_-\ne(0)$ since $\bM_\mu\not\cong\bM(0,1)^n$ and $\deg\bM_\mu=0$ by Theorem~\ref{Thm1.13} and Lemma~\ref{LemmaKottwitzCond}. Likewise we write $\bD=\bM(1,1)^{\oplus h_1}\oplus\bD_+\oplus\bM(0,1)^{\oplus h_0}$ such that the weights of all components of $\bD_+$ lie strictly between $0$ and $1$. Then the $\phi$-module $t^{-1}\bD=\bM(0,1)^{\oplus h_1}\oplus t^{-1}\bD_+\oplus\bM(-1,1)^{\oplus h_0}$. We consider the inclusions $\bD\into\bM_\mu\into t^{-1}\bD$ whose composition is the natural inclusion $\bD\subset t^{-1}\bD$, and the induced morphisms $\bM(0,1)^{\oplus h_0}\to\bM_0\to\bM(0,1)^{\oplus h_1}$ whose composition is the zero map. Since $\End_\phi\bigl(\bM(0,1)\bigr)=\BQ_p$ by Proposition~\ref{Prop1.6}\ref{Prop1.6_c}, we can decompose the $\phi$-modules $\bM_\mu\,=\,\bM_+\oplus\bM_0^1\oplus\bM_0^2\oplus\bM_0^3\oplus\bM_-$, and $\bD\,=\,\bM(1,1)^{\oplus\tilde h_1}\oplus\bM(1,1)^{\oplus h'_1}\oplus\bD_+\oplus\bM(0,1)^{\oplus h'_0}\oplus\bM(0,1)^{\oplus\tilde h_0}$ such that with respect to appropriate bases of all these $\phi$-modules (and the induced basis for $t^{-1}\bD$) the inclusions are given by block matrices
\[
\begin{array}{c|ccccc}  
\bM_\mu\into t^{-1}\bD&\bM_+&\bM_0^1&\bM_0^2&\bM_0^3&\bM_-\\
\hline\\[-3mm]
\bM(0,1)^{\oplus\tilde h_1}&A_1&0&0&0&0\\[1mm]
\bM(0,1)^{\oplus h'_1}&A_2&\Id&0&0&0\\[1mm]
t^{-1}\bD_+&D&E_1&E_2&E_3&F\\[1mm]
\bM(-1,1)^{\oplus h'_0}&G_1&H_{11}&H_{12}&H_{13}&I_1\\[1mm]
\bM(-1,1)^{\oplus\tilde h_0}&G_2&H_{21}&H_{22}&H_{23}&I_2
\end{array}
\]
and
\[
\begin{array}{c|ccccc}
\bD\into\bM_\mu&\bM(1,1)^{\oplus\tilde h_1}&\bM(1,1)^{\oplus h'_1}&\bD_+&\bM(0,1)^{\oplus h'_0}&\bM(0,1)^{\oplus\tilde h_0}\\
\hline
\bM_+&a_1&a_2&c&0&0\\[1mm]
\bM_0^1&d_{11}&d_{12}&e_1&0&0\\[1mm]
\bM_0^2&d_{21}&d_{22}&e_2&0&0\\[1mm]
\bM_0^3&d_{31}&d_{32}&e_3&\Id&0\\[1mm]
\bM_-&f_1&f_2&g&i_1&i_2
\end{array}
\]
We put $r_+:=\rk\bM_+\,,\,r_0^i:=\rk\bM_0^i\,,\,r_-:=\rk\bM_-$, and $h:=\rk\bD_+$. Here and in the sequel superscripts are indices and do not mean powers. We must have 
\begin{equation}\label{Eq8.1}
r_0^1\;=\;h_1'\,,\qquad r_0^3\;=\;h_0'\,,\qquad\tilde h_0\;\le\; r_-\;, \qquad\text{and}\qquad r_0^2+r_0^3+r_-\;\le\; h+h_0
\end{equation}
since the maps $\bM(0,1)^{\oplus\tilde h_0}\into\bM_-$ and $\bM_0^2\oplus\bM_0^3\oplus\bM_-\into t^{-1}\bD_+\oplus\bM(-1,1)^{\oplus h_0}$ are injective. The inclusion $\bD\into t^{-1}\bD$ is represented by the product of the two big block matrices which hence equals
\[
t\Id_n\es=\es
\left(\begin{array}{ccccc}A_1a_1&A_1a_2&A_1c&0&0\\
A_2a_1+d_{11}&A_2a_2+d_{12}&A_2c+e_1&0&0\\
 &  &  & E_3+Fi_1&Fi_2\\
 &  &  & H_{13}+I_1i_1&I_1i_2\\
 &  &  & H_{23}+I_2i_1&I_2i_2
\end{array}\right)
\]
This implies $A_1c=0\,,\,A_2c+e_1=0\,,\,E_3+Fi_1=0$\,, and $Fi_2=0$. Considering the entries of these matrices as elements of the fraction field $\CQ$ of the integral domain $\wt\bB^\dagger_\rig$ we have $\dim_\CQ\ker A_1\ge\rk c$ and $\dim_\CQ\ker F\ge\rk i_2=\tilde h_0$. Thus
\begin{equation}\label{Eq8.2}
\tilde h_1\;=\;\rk A_1\;\le\;r_+-\rk c\qquad\text{and}\qquad r_-\;=\;\rk F+\dim_\CQ\ker F\;\ge\;\tilde h_0+\rk F\,.
\end{equation}
Now we distinguish several cases.

\medskip
\noindent
\ul{Case $c=0$ :} This implies $e_1=0$ and we obtain inclusions
\[
\bD'\;:=\;\bD_+\oplus\bM(0,1)^{\oplus h_0}\;\into\;\bM_0^2\oplus \bM_0^3\oplus\bM_-\;\into\;t^{-1}\bD_+\oplus\bM(-1,1)^{\oplus h_0}
\]
Therefore $\bD'\subset\bD$ comes from a sub-isocrystal $\ulD'\subset\ulD$ with $\bM(\ulD')\supset\bM_0^2\oplus \bM_0^3\oplus\bM_-$. This implies $t_N(\ulD')-t_H(\ulD')=\deg\bM(\ulD')\le\deg(\bM_0^2\oplus \bM_0^3\oplus\bM_-)<0$ by Theorem~\ref{Thm1.13} and \cite[Lemma 3.4.2]{Kedlaya}. So $\mu$ is not weakly admissible, $\mu\notin\breve\CF_b^{wa}$.

\medskip
\noindent
\ul{Case $F=0$ :} This implies $E_3=0$ and we obtain inclusions
\[
\bD'\;:=\;\bM(0,1)^{\oplus h_0}\;\into\;\bM_0^3\oplus\bM_-\;\into\;\bM(-1,1)^{\oplus h_0}
\]
Therefore $\bD'\subset\bD$ comes from a sub-isocrystal $\ulD'\subset\ulD$ with $\bM(\ulD')\supset\bM_0^3\oplus\bM_-$. This again implies $t_N(\ulD')-t_H(\ulD')<0$ by Theorem~\ref{Thm1.13} and \cite[Lemma 3.4.2]{Kedlaya}, and so $\mu$ is not weakly admissible, $\mu\notin\breve\CF_b^{wa}$.

\medskip
\noindent
\ul{Case $c\ne0$ and $\rk c=h$ :} Then \eqref{Eq8.2} implies $\tilde h_1+h\le r_+$ and 
\[
\tilde h_0\;=\;n-h_1'-\tilde h_1-h-h_0'\;\ge\; n-r_0^1-r_+-r_0^3\;=\;r_0^2+r_-\;\ge\;r_-\;\ge\;\tilde h_0\,,
\]
implies $r_-=\tilde h_0$. Since $t^{-1}I_2i_2=\Id_{r_-}$ we have $i_2\in\GL_{r_-}(\CQ)$. Thus $Fi_2=0$ yields $F=0$ in which case we just saw that $\mu\notin\breve\CF_b^{wa}$.

\medskip
\noindent
\ul{Case $F\ne0$ and $\rk F=h$ :} Then \eqref{Eq8.2} implies $\tilde h_0+h\le r_-$ and 
\[
h_0+h\;=\;h_0'+\tilde h_0+h\;\le\; r_0^3+r_-\;\le\; r_0^2+r_0^3+r_-\;\le\; h_0+h\,,
\]
implies $h_0+h=r_0^2+r_0^3+r_-$ and $\tilde h_1=r_+$. Since $t^{-1}A_1a_1=\Id_{r_+}$ we have $A_1\in\GL_{r_+}(\CQ)$ and hence $A_1c=0$ yields $c=0$. This case was treated above.

\medskip
\noindent
\ul{Case $c\ne0$ and $F\ne0$ :} To treat this case we use the special form $\bD_+$ has under the conditions (a), (b), and (c). 

Under condition (a) $\bD_+\cong\bM(1,h)$. The $\phi$-module which is the image of $c:\bD_+\to\bM_+$ is isomorphic to $\bigoplus_i\bM(s_i,n_i)$. Since it is a quotient of $\bD_+$ and contained in $\bM_+$ we have $\frac{1}{h}\ge\frac{s_i}{n_i}>0$, whence $n_i\ge s_ih\ge h$ for all $i$. This means $\rk c=h$ and in fact, we are in a case treated above.

Under condition (b) $\bD_+\cong\bM(h-1,h)$. The image of $F:\bM_-\to t^{-1}\bD_+$ is isomorphic to $\bigoplus_i\bM(s_i',n_i')$. Since it is a quotient of $\bM_-$ and contained in $t^{-1}\bD_+\cong\bM(-1,h)$ we have $0>\frac{s_i'}{n_i'}\ge\frac{-1}{h}$, whence $n_i'\ge -s_i'h\ge h$ for all $i$. This means $\rk F=h$ and we are again in a case treated above.

Under condition (c) $\bD_+\cong\bM(1,2)^{\oplus 2}$. Again the image of $c:\bD_+\to\bM_+$ is isomorphic to $\bigoplus_i\bM(s_i,n_i)$ with $\frac{1}{2}\ge\frac{s_i}{n_i}>0$, whence $n_i\ge 2s_i\ge 2$ and $\rk c\ge2$. Furthermore the image of $F:\bM_-\to t^{-1}\bD_+$ is isomorphic to $\bigoplus_i\bM(s_i',n_i')$ with $0>\frac{s_i'}{n_i'}\ge\frac{-1}{2}$, whence $n_i'\ge -2s_i'\ge 2$ and $\rk F\ge2$. Then
\[
n\;=\;h_0+4+h_1\;\le\;h'_0+\tilde h_0+\rk F+\rk c+\tilde h_1+h'_1\;\le\;r_0^3+r_-+r_++r_0^1\;=\;n-r_0^2\;\le\;n
\]
yields $r_0^2=0\,,\,\rk c=2=\rk F$\, and $r_-=\tilde h_0+\rk F=\tilde h_0+2$. This implies $n'_i\le 2$, whence $n'_i=2$, $s'_i=-1$ and $\im(F)\cong\bM(-1,2)$. We may decompose $t^{-1}\bD_+=t^{-1}\bD_+^1\oplus t^{-1}\bD_+^2$ with $\im(F)= t^{-1}\bD_+^2$. With respect to a basis suited to this decomposition we write $F={0\choose F_2}\,,\,E_3={E_{31}\choose E_{32}}\,,c=(c_1|c_2)$\, and $e_1=(e_{11}|e_{12})$. Since the zero map 
\[
\bM_-\xrightarrow{\TS{0\choose F_2}} t^{-1}\bD_+\xrightarrow{(c_1|c_2)}t^{-1}\bM_+
\]
equals $c_2F_2$ and $\rk F_2=2$ we find $c_2=0$. Then ${E_{31}\choose E_{32}}+{0\choose F_2}i_1=0$ implies $E_{31}=0$ and $A_2(c_1|c_2)+(e_{11}|e_{12})=0$ implies $e_{12}=0$. We obtain inclusions
\[
\bD'\;:=\;\bD_+^2\oplus\bM(0,1)^{\oplus h_0}\;\into\;\bM_0^3\oplus\bM_-\;\into\;t^{-1}\bD_+^2\oplus\bM(-1,1)^{\oplus h_0}
\]
The $\phi$-submodule $\bD_+^2$ is inherited from a  sub-isocrystal \mbox{$\ulD_+^2:=(K_0^n,b\!\cdot\!\phi)\cap\bD_+^2$} of the isocrystal $\ulD_+=(K_0^n\,,\,b\!\cdot\!\phi)\cap\bD_+$ because all slopes of $\im(F)$ and $t^{-1}\bD_+$ are equal, and $\End_\phi(\bM(-1,2))$ equals the central division algebra over $\BQ_p$ of dimension $4$ and Hasse invariant $\frac{1}{2}$ by Proposition~\ref{Prop1.6}. Therefore $\bD'\subset\bD$ comes from a sub-isocrystal $\ulD'\subset\ulD$ with $\bM(\ulD')\supset\bM_0^3\oplus\bM_-$. This implies $t_N(\ulD')-t_H(\ulD')<0$ by Theorem~\ref{Thm1.13} and \cite[Lemma 3.4.2]{Kedlaya}, and so $\mu$ is not weakly admissible, $\mu\notin\breve\CF_b^{wa}$.

Altogether we have proved that under conditions (a), (b), and (c) the equality $\breve\CF_b^{wa}=\breve\CF_b^0$ holds, that is ``weakly admissible implies admissible''.

\bigskip

To prove the converse we assume that conditions (a), (b), or (c) are not satisfied. By Lemma~\ref{Lemma8.1} we only need to find a direct summand of the isocrystal $(K_0^n\,,\,b\!\cdot\!\phi)$ for which ``weakly admissible does not imply admissible''. Considering simple sub-isocrystals it hence suffices to treat the cases where $\bD=\bD(\ulD)$ has the following form
\begin{enumerate}
\def\labelenumi{\theenumi} 
\def\theenumi{\arabic{enumi}.} 
\item\label{8.2i} $\bD\cong\bM(c,h)$ and $t^{-1}\bD\cong\bM(-d,h)$ for $c,h\in\BN_{>0}$ relatively prime, $d=h-c$, with $2\le c,d\le h-2$ and $h\ge 5$.
\item \label{8.2ii} $\bD\cong\bM(1,h_1)\oplus\bM(1,h_2)$ with $2\le h_1\le h_2$ and $2<h_2$.
\item \label{8.2iii} $\bD\cong\bM(h_1-1,h_1)\oplus\bM(h_2-1,h_2)$ with $2\le h_1\le h_2$ and $2<h_2$.
\item \label{8.2iv} $\bD\cong\bM(1,h_1)\oplus\bM(h_2-1,h_2)$ with $3\le h_1,h_2$.
\item \label{8.2v} $\bD\cong\bM(1,2)^{\oplus 3}$.
\end{enumerate}

Since the isocrystals in cases \ref{8.2ii} and \ref{8.2iii} are dual to each other, case \ref{8.2iii} follows from case \ref{8.2ii} by Lemma~\ref{Lemma8.1'}. Also note that a special instance of case \ref{8.2i} was discussed in Example~\ref{Example6.4}. However, the beginning of our argument treats cases \ref{8.2i}, \ref{8.2ii}, \ref{8.2iv}, and \ref{8.2v} simultaneously. We set $h:=\rk\bD,\,c:=t_N(\ulD)=\deg\bD$, and $d=h-c$. In case \ref{8.2i} the conditions on $c,d,h$ imply $cd>h$ and we define $e=\lfloor\frac{cd}{h}\rfloor\ge1$ such that $\frac{-e}{c}\ge\frac{-d}{h}$. In the other three cases we set $e=1$. By Proposition~\ref{Prop1.6} there exists $0\ne A\in\Hom_\phi\bigl(\bM(-e,c),t^{-1}\bD\bigr)$, such that depending on the case
\begin{description}
\item[{\rm Case \ref{8.2i}}] the map $A:\bM(-e,c)\to t^{-1}\bD$ is arbitrary but non-zero. 
\item[{\rm Case \ref{8.2ii}}] the map $A:\bM(-e,c)=\bM(-1,2)\to t^{-1}\bD\cong\bM(1-h_1,h_1)\oplus\bM(1-h_2,h_2)$ induces non-zero maps $A_1$ and $A_2$ into the first, resp.\ second summand of $t^{-1}\bD$. If $h_1=h_2>2$ we require in addition that $A_1,A_2$ are linearly independent over the central division algebra $\End_\phi\bigl(\bM(1-h_1,h_1)\bigr)$ of dimension $(h_1)^2$ over $\BQ_p$. This is possible since $\dim_{\BQ_p}\Hom_\phi\bigl(\bM(-1,2),\bM(1-h_1,h_1)\bigr)=\infty$ in this case by Proposition~\ref{Prop1.6}.
\item[{\rm Case \ref{8.2iv}}] the map $A:\bM(-e,c)=\bM(-1,h_2)\to t^{-1}\bD\cong\bM(1-h_1,h_1)\oplus\bM(-1,h_2)$ induces a non-zero map into the first summand of $t^{-1}\bD$ and is the identity onto the second summand.
\item[{\rm Case \ref{8.2v}}] the map $A:\bM(-e,c)=\bM(-1,3)\to t^{-1}\bD\cong\bM(-1,2)^{\oplus3}$ induces maps $A_i$ into the $i$-th summand of $t^{-1}\bD$ which are linearly independent over the central division algebra $\End_\phi\bigl(\bM(-1,2)\bigr)$ of dimension $4$ over $\BQ_p$. Again this is possible since $\dim_{\BQ_p}\Hom_\phi\bigl(\bM(-1,3),\bM(-1,2)\bigr)=\infty$ by Proposition~\ref{Prop1.6}.
\end{description}
Since both $\bM(-e,c)$ and $t^{-1}\bD$ are represented over $\wt\bB^{]0,1]}(\BC_p)$, we may by Proposition~\ref{Prop1.9b} choose a $K_0$-basis of $\ulD$ and consider $A$ as an $(h\times c)$-matrix with entries in $\wt\bB^{]0,1]}$. There is a $\mu$ defined over $\BC_p$ in the conjugacy class $\{\mu_0\}$, such that the $c$-dimensional space $Fil^1_\mu \BC_p^n$ contains the columns of the matrix $\theta_{t^{-1}\bD}(A)$. By Definition~\ref{Def1.12b} (compare Example~\ref{Example6.4}) the non-zero map $\bM(-e,c)\to t^{-1}\bD$ factors through $\bM_\mu$ for this $\mu$, and therefore $\bM_\mu\not\cong\bM(0,1)^{\dim D}$ and $\mu\notin\breve\CF_b^0$. It remains to show that $\mu\in\breve\CF_b^{wa}$. We distinguish the cases

\smallskip
\ref{8.2i} In this case $\ulD$ is simple and hence $\breve\CF_b^{wa}=\breve\CF^\an$. This implies $\mu\in\breve\CF_b^{wa}$.

\smallskip
\ref{8.2ii} Assume that there is a sub-isocrystal $(0)\ne\ulD'\subsetneq\ulD$ that contradicts weak admissibility. Then $t^{-1}\bD':=t^{-1}\ulD'\otimes_{K_0}\wt\bB^\dagger_\rig$ is the kernel of a non-zero map $E:t^{-1}\bD\to\bM(1-h_i,h_i)$ for $i=1$ or $i=2$. If $\mu\notin\breve\CF_b^{wa}$ then 
\[
1\;=\;t_N(\ulD')\;<\;t_H(\ulD')\;=\;\dim_{\BC_p}(D'_{\BC_p}\cap Fil^1_\mu D_{\BC_p})\;\le\; \dim_{\BC_p}Fil^1_\mu D_{\BC_p}\;=\;2\,,
\]
and this would imply $Fil^1_\mu D_{\BC_p}\subset D'_{\BC_p}$, whence $\theta_{\bM(1-h_i,h_i)}(EA)=0$. Thus $EA=T\cdot B$ for some map $B\in\Hom_\phi\bigl(\bM(-1,2)\,,\,\bM(1,h_i)\bigr))$ by Corollary~\ref{Cor1.6b}. Then Proposition~\ref{Prop1.6} shows that $B=0$, whence $EA=0$. However, this is not the case by our assumptions on $A$, since for $h_1\ne h_2$ the map $E$ is the projection onto the summand $\bM(1-h_i,h_i)$ of $t^{-1}\bD$. Likewise for $h_1=h_2$ we have $E=(E_1|E_2):t^{-1}\bD\to\bM(1-h_1,h_1)$ with $E_i\in\End_\phi\bigl(\bM(1-h_1,h_1)\bigr)$ and $E_1A_1+E_2A_2=0$ contradicts the linear independence of $A_1,A_2$ over $\End_\phi\bigl(\bM(1-h_1,h_1)\bigr)$. Therefore $\mu\in\breve\CF_b^{wa}$.

\smallskip
\ref{8.2iv} In this case $A={A_1\choose \Id_{h_2}}$ with $0\ne A_1\in\Hom_\phi\bigl(\bM(-1,h_2)\,,\,\bM(1-h_1,h_1)\bigr)$, and $\theta_{t^{-1}\bD}(A)={\theta(A_1)\choose\Id_{h_2}}\ne{0\choose\Id_{h_2}}$ because $\theta_{\bM(1-h_1,h_1)}(A_1)=0$ would by Corollary~\ref{Cor1.6b} imply that $A_1=T\cdot B$ for $0\ne B\in\Hom_\phi\bigl(\bM(-1,h_2)\,,\,\bM(1,h_1)\bigr)$ which is impossible by Proposition~\ref{Prop1.6}. In particular, $\theta_{t^{-1}\bD}(A)$ has rank $c$ and $Fil^1_\mu D_{\BC_p}$ equals $\im\bigl(\theta_{t^{-1}\bD}(A)\bigr)$.

There are exactly two sub-isocrystals $(0)\ne\ulD'\subsetneq\ulD$ of $\ulD$. For the first one we have $\bD':=\ulD'\otimes_{K_0}\wt\bB^\dagger_\rig=\bM(1,h_1)$ and consequently $t_N(\ulD')=1$. Furthermore, we have $D'_{\BC_p}\cap Fil^1_\mu D_{\BC_p}=\{{v\choose 0}\in\im{\theta(A_1)\choose\Id_{h_2}}\}=(0)$ and $t_H(\ulD')=0$.
For the second sub-isocrystal we have $t^{-1}\bD'=\bM(-1,h_2)$ and $t_N(\ulD')=h_2-1$. We claim that $t_H(\ulD')\le t_N(\ulD')$. Namely, if 
\[
h_2-1\;=\;t_N(\ulD')\;<\;t_H(\ulD')\;=\;\dim_{\BC_p}(D'_{\BC_p}\cap Fil^1_\mu D_{\BC_p})\;\le\;\dim_{\BC_p}Fil^1_\mu D_{\BC_p}\;=\;h_2
\]
then $\im\bigl(\theta_{t^{-1}\bD}(A)\bigr)=Fil^1_\mu D_{\BC_p}=D'_{\BC_p}=\im{0\choose\Id_{h_2}}$ in contradiction to $\theta_{\bM(1-h_1,h_1)}(A_1)\ne0$. This proves that $\mu$ is weakly admissible.

\smallskip
\ref{8.2v} Here the sub-isocrystals $(0)\ne\ulD'\subsetneq\ulD$ are of two kinds. For the first kind $t^{-1}\bD'$ has rank $4$ and is the kernel of a non-zero map $E=(E_1|E_2|E_3):t^{-1}\bD\to\bM(-1,2)$ for $E_i\in\End_\phi\bigl(\bM(-1,2)\bigr)$. Then $t_N(\ulD')=2$ and we claim that $t_H(\ulD')\le 2$. Indeed 
\[
2\;<\;t_H(\ulD')\;=\;\dim_{\BC_p}(D'_{\BC_p}\cap Fil^1_\mu D_{\BC_p})\;\le\; \dim_{\BC_p}Fil^1_\mu D_{\BC_p}\;=\;3
\]
would imply that $Fil^1_\mu D_{\BC_p}\subset D'_{\BC_p}$ and $\theta_{\bM(-1,2)}(EA)=0$. Thus $EA=T\cdot B$ for some morphism $B\in\Hom_\phi\bigl(\bM(-1,3)\,,\,\bM(1,2)\bigr)=(0)$ and this implies $E_1A_1+E_2A_2+E_3A_3=0$ in contradiction to the linear independence of $A_1,A_2,A_3$ over $\End_\phi\bigl(\bM(-1,2)\bigr)$. 

For the second kind of sub-isocrystal $t^{-1}\bD'\cong\bM(-1,2)$ and $t_N(\ulD')=1$. Again we claim that $t_H(\ulD')\le1$. Indeed 
\[
1\;<\;t_H(\ulD')\;=\;\dim_{\BC_p}(D'_{\BC_p}\cap Fil^1_\mu D_{\BC_p})\;\le\; \dim_{\BC_p}D'_{\BC_p}\;=\;2
\]
implies that $D'_{\BC_p}\subset Fil^1_\mu D_{\BC_p}$ and hence $t^{-1}\bD'\subset\bM_\mu\subset t^{-1}\bD=\bM(-1,2)^{\oplus3}$. Therefore $\bM_\mu\cong\bM(-1,2)\oplus\bigoplus_i\bM(s_i,n_i)$ for $-\frac{1}{2}\le\frac{s_i}{n_i}\le\frac{1}{2}$. We had defined $\mu$ by constructing a morphism $\bM(-1,3)\to\bM_\mu$. If the induced morphism $\bM(-1,3)\to\bM(s_i,n_i)$ were non-zero for some $i$, then $\frac{s_i}{n_i}\le-\frac{1}{3}$ and $s_i\le-1$ and $\,n_i\ge-2 s_i\ge2$. The constraints on the remaining $\frac{s_i}{n_i}$ imply that $\deg\bM_\mu\le-1$, a contradiction to Theorem~\ref{Thm1.13} and Lemma~\ref{LemmaKottwitzCond}. Therefore $A:\bM(-1,3)\to t^{-1}\bD$ factors through $t^{-1}\bD'$. But $t^{-1}\bD'$ lies in the kernel of an appropriate morphism $0\ne(E_1,E_2,E_3):t^{-1}\bD\to\bM(-1,2)$ for $E_i\in\End_\phi\bigl(\bM(-1,2)\bigr)$. Again this means $E_1A_1+E_2A_2+E_3A_3=0$ in contradiction to the linear independence of the $A_i$. This shows that $t_H(\ulD')\le t_N(\ulD')$ and so $\mu$ is weakly admissible.

\medskip
Therefore we have established in all cases that $\mu\in\breve\CF_b^{wa}\setminus\breve\CF_b^0$ and the theorem is proved.
\end{proof}

%
%

\section{Conjectures} \label{SectConjectures}
\setcounter{equation}{0}

As stated in Conjectures~\ref{Conj5.3b} and \ref{Conj5.3c} we expect that the open subspace $\breve\CF^0_b\subset\breve\CF^{wa}_b$ is arcwise connected and satisfies $\breve\CF^0_b(\LL)=\breve\CF^{wa}_b(\LL)$ if the value group of $\LL$ is finitely generated. This belief comes from the function field analog \cite{HartlPSp}, \cite[\S\S5.3 and 6.2]{HartlDict} of the theory developed here. In that analog we were able to prove the assertions corresponding to these conjectures; see \cite[Theorems 2.5.3 and 3.2.5]{HartlPSp}.

Unfortunately we have little more to support Conjecture~\ref{Conj5.3b}. However, the proof of Theorem~\ref{ThmC} and the computations in Example~\ref{Example6.4} and Theorem~\ref{Thm8.2} suggest that the fields $\LL$ with $\breve\CF^0_b(\LL)\subsetneq\breve\CF^{wa}_b(\LL)$ must contain elements of the form 
\[
\sum_{j=0}^{c-1}p^j\sum_{\nu\in\BZ} p^{c\nu}u_{j}^{p^{-e\nu}}
\]
for $0<c\le e$ and $u_j\in C$, $|u_j|<1$. Maybe one could show that fields containing such elements cannot have finitely generated value group.

\subsubsection*{Evidence for the connectedness of $\breve\CF^0_b$}

We suggest a strategy to prove that $\breve\CF^0_b$ is arcwise connected by noting that this follows from Conjecture~\ref{Conj5.3b}. In part \ref{Prop1.2b_1} of the following proposition we also add the proof of a fact mentioned on page~\pageref{EqDecent}.

\begin{proposition}\label{Prop1.2b}
\begin{enumerate}
\item \label{Prop1.2b_1}
There exists an arcwise connected paracompact open $\breve E$-analytic subspace $\breve\CF^{wa}_b$ of $\breve\CF^\an$ whose associated rigid analytic space is the period space $(\breve\CF^{wa}_b)^\rig$.
\item \label{Prop1.2b_2}
If $(G,b,\{\mu\})$ satisfies \eqref{Eq2.2} and Conjecture~\ref{Conj5.3b} is true then $\breve\CF^0_b$ is arcwise connected.
\item \label{Prop1.2b_3}
If $(G,b,\{\mu\})$ satisfies \eqref{Eq2.2} and $\breve\CF^0_b$ is connected then any open $\breve E$-analytic subspace $\breve\CF^a_b$ with $\breve\CF^0_b\subset\breve\CF^a_b\subset\breve\CF^{wa}_b$ is also arcwise connected.
\end{enumerate}
\end{proposition}

\begin{proof}[Proof of Proposition~\ref{Prop1.2b}]
\ref{Prop1.2b_1} Since $\breve\CF^\an$ has a finite covering by affine spaces, all its open $\breve E$-analytic subspaces are paracompact by Lemma~\ref{LemmaParacompact}. From its construction in \cite[Proposition 1.34]{RZ} the period space $(\breve\CF^{wa}_b)^\rig$ possesses an admissible covering $\{X_i^\rig\}_{i\in\BN}$ by admissible subsets $X_i^\rig\subset\breve\CF^\rig$ such that $\breve\CF^\rig\setminus X_i^\rig$ is a finite union $\bigcup_j\Spm B_{i,j}$ of affinoid subdomains and $X_i^\rig\subset X_{i+1}^\rig$. The $X_i^\rig$ correspond to open $\breve E$-analytic subspaces $X_i:=\breve\CF^\an\setminus\bigcup_j\CM(B_{i,j})$ of $\breve\CF^\an$. Moreover, it follows from \cite[Proposition 1.34]{RZ} that the closure $\ol X_i$ of $X_i$ in $\breve\CF^\an$ is a finite union of $\breve E$-affinoid subdomains and $\ol X_i\subset X_{i+1}$. Thus the union $\breve\CF^{wa}_b:=\bigcup_{i\in\BN}X_i$ is an open $\breve E$-analytic subspace of $\breve\CF^\an$ whose associated rigid analytic space equals $(\breve\CF^{wa})^\rig$. (Use \cite[Lemma 1.6.2]{Berkovich2} to see that the Grothendieck topology induced from $\breve\CF^{wa}_b$ on $\{\,\mu\in\breve\CF^{wa}_b:\;\CH(\mu)/\breve E\text{ is finite}\,\}$ coincides with the Grothendieck topology of $(\breve\CF^{wa}_b)^\rig$.)

To prove the connectedness of $\breve\CF^{wa}_b$ note that every $\breve E$-analytic space is locally arcwise connected by \cite[Theorem 3.2.1]{Berkovich1}. The flag variety $\breve\CF^\an$ has a finite covering by affine spaces, hence a countable affinoid covering by polydiscs $U\cong\CM\bigl(\breve E\langle\frac{x_1}{\zeta},\ldots,\frac{x_n}{\zeta}\rangle\bigr)$ for varying $\zeta\in\breve E$. Since the points $x\in\breve\CF^{wa}_b$ with $\CH(x)$ finite over $\breve E$ lie dense in $\breve\CF^{wa}_b$ by \cite[Proposition 2.1.15]{Berkovich1} it suffices to exhibit for every such point $x\in\breve\CF^{wa}_b\cap U$ a continuous map $\alpha$ from the compact interval $[0,|\zeta|\,]$ into $\breve\CF^{wa}_b\cap U$ such that $\alpha(0)=x$ and $\alpha(|\zeta|)$ is the point corresponding to the Gau{\ss} norm on $U$.

So let $x\in\breve\CF^{wa}_b\cap U$ with $\CH(x)$ finite over $\breve E$ be the point with coordinates $x_i=c_i\in\CH(x)$, $|c_i|\le|\zeta|$. To define $\alpha:[0,|\zeta|\,]\to\breve\CF^{wa}_b\cap U$ we expand each $f\in \breve E\langle\frac{x_1}{\zeta},\ldots,\frac{x_n}{\zeta}\rangle$ around $x$
\[
f\es=\es\sum_{\ul i\in\BN_0^{\,n}} a_{\ul i}\, (x_1-c_1)^{i_1}\cdots(x_n-c_n)^{i_n}\qquad\text{with}\quad a_{\ul i}\in \breve E(c_1,\ldots,c_n)\,.
\]
For every element $t\in[0,|\zeta|\,]$ we obtain an analytic point $P=P(x,t)$ in $U$ by setting
\[
|f|_{\SSC P}\es:=\es \sup\{\,|a_{\ul i}|\,t^{i_1+\ldots+i_n}:\es\ul i\in\BN_0^{\,n}\,\}\,.
\]
If $x$ is fixed it is easy to see that this defines a continuous and injective map $\alpha:[0,|\zeta|\,]\to U,\,t\mapsto P(x,t)$ with $\alpha(0)=x$ and $\alpha(|\zeta|)$ being the point corresponding to the Gau{\ss} norm on $U$. It remains to show that $\alpha(t)$ lies in $\breve\CF^{wa}_b$. For $t>0$ the prime ideal $\ker|\,.\,|_{\alpha(t)}:=\{\,f\in\breve E\langle\frac{x_1}{\zeta},\ldots,\frac{x_n}{\zeta}\rangle:|f|_{\alpha(t)}=0\,\}$ corresponding to $\alpha(t)$ is the zero ideal. Since by construction in \cite[Proposition 1.36]{RZ} the set $U\setminus\breve\CF^{wa}_b$ is a union of Zariski closed subsets of $U$ and since it does not contain the point $x$, we obtain $\alpha(t)\in\breve\CF^{wa}_b$. This implies that $\breve\CF^{wa}_b$ is arcwise connected.

\smallskip
\noindent
\ref{Prop1.2b_2} To prove the connectedness of $\breve\CF^0_b$ under the assumption that Conjecture~\ref{Conj5.3b} is true, we show that even $\alpha\bigl([0,|\zeta|\,]\bigr)\subset\breve\CF^0_b\cap U$. Namely, the value group of $\CH\bigl(\alpha(t)\bigr)$ is generated by $\bigl|\breve E(c_1,\ldots,c_n)\mal\bigr|$ and $t$ and hence is finitely generated. So Conjecture~\ref{Conj5.3b} will imply $\alpha(t)\in\breve\CF^0_b\cap U$ and so $\breve\CF^0_b$ is arcwise connected.

\smallskip
\noindent
\ref{Prop1.2b_3} To prove the connectedness of $\breve\CF^a_b$ we note that $\breve\CF^a_b(\LL)=\breve\CF^{wa}_b(\LL)$ for all finite extensions $\LL/\breve E$ because $\breve\CF^0_b(\LL)=\breve\CF^{wa}_b(\LL)$. Since the points $x\in\breve\CF^a_b$ with $\CH(x)$ finite over $\breve E$ lie dense in $\breve\CF^a_b$ by \cite[Proposition 2.1.15]{Berkovich1}, $\breve\CF^0_b$ is a dense connected subset of $\breve\CF^a_b$.
\end{proof}

\subsubsection*{Evidence for the maximality of $\breve\CF^0_b$}

We conjectured in Conjecture~\ref{Conj5.10} that the open $\breve E$-analytic subspace $\breve\CF^0_b\subset\breve\CF^{wa}_b$ and the local system $\CV$ from Theorem~\ref{ThmDOR} are maximal. From this also Conjectures~\ref{Conj6.5} and \ref{Conj7.6} would follow. For this maximality it remains to prove that the local system does not extend to a larger open subspace $\breve\CF^0_b\subsetneq Y\subset\breve\CF^{wa}_b$. This is suggested by the following results of Andreatta and Brinon~\cite{Andreatta,AB,AB2}. Assume there exists a local system $\CV$ of $\BQ_p$-vector spaces on $Y$ as in Theorem~\ref{ThmDOR} and let $\wt Y$ be the space of $\BZ_p$-lattices inside $\CV$. 

Consider a morphism $f:U\to\wt Y$ of $\breve E$-analytic spaces such that $U$ possesses an affine formal model $\Spf R$ that is \'etale over $\Spf W(\BF_p^{\,\alg})\langle T_1,T_1^{-1},\ldots,T_d,T_d^{-1}\rangle$. (This hypothesis on $U$ can even be weakened slightly; see the conditions in~\cite{Andreatta,AB}.) The pullback of $\CV$ to $U$ corresponds by Proposition~\ref{Prop1.2c} to a representation $\rho_V:\pi_1^\et(U,\bar x)\to\GL(V)(\BQ_p)$ which stabilizes a $\BZ_p$-lattice in $V=\CV_{\bar x}$. Thus $\rho_V$ factors through
\[
\pi_1^\et(U,\bar x)\;\to\;\pi_1^\alg(U,\bar x)\;\to\;\GL(V)(\BZ_p)\;\subset\;\GL(V)(\BQ_p), 
\]
where $\pi_1^\alg(U,\bar x)$ is the algebraic fundamental group classifying \emph{finite} \'etale coverings; see \cite{dJ}. By \cite[Theorem 7.11]{Andreatta} $\rho_V$ corresponds to an \'etale $(\phi,\Gamma)$-module $\bM$ over $\bA_R$ which is overconvergent by \cite{AB}. If $\mu\in U$ is an analytic point and $C$ is the completion of an algebraic closure of $\CH(\mu)$ the fiber $\bM_\mu\otimes\wt\bB(C)$ at $\mu$ is therefore of the form $\bM_\mu^\dagger\otimes_{\wt\bB^\dagger(C)}\wt\bB(C)$ for a $\phi$-module $\bM^\dagger_\mu$ over $\wt\bB^\dagger(C)$. By \cite[Proposition 5.5.1]{Kedlaya} the $\phi$-module $\bM_\mu^\dagger\otimes_{\wt\bB^\dagger(C)}\wt\bB^\dagger_\rig(C)$ over $\wt\bB^\dagger_\rig(C)$ is isomorphic to $\bM(0,1)\otimes_{\BQ_p}V$. Due to (\ref{EqRZ}) we have $\bM_\mu^\dagger\otimes_{\wt\bB^\dagger(C)}\wt\bB^\dagger_\rig(C)\cong\bM\bigl((V_{K_0},b\!\cdot\!\phi,Fil^\bullet_{f(\mu)} V_{\CH(\mu)})\bigr)$. Thus the image of $U$ in $Y$ must be contained in $\breve\CF_b^0$. 

However, at present this approach does not prove Conjecture~\ref{Conj5.10} because unfortunately $Y$ cannot be covered by $\breve E$-analytic spaces $U$ allowed by \cite{Andreatta,AB}. For instance consider on an annulus $Y_1=\CM\bigl(\breve E\langle T,\frac{p}{T}\rangle\bigr)$ the point $y_1$ with $|\sum_ia_iT^i|_y=\sup\{|a_i|\,|p|^{ri}:i\in\BZ\}$ for $r\in[0,1]\setminus\BQ$. Furthermore, consider a point $y_2$ with $\CH(y_2)=\BC_p$ in a small polydisc $Y_2=\CM\bigl(\breve E\langle\frac{T_1}{\zeta},\ldots,\frac{T_n}{\zeta}\rangle\bigr)$ for $\zeta\in\breve E$, $0<|\zeta|\ll 1$ (such points exist; see \cite[Example A.2.2 (d)]{HartlPSp}). If $Y$ contains $Y_1\times Y_2$ then the point $y=(y_1,y_2)\in Y$ does not lie in the image of any $U$ as above. So this method does not yield $y\in\breve\CF^0_b$. Also Conjecture~\ref{Conj5.3b} does not give $y\in\breve\CF^0_b$ because the value group $|\CH(y)\mal|\supset|\CH(y_2)\mal|\cong\BQ$ is not finitely generated. 

\begin{remark}\label{RemarkKedlaya}
Due to these restrictions of Andreatta's and Brinon's theory \cite{Andreatta,AB}, Kedlaya uses a general new foundation of relative $p$-adic Hodge theory \cite{KedlayaLiu12,KedlayaLiu13} to obtain the generalizations of our results, which he announced in \cite{Kedlaya10}.
\end{remark}

%
%

\begin{appendix}
\section{Berkovich's Analytic Spaces} \label{AppBerkovichSpaces}
\setcounter{equation}{0}

Let $\CO_\LF$ be a (not necessarily discrete) complete valuation ring of rank one and let $\LF$ be its fraction field. We briefly recall Berkovich's~\cite{Berkovich1,Berkovich2} theory of \emph{$\LF$-analytic spaces}. Let $B$ be an affinoid $\LF$-algebra in the sense of \cite[Chapter 6]{BGR} with $\LF$-Banach norm $|\,.\,|$. Berkovich calls these algebras \emph{strictly $\LF$-affinoid}.

\begin{definition} \label{DefAnalyticPoint} 
An {\em analytic point\/} $x$ of $B$ is a semi-norm $|\,.\,|_x:B \to \BR_{\geq 0}$ which satisfies: 
\begin{enumerate} 
\item 
$|f+g|_x \leq \max\{\,|f|_x,|g|_x\,\}$ \quad for all $f,g \in B$, 
\item  
$|fg|_x = |f|_x \,|g|_x$ \quad for all $f,g \in B$, 
\item 
$|\lambda|_x = |\lambda|$ \quad for all $\lambda \in \LF$, 
\item 
$|\,.\,|_x:B \to \BR_{\geq 0}$ is continuous with respect to the norm $|\,.\,|$ on $B$. 
\end{enumerate} 
The set of all analytic points of $B$ is denoted~$\CM(B)$. 
On $\CM(B)$ one considers the coarsest topology such that for every $f\in B$ the map $\CM(B) \to \BR_{\geq0}$ given by $x \mapsto |f|_x$ is continuous. Equipped with this topology, $\CM(B)$ is a compact Hausdorff space; see~
\cite[Theorem 1.2.1]{Berkovich1}. Such a space is called a \emph{strictly $\LF$-affinoid space}.

Every morphism $\alpha:B\to B'$ of affinoid $\LF$-algebras is automatically continuous and hence induces a continuous morphism $\CM(\alpha):\CM(B') \to \CM(B)$ by mapping the semi-norm $B' \to \BR_{\geq 0}$ to the composition $B\to B'\to \BR_{\geq 0}$. By definition the $\CM(\alpha)$ are the \emph{morphisms} in the category of strictly $\LF$-affinoid spaces. In particular, for an affinoid subdomain $\Spm B'\subset \Spm B$ this morphism identifies $\CM(B')$ with a closed subset of $\CM(B)$. 
 \end{definition} 
 
For every analytic point $x\in\CM(B)$ we let $\ker|\,.\,|_x\,:=\,\{\,b\in B: |b|_x=0\,\}$. It is a prime ideal in $B$. We define the {\em (complete) residue field of $x$} as the completion with respect to $|\,.\,|_x$ of the fraction field of $B/\ker|\,.\,|_x$. It will be denoted $\CH(x)$. There is a natural continuous homomorphism $B \to \CH(x)$ of $\LF$-algebras. Conversely let $\LL$ be a \emph{complete extension} of $\LF$, by which we mean a field extension of $\LF$ equipped with an absolute value $|\,.\,|:\LL\to\BR_{\geq0}$ which restricts on $\LF$ to the norm of $\LF$ such that $\LL$ is complete with respect to $|\,.\,|$. Any continuous $\LF$-algebra homomorphism $B\to \LL$ defines on $B$ a semi-norm which is an analytic point. 

\medskip

In \cite[\S1.2]{Berkovich2} Berkovich defines the category of \emph{strictly $\LF$-analytic spaces}. These spaces are topological spaces which admit an atlas with strictly $\LF$-affinoid charts. 

The spaces $\CM(B)$ for affinoid $\LF$-algebras $B$ are examples for strictly $\LF$-analytic spaces. 
Other examples arise from schemes $Y$ which are locally of finite type over $\LF$, see \cite[\S3.4]{Berkovich1}. Namely if $Y=\BA_\LF^n$ is affine $n$-space over $\LF$ the associated strictly $\LF$-analytic space $Y^\an=(\BA_\LF^n)^\an$ consists of all semi-norms on the polynomial ring $\LF[y_1,\ldots,y_n]$ as in Definition~\ref{DefAnalyticPoint} (a)--(c). The topology on $(\BA^n_\LF)^\an$ is the coarsest topology such that for all $f\in \LF[y_1,\ldots,y_n]$ the map $(\BA^n_\LF)^\an\to\BR_{\ge0}, x\mapsto |f|_x$ is continuous. The space $(\BA_\LF^n)^\an$ is the union of the increasing sequence of compact polydiscs $\CM\bigl(\LF\langle \frac{y_1}{\zeta^m},\ldots,\frac{y_n}{\zeta^m}\rangle\bigr)$ of radii $(|\zeta^{m}|,\ldots,|\zeta^{m}|)$ for $m\in\BN_0$ where $\zeta\in \LF$ has $|\zeta|>1$. If $Y\subset\BA_\LF^n$ is a closed subscheme of affine $n$-space with coherent ideal sheaf $\CJ$, the ideal sheaf $\CJ\CO_{(\BA_\LF^n)^\an}$ defines a closed strictly $\LF$-analytic subspace $Y^\an$ of $(\BA_\LF^n)^\an$. Finally, if $Y$ is arbitrary and $\{Y_i\}_i$ is a covering of $Y$ by affine open subschemes then one can glue the associated strictly $\LF$-affinoid spaces $Y_i^\an$ to the $\LF$-analytic space $Y^\an$. Moreover $Y^\an$ is Hausdorff if and only if the scheme $Y$ is separated, see \cite[Theorems 3.4.1 and 3.4.8]{Berkovich1}.

\medskip

The relation between strictly $\LF$-analytic spaces, rigid analytic spaces, and formal schemes is as follows.
To every strictly $\LF$-analytic space $X$ which is Hausdorff one can associate a quasi-separated rigid analytic space
\[
X^\rig\es:=\es\{\,x\in X:\es\CH(x) \text{ is a finite extension of }\LF\,\},
\]
see \cite[\S1.6]{Berkovich2}. Recall that a rigid analytic space is called \emph{quasi-separated} if the intersection of any two affinoid subdomains is a finite union of affinoid subdomains. To describe the subcategories on which the functor $X\mapsto X^\rig$ is an equivalence we need the following terminology. A topological Hausdorff space is called \emph{paracompact} if every open covering $\{U_i\}_i$ has a locally finite refinement $\{V_j\}_j$, where \emph{locally finite} means that every point has a neighborhood which meets only finitely many of the $V_j$. On the other hand an admissible covering of a rigid analytic space is said to be \emph{of finite type} if every member of the covering meets only finitely many of the other members. A rigid analytic space over $\LF$ is called \emph{quasi-paracompact} if it possesses an admissible affinoid covering of finite type. Similarly we define the notions \emph{of finite type} and \emph{quasi-paracompact} also for (an open covering of) an admissible formal $\CO_\LF$-scheme in the sense of Raynaud~\cite{Raynaud}; see also \cite{FRG1,FRG2,Bosch}.

\begin{theorem}\label{ThmFormalRigBerkovich}
The following three categories are equivalent:
\begin{enumerate}
\item 
the category of paracompact strictly $\LF$-analytic spaces,
\item 
the category of quasi-separated quasi-paracompact rigid analytic spaces over $\LF$, and
\item 
the category of quasi-paracompact admissible formal $\CO_\LF$-schemes, localized by admissible formal blowing-ups.
\end{enumerate}
\end{theorem}

\begin{proof}
It is shown in \cite[Theorem 1.6.1]{Berkovich2} that $X\mapsto X^\rig$ is an equivalence between (a) and (b). The equivalence of (c) with (b) is due to Raynaud. See \cite[Theorem 2.8/3]{Bosch} for a proof.
\end{proof}

Regarding paracompactness the following result was proved in \cite[Lemma A.2.6]{HartlPSp}.

\begin{lemma} \label{LemmaParacompact}
Let $\LF$ be discretely valued with countable residue field and let $X$ be a strictly $\LF$-analytic space. Assume that $X$ is Hausdorff and admits a countable covering by strictly $\LF$-affinoid spaces. Then every open subset of $X$ is a paracompact strictly $\LF$-analytic space.
\end{lemma}
This applies in particular if $X=Y^\an$ for a separated scheme $Y$ of finite type over $\LF$.

\medskip

We will also need the following well known fact which we were unable to find in the literature. 

\begin{proposition}\label{PropA.3}
Let $f:X\to Y$ be an \'etale morphism of $\LF$-analytic spaces (\cite[\S3.3]{Berkovich2}) such that $f$ is bijective on $C$-valued points for any algebraically closed complete extension $C$ of $\LF$. Then $f$ is an isomorphism.
\end{proposition}

\begin{proof}
Let $x\in X$ and $y=f(x)\in Y$. Since $f$ is quasi-finite there are open neighborhoods $V=V_x\subset X$ of $x$ and $U=U_y\subset Y$ of $y$ with $f|_{V}:V\to U$ finite \'etale. Since $f$ is open by \cite[Proposition 3.2.7]{Berkovich2} we may replace $U$ by $f(V)$. For any affinoid $\CM(A)\subset U$ the preimage $(f|_{V})^{-1}(\CM(A))=\CM(B)$ is affinoid with $B$ finite flat over $A$ by \cite[Proposition 3.2.3]{Berkovich2}. Hence $B$ is finite locally free since $A$ is noetherian by \cite[Proposition 2.1.3]{Berkovich1}. For any $C$-valued point of $\CM(A)$ the fiber $\CM(B\otimes_AC)\to\CM(C)$ is finite \'etale by \cite[Corollary 3.3.8]{Berkovich2} and bijective by our assumption. Since $C$ is algebraically closed, $B\otimes_AC=C$ and so $\rk_A B=1$. Therefore the inverse map $(f|_{V})^{-1}:U\to V$ exists by \cite[Proposition 1.2.15(i)]{Berkovich2}. Now our assumption on bijectivity on $C$-valued points implies that the local inverses $(f|_{V})^{-1}:U\isoto V$ glue by \cite[Proposition 1.3.2]{Berkovich2} to the global inverse of $f$.
\end{proof}

\end{appendix}

%
%

\vfill

\noindent
Urs Hartl\\
University of Muenster\\
Institute of Mathematics\\
Einsteinstr.~62\\
D -- 48149 Muenster\\
Germany\\[1mm]
http:/\!/www.math.uni-muenster.de/u/urs.hartl/

\end{document}